\documentclass[preprint,12pt]{elsarticle}

\usepackage{amssymb}
\usepackage{amsmath}
\usepackage{mathrsfs}
\usepackage{subfig}
\usepackage{graphicx} 
\usepackage{epstopdf}
\usepackage{latexsym}
\usepackage{color}
\usepackage{lineno} 

\makeatletter
\newcommand{\rmnum}[1]{\romannumeral #1}
\newcommand{\Rmnum}[1]{\expandafter\@slowromancap\romannumeral #1@}
\makeatother

\newcommand{\reffig}[1]{Figure \ref{#1}}

\newtheorem{thm}{Theorem}
\newtheorem{lem}[thm]{Lemma}
\newtheorem{definition}[thm]{Definition}

\newdefinition{rmk}{Remark}
\newdefinition{example}{Example}

\newproof{pf}{\emph{Proof}}
\newproof{pot}{Proof of Theorem \ref{thm2}}

\usepackage[ruled,vlined,linesnumbered]{algorithm2e}

\begin{document}
\begin{frontmatter}



\title{Modified PHT-splines}


\author[a]{Qian Ni}
\author[b]{Xuhui Wang}
\author[a]{Jiansong Deng\corref{cor1}}
\cortext[cor1]{Corresponding author.}
\ead{dengjs@ustc.edu.cn}
\address[a]{School of Mathematical Sciences, University of Science and Technology of China, Hefei, Anhui 230026, P. R. China}
\address[b]{School of Mathematics, Hefei University of Technology, Hefei, Anhui 230009, P. R. China}
\begin{abstract}
The local refinement of PHT-splines (polynomial splines over hierarchical T-meshes) is achieved by a simple cross insertion, which may introduce superfluous control points or coefficients.
By allowing split-in-half in mesh refinement, modified hierarchical T-meshes are defined. Using this approach, polynomial splines defined over the modified hierarchical T-meshes (modified PHT-splines) are introduced to increase the flexibility of PHT-splines.
Numerical examples demonstrate the advantages of our new splines when applied to surface fitting and isogeometric analysis problems with anisotropic features.
\end{abstract}


\end{frontmatter}


\section{Introduction}

The multivariate splines used in geometric design (GM) and isogeometric analysis (IGA) are often based on the tensor-product of B-splines or NURBS. However, the standard tensor-product mechanism prevents a localized editing of the mesh. As such, a large number of superfluous control points or coefficients may be contained under the tensor-product representations when we deal with problems having anisotropic features. To overcome this weakness and provide more flexibility, locally refinable splines are introduced,  which generalize the B-spline model in different ways. Hence the development of locally refinable splines has recently become an active research topic in geometric modeling and isogeometric analysis.
There are a number of different locally refinable splines, and among them, a) T-splines,  b) HB-splines and THB-splines, c) LR-splines, d) PHT-splines are very popular. We shall review these splines briefly below.

1) By allowing T-junctions in the control meshes, \emph{T-splines} \cite{tom2003,tom2004} were introduced to remove the redundant control points in geometric modeling. The potential of T-splines in isogeometric analysis was reported in \cite{bazilevs-2010-iga,dorfel-2009-iga}. Additionally, AST-splines (analysis-suitable T-splines) were introduced to fix the linear dependency problem of T-splines \cite{li2012,Scott-2012-iga}, and are important for their applications in isogeometric analysis. Recently, AS++ T-splines (analysis-suitable++ T-splines), which include analysis-suitable (AS) T-splines as a special case, were also introduced  \cite{li2018}.

2) As a classical extension of tensor-product representations, \emph{hierarchical B-splines} (HB-splines) allow an effective local control of refinement \cite{forsey-1998-HB}. The basis functions of HB-splines are locally supported, linearly independent and non-negative.
The application of HB-splines in geometric modeling and isogeometric analysis was reported in \cite{evans-2015-hb-gm,greiner-1997-HB-gm,schillinger-2012-hb-gm,vuong-2011-HB-iga}. Unfortunately, the standard hierarchical construction does not preserve the partition of unity property.
By adequately truncating basis functions in the construction of HB-splines, \emph{truncated hierarchical B-splines} (THB-splines), which form a convex partition of unity, were introduced by Giannelli et al. in \cite{giannelli-2012-thb,giannelli-2014-thb}. In addition, based on their good stability \cite{giannelli-2014-thb} and approximation abilities \cite{speleers-2017},  THB-splines are well suited for applications in geometric modeling, computer aided design, and isogeometric analysis \cite{giannelli-2016-thb,giannelli-2014-thb, kiss-2014-thb-cad,speleers-2016-thb-gm}. In addition, based on the techniques of truncation, the normalized quasi-hierarchical basis was proposed for hierarchical Powell-Sabin splines in \cite{speleers-2009-Powell-sabin}. Recently, truncated T-splines are also introduced in \cite{Weixd2017-iga-truncated}.

3) Based on the idea of splitting basis functions, \emph{locally refined} (LR) splines were proposed by Dokken et al. in \cite{dokken-2013-lr}.  LR splines possess the nested spaces property while they have difficulties with linear independence. The properties of LR splines and their hierarchical construction in 2D were further discussed in \cite{bressan-2013-lr,bressan-2015-lr-2d}. LR splines were also applied in isogeometric analysis \cite{Johannessen-2014-lr-iga}.
Recently,  Johannessen et al. analyzed the
corresponding stiffness and mass matrices for HB-splines, THB-splines and LR-splines, in terms of sparsity patterns and conditioning numbers \cite{johannessen-2015-compare}.
The comparison shows that the basis functions of these splines
in general do not span the same space, and that conditioning numbers are comparable.


4) \emph{Polynomial splines over hierarchical T-meshes} (PHT-splines), i.e., bi-cubic splines over hierarchical T-meshes with reduced $C^1$ regularity, use a modification mechanism to obtain nested spline spaces \cite{deng2008}. PHT-splines modify their basis functions in coarse levels directly by resetting the B{\'e}zier ordinates corresponding to the new basis vertices. The modification mechanism ensures the properties like partition of unity, nonnegativity, and linear independence. Therefore, PHT-splines can be applied successfully in surface model reconstruction  \cite{deng2008,li2007,wang2010}, adaptive finite element method for elliptic equations \cite{tian2011}, and isogeometric analysis for solving elastic problems \cite{n2011a,n2011b,pingwang2011}.
The local refinement for hierarchical T-meshes of PHT-splines is achieved by cross insertion \cite{deng2008}, i.e., splitting each candidate cell into four subcells by inserting a cross. This simple refinement rule may still introduce redundant control points or coefficients when applied to the models or problems with anisotropic features, i.e., the problems or models are directionally dependent.

The present paper is devoted to increasing flexibility to PHT-splines, which can handle models or problems with anisotropic features efficiently. Instead of performing refinement exclusively by cross insertion, modified hierarchical T-meshes are obtained by allowing for the splitting of cells into  halves (horizontal or vertical). Similar to the construction of classical PHT-splines, polynomial splines over modified hierarchical T-meshes (\emph{modified PHT-splines}) are also presented. Numerical experiments show that our new splines have advantages when applied to problems with anisotropic features.  It is worth noting that using anisotropic refinement for surface construction has been reported recently in \cite{engleitner-2017-aniso,engleitner-2017-aniso2}.

%

The remainder of the paper is organized as follows.
 Section \ref{section-pht} reviews the definition of T-meshes and hierarchical T-meshes. The definition of modified hierarchical T-meshes and the dimension formula for polynomial spline spaces over modified hierarchical T-meshes are also provided. In Section \ref{section-basis-constuct}, we discuss the construction of basis functions for modified PHT-splines. Section \ref{section-fit-open-mesh} presents three examples for fitting open meshes with modified PHT-splines. In Section \ref{section-iga}, the adaptive isogeometric analysis based on modified PHT-splines is provided. Section \ref{section-conclusion} concludes the paper with a summary and some future work.

\section{Polynomial splines over T-meshes}\label{section-pht}
In this section, we briefly review some notations of T-meshes and hierarchical T-meshes. The definition of modified hierarchical T-meshes is provided and the dimension formula for polynomial splines space over modified T-meshes is also given.

\subsection{T-meshes}
\begin{figure}[htbp!]
\centering
\subfloat[A T-mesh]{
\includegraphics[width=1.2in]{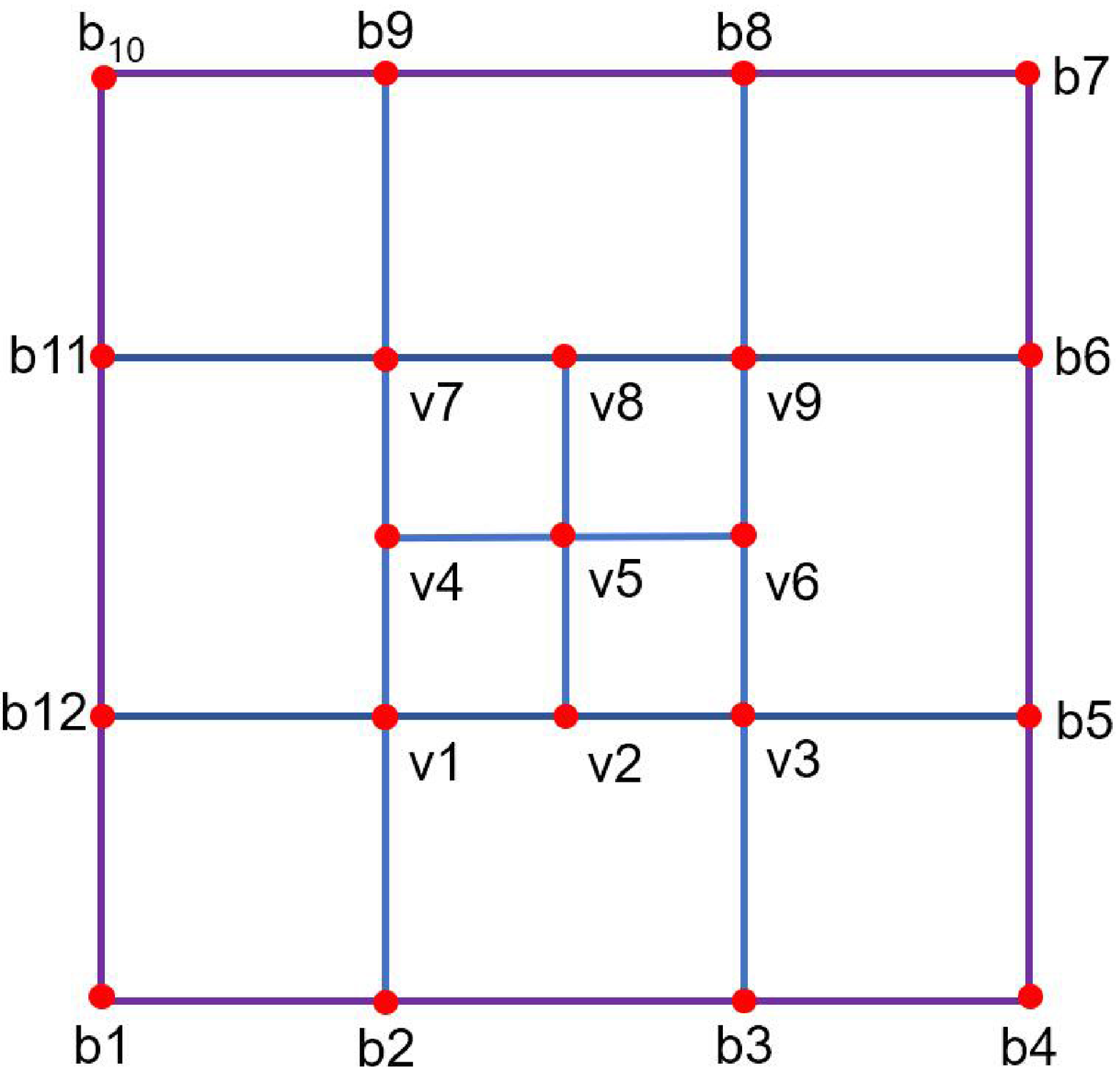}}\hspace{2cm}
\subfloat[A non-T-mesh]{
\includegraphics[width=1.2in]{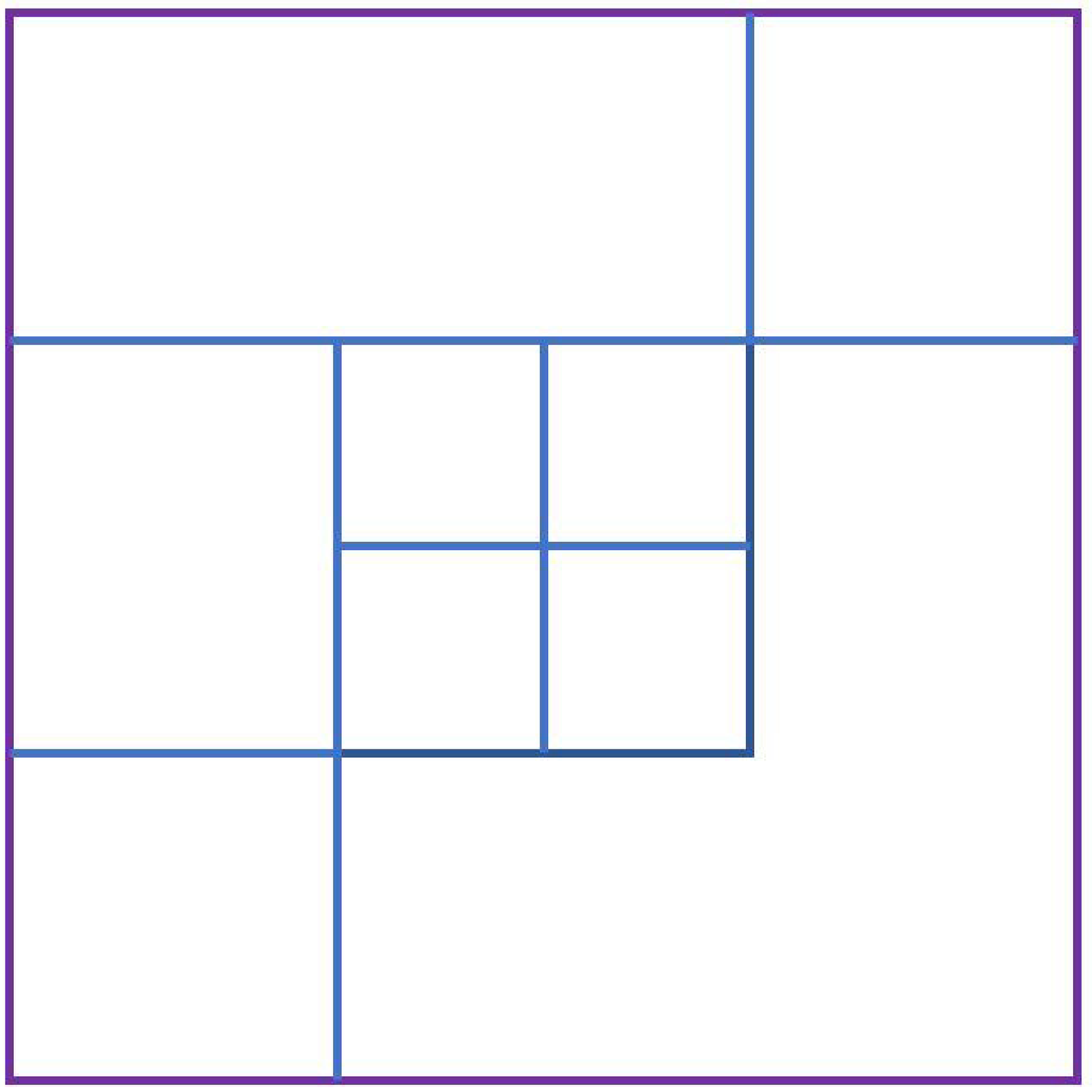}}
\caption{A T-mesh and a non-T-mesh.}
\label{fig:tmesh}
\end{figure}

A \emph{T-mesh} is a rectangular grid that allows T-junctions \cite{deng2006,tom2003}.
It is assumed that except for the corner case, the endpoints of each grid line in the T-mesh must be on two other grid lines,
and each cell or facet in the grid must be a rectangle. \reffig{fig:tmesh}(a) shows an example of a T-mesh,
and \reffig{fig:tmesh}(b) shows an example of a non-T-mesh. A grid point in a T-mesh is also called a \emph{vertex} of the T-mesh.
If a vertex is on the boundary of the domain, then it is called a \emph{boundary vertex}. Otherwise, it is called an \emph{interior vertex}.
For example, $\mathbf{b}_{i}$, $i = 1, \ldots, 12$ in \reffig{fig:tmesh}(a) are boundary vertices,
while $\mathbf{v}_{i}$, $i=1, \ldots, 9$ are interior vertices.
Interior vertices have two types: one is \emph{crossing vertex},
e.g., $\mathbf{v}_{1}$, $\mathbf{v}_{3}$, $\mathbf{v}_{5}$, $\mathbf{v}_{7}$, and $\mathbf{v}_{9}$ in \reffig{fig:tmesh}(a),
and the other is \emph{T-vertex}, e.g., $\mathbf{v}_{2}$, $\mathbf{v}_{4}$, $\mathbf{v}_{6}$, and $\mathbf{v}_{8}$ in \reffig{fig:tmesh}(a).
The line segment connecting two adjacent vertices on a grid line is called an \emph{edge} of the T-mesh.
If an edge is on the boundary of the T-mesh, it is called a \emph{boundary edge}.
Otherwise, it is called an \emph{interior edge}.
A cell is called an \emph{interior cell} if all four edges are interior. Otherwise, it is called a \emph{boundary cell}.

In terms of neighbor relations of cells in a T-mesh, a pair of  cells are called \emph{adjacent} if they share one common edge.
Moreover, a pair of adjacent cells are called \emph{horizontally (vertically) aligned-adjacent} if they share only one vertical (horizontal) edge,
and the knot intervals of these two cells are identical in the vertical (horizontal) direction.
A pair of cells are \emph{aligned-adjacent} if they are horizontally or vertically aligned-adjacent.

\subsection{Hierarchical T-meshes and modified hierarchical T-meshes}

As a special type of T-mesh, \emph{hierarchical T-mesh} $\mathbb{T}$ is defined in a hierarchical manner \cite{deng2008}.
To describe the level of hierarchical T-mesh, we fix the following notations in the rest of the paper.
$\mathbb{T}_k$ denotes the T-mesh at level $k$. Let $\mathcal{F}_k$ be the set consisting of all cells of level $k$.
Note that the level of cells and the level of mesh are defined separately; not every cell in $\mathbb{T}_k$ is of level $k$.
In this paper, we focus on a special type of hierarchical T-mesh, i.e., the cells to be subdivided are chosen from $\mathcal{F}_k$.
Let $\Theta_k \subset \mathcal{F}_k$ be the set consisting of all cells in $\mathcal{F}_k$ that are to be subdivided at level $k$.
Then, a hierarchical T-mesh $\mathbb{T}$ is defined as follows:
\begin{enumerate}[(i)]
  \item Start with a tensor product mesh $\mathbb{T}_0$. Let $\mathcal{F}_0$  be the set of all cells in $\mathbb{T}_0$.
  \item Recursive steps: Subdivide each cell in $\Theta_k$ by inserting a cross to get a new T-mesh $\mathbb{T}_{k+1}, k=0,1,\ldots, N-1$.
  Let $\mathcal{F}_{k+1}$ be the set consisting of all new cells, which emerge at level $k+1$ after subdivision.
  \item $\mathbb{T}=\mathbb{T}_N$.
\end{enumerate}
\begin{rmk}\label{remark-celllevel}
Based on the above (\rmnum{2}), if a cell $\theta \notin \Theta_k$, i.e., $\theta$ is not subdivided at level $k$, the cell $\theta$ will be excluded from subdivision henceforth.
\end{rmk}
The level of a vertex is defined as follows.
Given a vertex $\mathbf{v}$ in $\mathbb{T}$,
if $\mathbf{v}\in \mathbb{T}_k$ but $\mathbf{v}\not\in \mathbb{T}_{l}, l=0,\ldots,k-1$,
then it is called a vertex of level $k$.
Specifically, every vertex in $\mathbb{T}_0$ is a vertex of level 0.

\begin{figure}[htbp!]
\centering
\subfloat[$\mathbb{T}_0$]{
\includegraphics[width=1.2in]{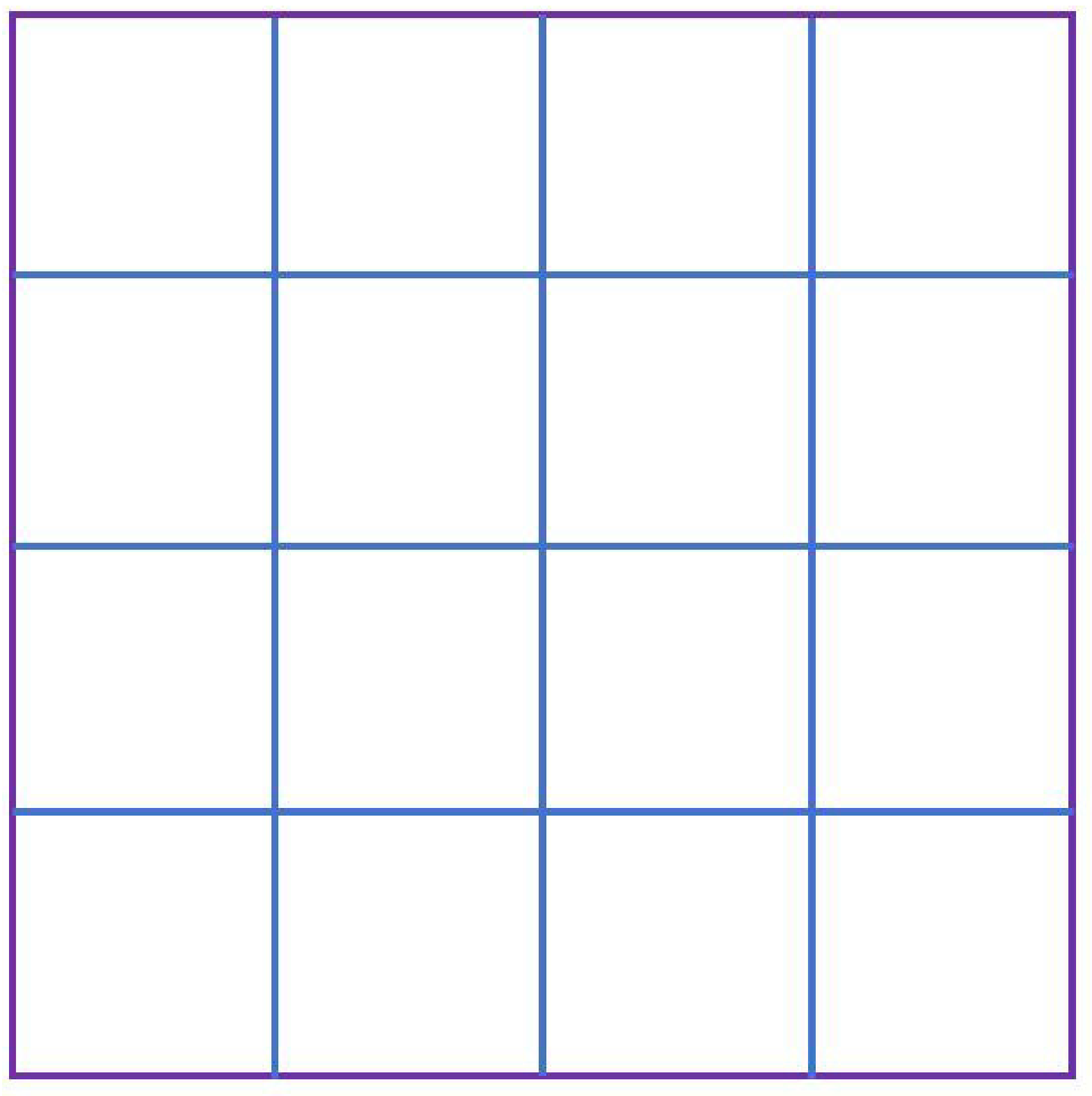}}\hspace{1cm}
\subfloat[$\mathbb{T}_1$]{
\includegraphics[width=1.2in]{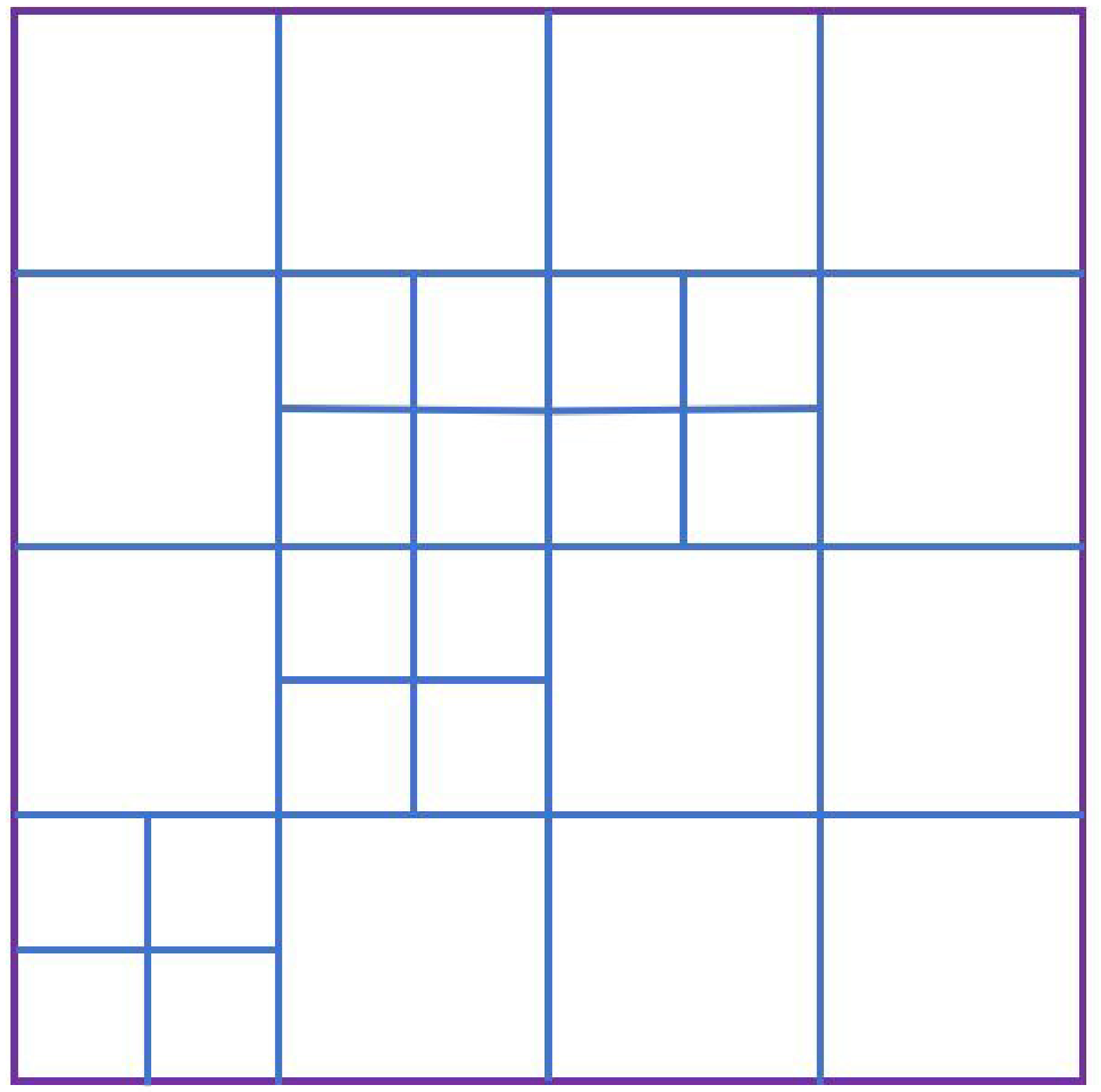}}\hspace{1cm}
\subfloat[$\mathbb{T}_2$]{
\includegraphics[width=1.2in]{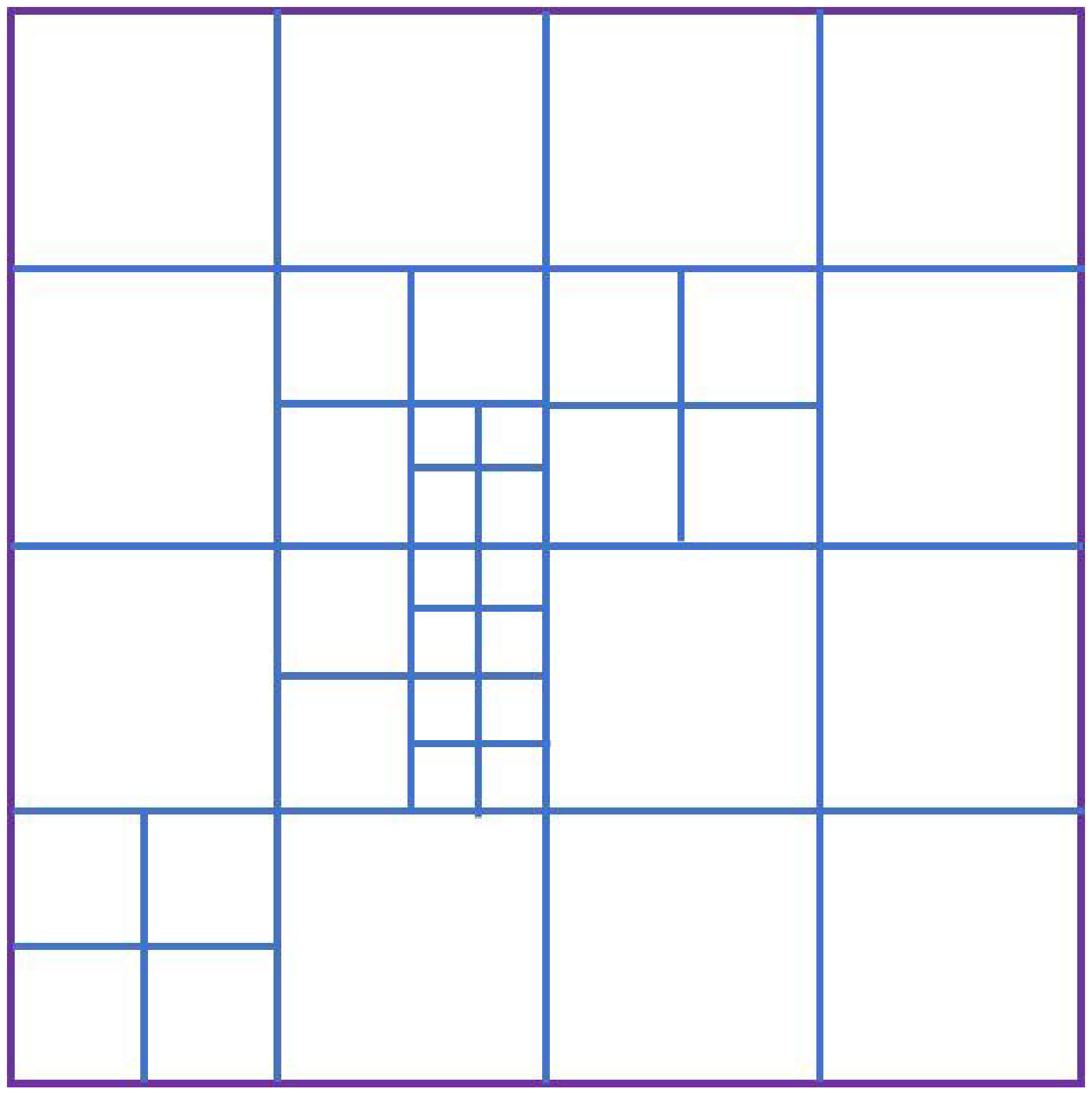}}
\caption{An example of hierarchical T-mesh.}
\label{fig:ht}
\end{figure}

In \reffig{fig:ht}, an example is provided to illustrate the dynamic refinement process of a hierarchical T-mesh.

\begin{figure}[htbp!]
\centering
\subfloat[$\mathbb{T}_0$]{
\includegraphics[width=1.2in]{ht0.eps}}\hspace{1cm}
\subfloat[$\mathbb{T}_1$]{
\includegraphics[width=1.2in]{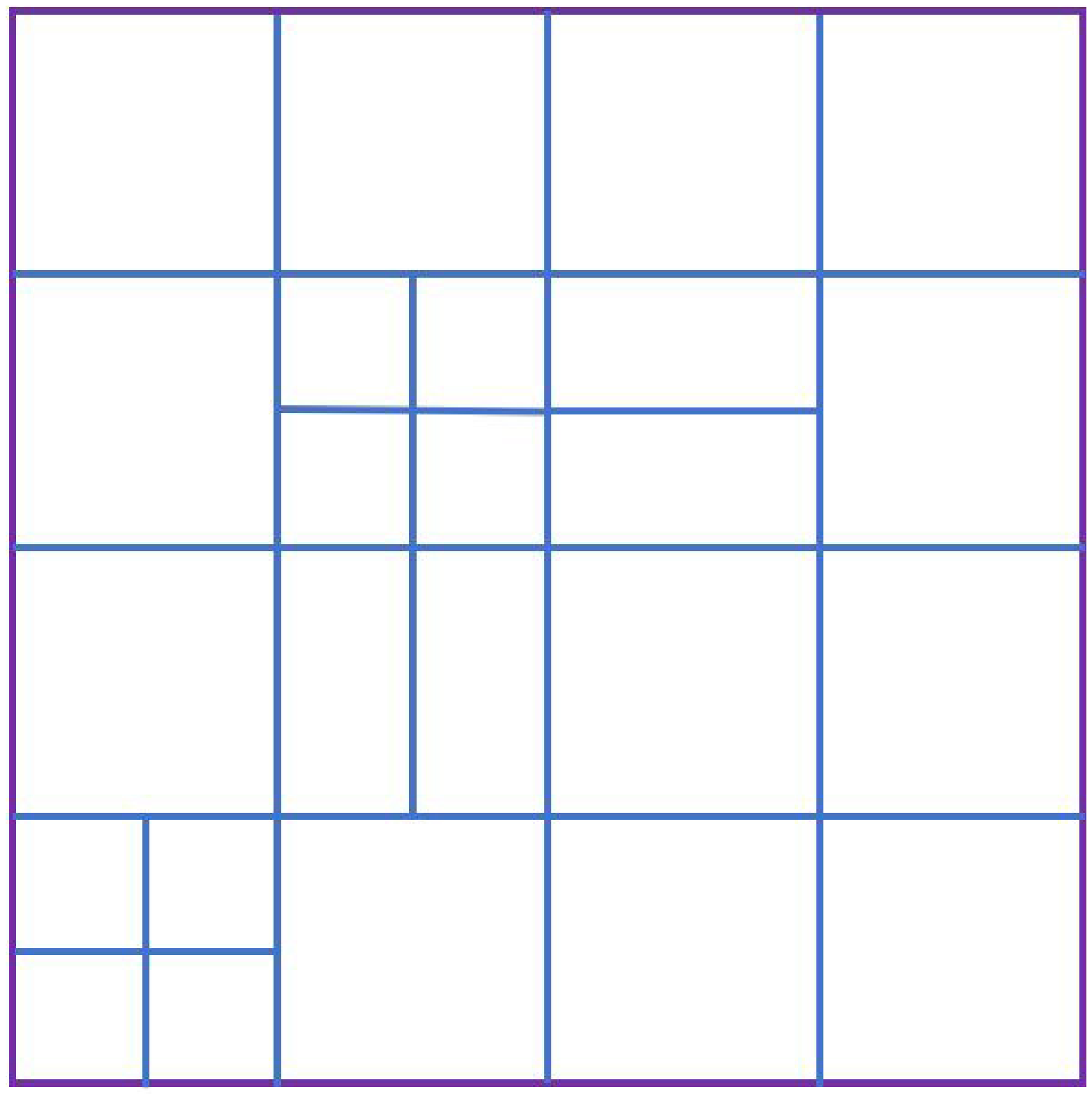}}\hspace{1cm}
\subfloat[$\mathbb{T}_2$]{
\includegraphics[width=1.2in]{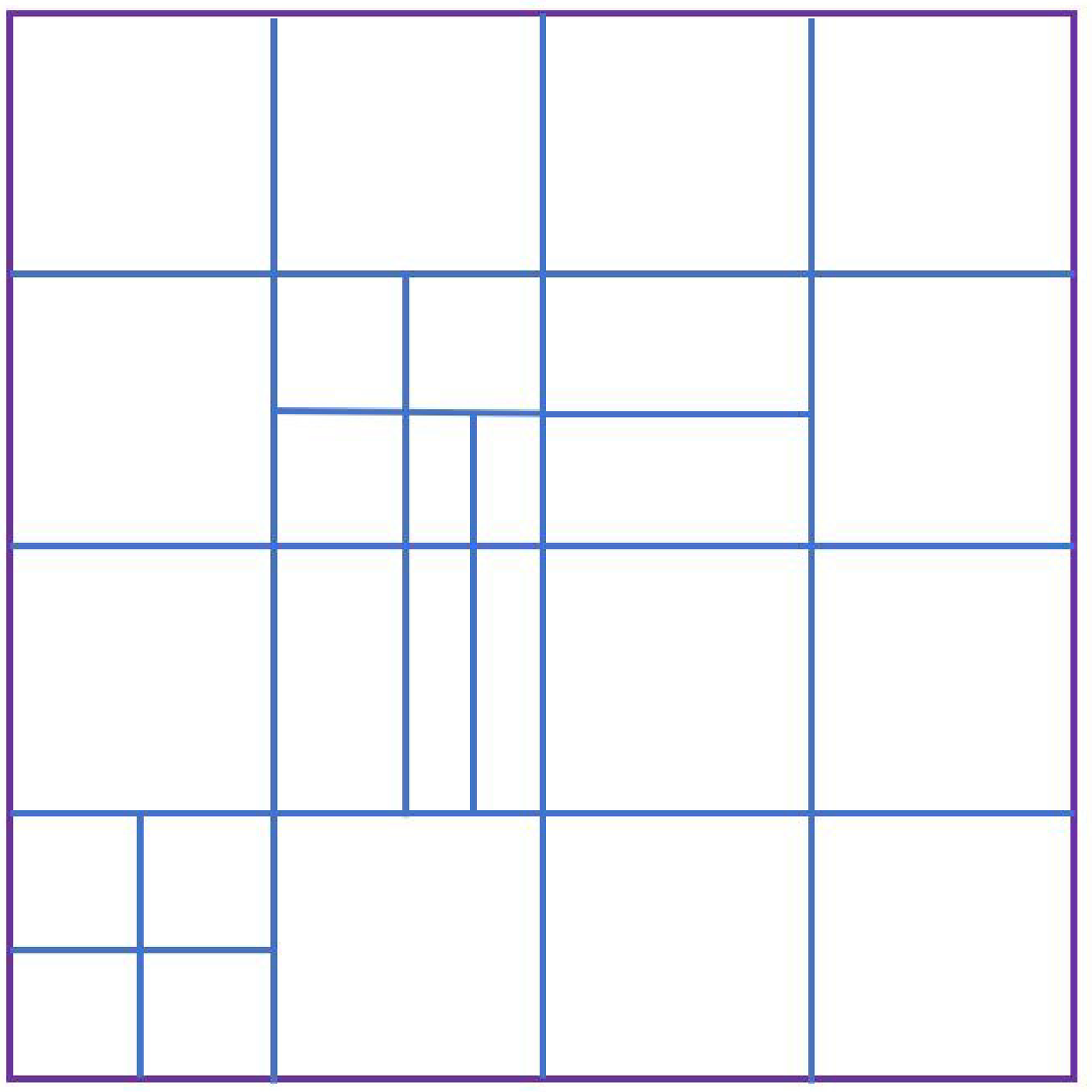}}
\caption{An example of modified hierarchical T-mesh.}
\label{fig:ght}
\end{figure}

To improve the flexibility of hierarchical T-meshes, we modify hierarchical T-meshes as follows.
\begin{definition}\label{def-modtmesh}
A modified hierarchical T-mesh is generated as follows.
In the recursive step (\rmnum{2}) of hierarchical T-meshes, besides cross insertion, cells in $\Theta_k$ can be subdivided by inserting a single (horizontal or vertical) edge in half.
\end{definition}
\reffig{fig:ght} is an example of the refinement process of a modified hierarchical T-mesh.

It should be noted that a strategy of type selection for cells will be provided in Section \ref{sub-ref-way}.

\subsection{Spline spaces over modified hierarchical T-meshes}
Given a T-mesh $\mathbb{T}$, $\mathcal{F}$ represents the set of all the cells in $\mathbb{T}$ and $\Omega$ represents the region occupied by $\mathcal{F}$. The polynomial spline space over $\mathbb{T}$ is defined as
\begin{align*}
\mathcal{S}(m,n,\alpha,\beta,\mathbb{T}):=\{s(x,y)\in C^{\alpha,\beta}(\Omega)\arrowvert s(x,y)\vert_{\phi}\in \mathbb{P}_{mn}\ \mbox{for any } \phi \in \mathcal{F}\}
\end{align*}
where $\mathbb{P}_{mn}$ is the space of all the polynomials with bi-degree $(m,n)$, and $C^{\alpha,\beta}(\Omega)$ is the space consisting of all the bivariate functions that are continuous in $\Omega$ with order $\alpha$ in the $x$-direction and with order $\beta$ in the $y$-direction.

An explicit formula of the spline space $\mathcal{S}(m,n,\alpha,\beta,\mathbb{T})$ in the case of $m\ge2\alpha+1$, $n\ge2\beta+1$ is provided in \cite{deng2006}. For the spline space $\mathcal{S}(3,3,1,1,\mathbb{T})$ with $\mathbb{T}$ being a modified hierarchical T-mesh, the dimension formula can be simplified into
\begin{equation}
\label{equdim}
\dim \mathcal{S}(3,3,1,1,\mathbb{T})=4(V^{b}+V^{+}),
\end{equation}
where $V^{b}$ and $V^{+}$ denote the number of boundary vertices and interior crossing vertices in $\mathbb{T}$ respectively. The dimension formula \eqref{equdim} implies that every boundary vertex or interior crossing vertex corresponds to four basis functions. Therefore, following the method in \cite{deng2008}, we call a boundary vertex or an interior crossing vertex a \emph{basis vertex}.

 For brevity, the polynomial spline in $\mathcal{S}(3,3,1,1,\mathbb{T})$ defined over a modified hierarchical T-mesh is called \emph{modified PHT-spline}.

\section{Basis functions of polynomial splines over modified hierarchical T-meshes}\label{section-basis-constuct}
Inspired by \cite{deng2008}, basis functions of the modified PHT-splines can be constructed in a level-by-level approach.


For the initial level $\mathbb{T}_0$, the standard bicubic $C^1$ continuous tensor-product B-splines are used as basis functions. Let the knot vector of a $C^1$ continuous cubic spline be
 \begin{align*}
&\Xi = [s_0, s_0, s_1, s_1, s_2, s_2, \ldots, s_{m-1}, s_{m-1}, s_{m}, s_{m}],&
\end{align*}
where $s_i<s_{i+1}, 1 \leq i \leq m-2$, $s_0 = s_1$, and $s_{m-1} = s_{m}$.  For every interior knot $s_i$, there are two basis functions with two knot vectors $[s_{i-1}, s_{i-1}, s_i, s_i, s_{i+1}]$, $[s_{i-1}, s_{i}, s_i, s_{i+1}, s_{i+1}]$. Hence these two basis functions have support $[s_{i-1}, s_{i+1}]$. Extending this fact to the tensor-product case of the surface, every vertex of $\mathbb{T}_0$ is associated with four basis functions with  knot vectors
\begin{align*}
&[s_{i-1}, s_{i-1}, s_i, s_i, s_{i+1}] \times [t_{i-1}, t_{i-1}, t_i, t_i, t_{i+1}], &\\
&[s_{i-1}, s_{i}, s_i, s_{i+1}, s_{i+1}] \times [t_{i-1}, t_{i-1}, t_i, t_i, t_{i+1}], &\\
&[s_{i-1}, s_{i-1}, s_i, s_i, s_{i+1}] \times [t_{i-1}, t_{i}, t_i, t_{i+1}, t_{i+1}], &\\
&[s_{i-1}, s_{i}, s_i, s_{i+1}, s_{i+1}] \times [t_{i-1}, t_{i}, t_i, t_{i+1}, t_{i+1}]&
\end{align*}
respectively.

\subsection{Mesh refinement at level $k$}
\label{sub-ref-way}
\begin{figure}[htpb!]
\centering
\subfloat[$\mathbb{T}_0$]{
\includegraphics[width=1.2in]{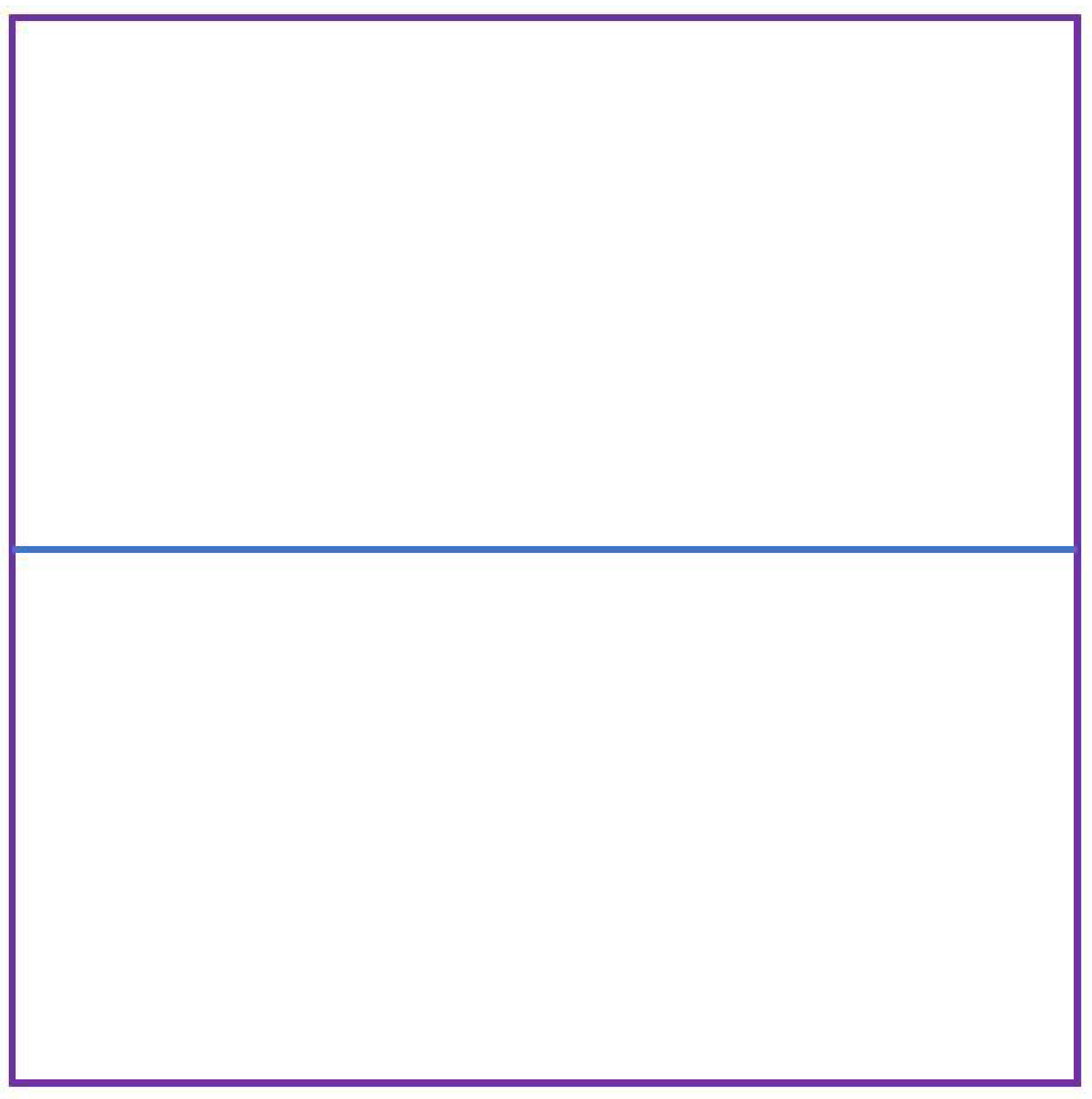}}\hspace{1cm}
\subfloat[$\mathbb{T}_1$]{
\includegraphics[width=1.2in]{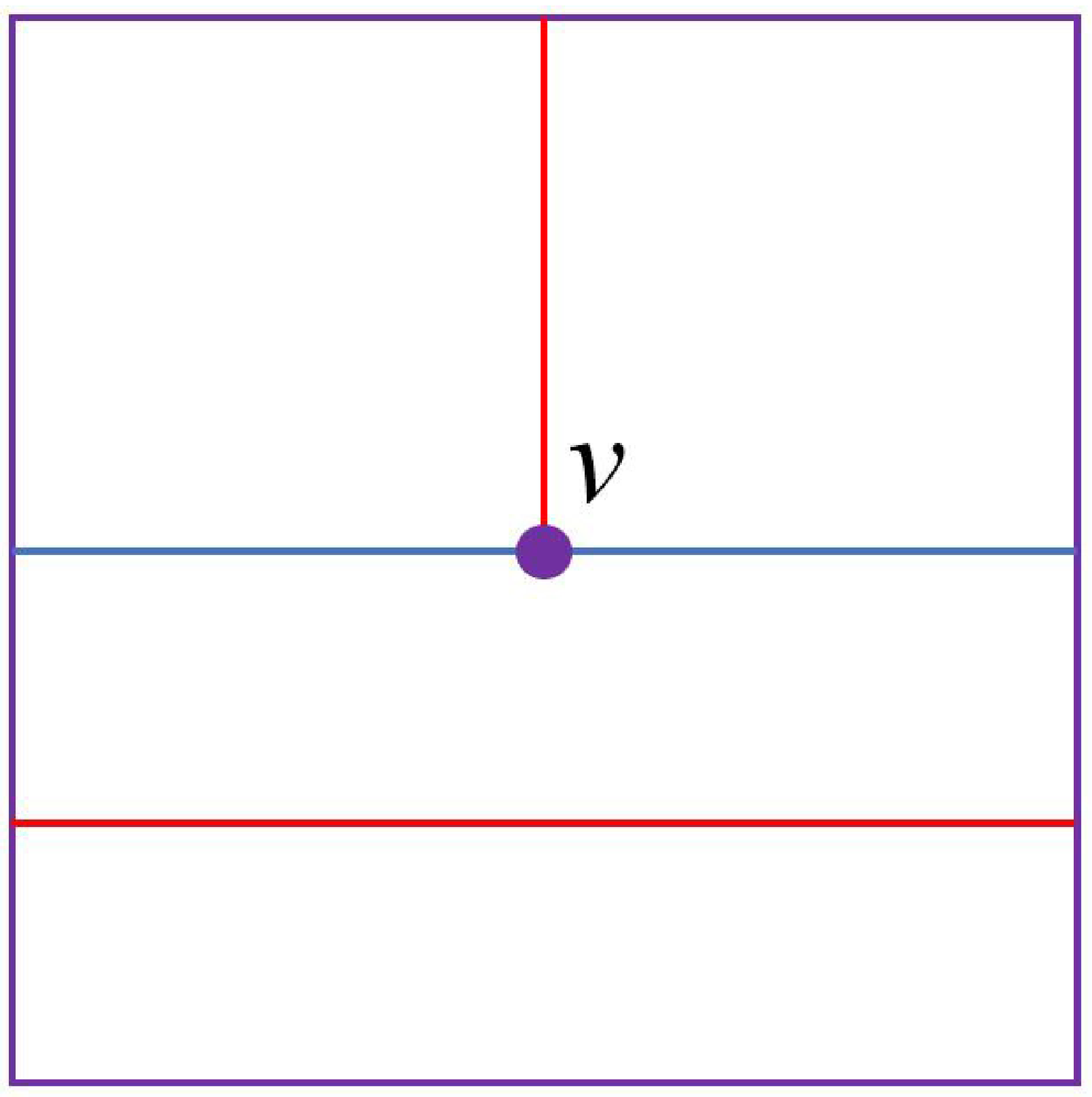}}\hspace{1cm}
\subfloat[$\mathbb{T}_2$]{
\includegraphics[width=1.2in]{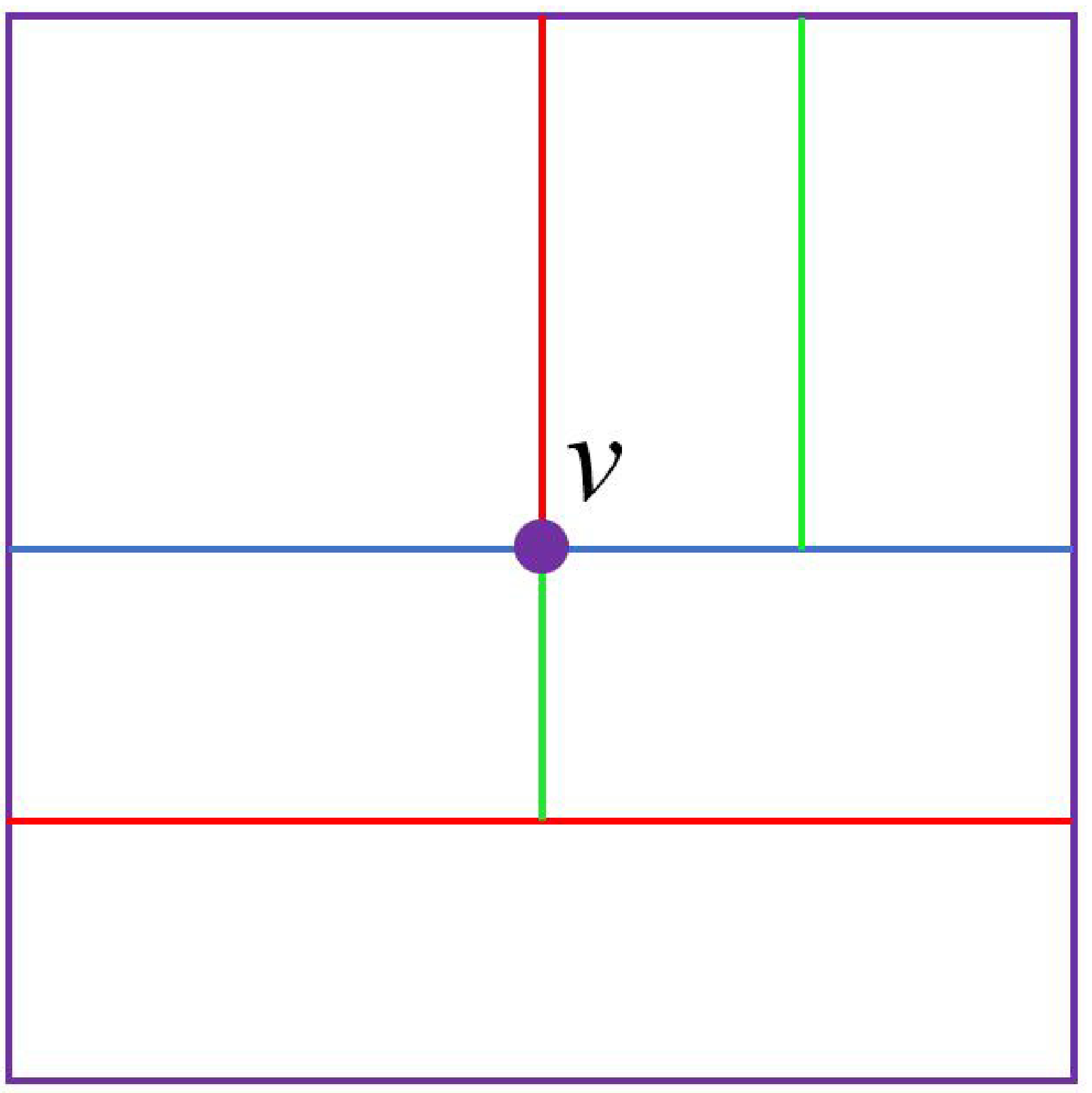}}
\caption{From level 1 to level 2, a T-vertex $\mathbf{v}$ changes into a crossing vertex. }
\label{fig:tr1}
\end{figure}

\begin{figure}[htpb!]
\centering
\subfloat[$\mathbb{T}_0$]{
\includegraphics[width=1.2in]{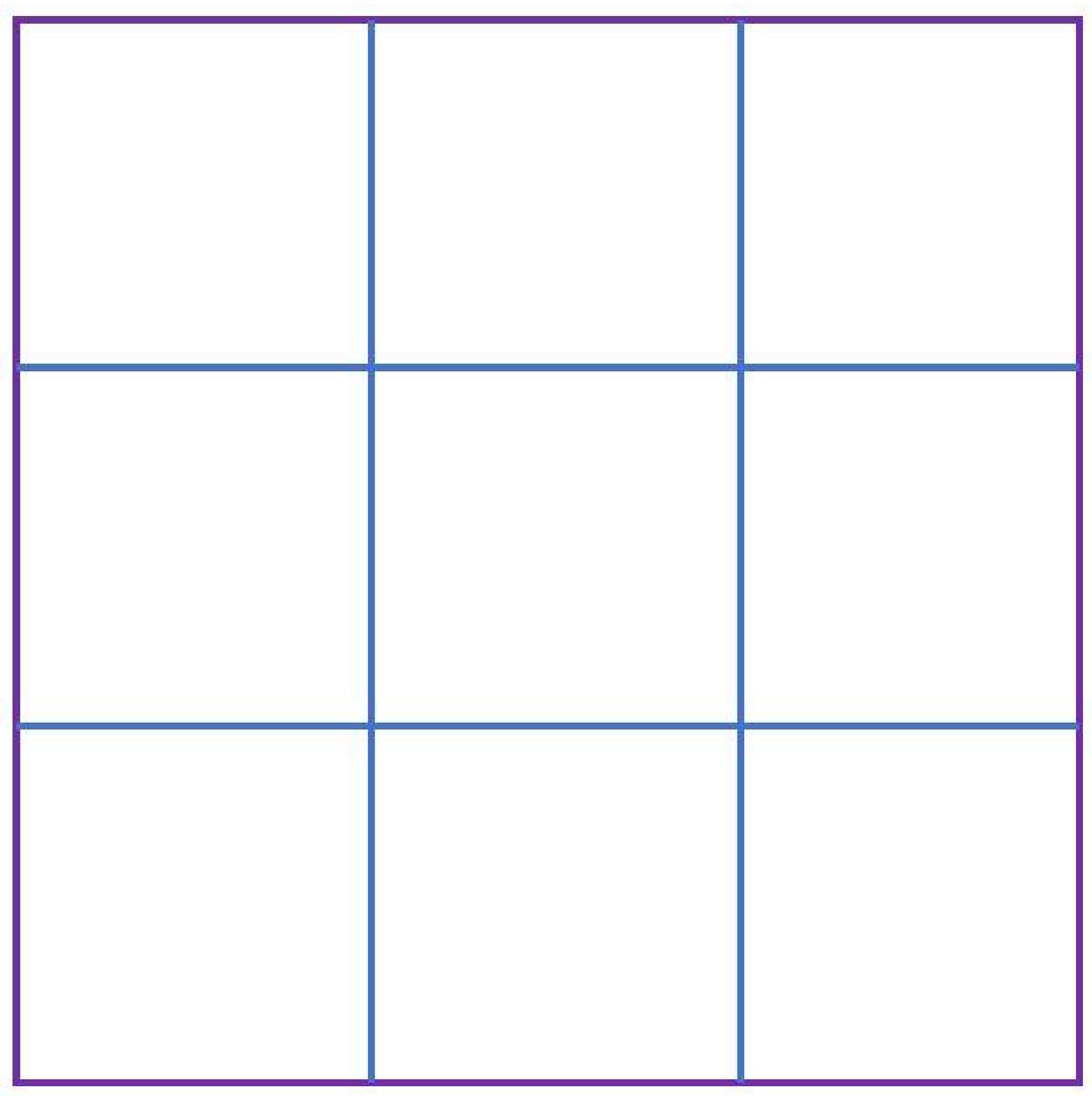}}\hspace{2cm}
\subfloat[$\mathbb{T}_1$]{
\includegraphics[width=1.2in]{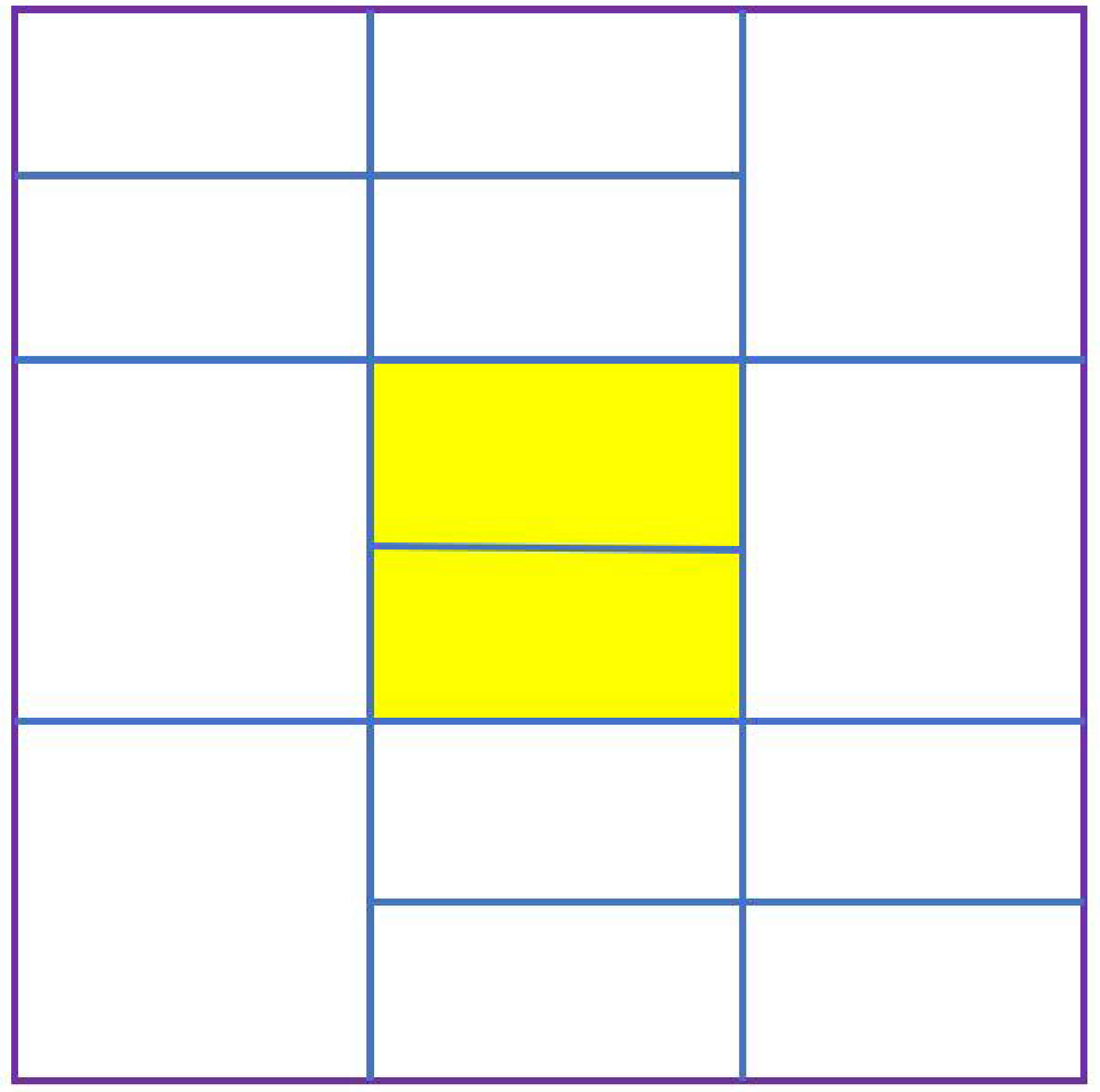}}
\caption{From level 0 to level 1, no basis vertex appears in the yellow cell.}
\label{fig:tr2}
\end{figure}

Assume $\Theta_k$ is given at level $k$ first. For each cell in $\Theta_k$, its anisotropic information is also provided as labels, which indicates the subdivision type of the cell that tends to be selected.
Namely, for each cell $\theta\in\Theta_k$, it has been labeled by `H', `V', or `C' based on anisotropic estimation. Here, `H', `V', and `C' represent horizontal subdivision, vertical subdivision, and cross insertion respectively, which are three different subdivision types that $\theta$ tends to choose. For brevity, the label of the cell $\theta$ is denoted by label$(\theta)$ in the rest of the paper.
For more details about the estimation of anisotropic information, see Section \ref{sub-label-cell} and Section \ref{sub-iga-label}. If we refine these cells directly as their labels indicate, then the following two issues may occur, which increase the complexity in the construction of the basis functions.
\begin{enumerate}[(i)]
  \item Some T-vertices in $\mathbb{T}_{l}$ $(l<k)$ change into crossing vertices in $\mathbb{T}_{k}$ (see \reffig{fig:tr1} for an example).
  \item No basis vertex appears for some cells belonging to $\Theta_k$ after subdivision (see \reffig{fig:tr2} for an example).
\end{enumerate}

With the first issue, the basis construction is slightly complicated,
and we need to revisit the levels at which these T-vertices first appear,
which is required in the modification of basis functions at a later stage.
With the second issue, invalid refinement will occur after subdivision, which contradicts the most common situation in geometric modeling and adaptive isogeometric analysis applications.

Hence, in the current paper,  we explore the construction of polynomial spline in $\mathcal{S}(3,3,1,1,\mathbb{T})$, which is defined over modified hierarchical T-meshes in the absence of the above two issues.
To resolve the above two issues, we provide a refinement strategy to guide subdivision in the following Algorithm 1. The refinement strategy is based on anisotropic features and neighbor relations of cells.

\medskip
\hrule
\medskip
\noindent
\textbf{Algorithm} 1.  Anisotropic Refinement Strategy\\
\textbf{Input}: $\mathbb{T}_k$, $\Theta_k$, $\{\mbox{label}(\theta)\mid \theta\in\Theta_k \}$\\
\textbf{Output}: $\mathbb{T}_{k+1}$ \\
\begin{enumerate}[1)]
 \item Classify $\Theta_k$ into several connected groups $G_1^k,\cdots, G_M^k$ (See Remark \ref{remark1} for more details).
 \item For each connected group $G_i^k$, check label$(\theta)$ for all $\theta\in G_i^k$.
     \begin{enumerate}[2.1)]
       \item When labels are identical, check whether new basis vertices will appear for each $\theta \in G_i^k$ or not. \\
           If no, for each cell $\theta_0\in G_i^k$ without new basis vertices, find all cells that are in $G_i^k$ and connected to $\theta_0$ by a sequence of vertically (or horizontally, if there exists) aligned-adjacent cells in $G_i^k$. Relabel them with `C'.
       \item When labels are nonidentical, check neighbor cells for each $\theta \in G_i^k$.
       \begin{enumerate}
         \item If $\theta$ only has horizontally aligned-adjacent cells $\overline{\theta} \in G_i^k$, $\mbox{label}(\theta) =$ `H';
         \item If $\theta$ only has vertically aligned-adjacent cells $\overline{\theta} \in G_i^k$, $\mbox{label}(\theta) =$ `V';
         \item if $\theta$ has both vertically and horizontally aligned-adjacent cells or has not any aligned-adjacent cell $\overline{\theta} \in G_i^k$, $\mbox{label}(\theta) = $ `C'.
       \end{enumerate}

     \end{enumerate}
 \item Subdivide cells as new labels indicates. For each cell $\theta \in G_i^k$,
 \begin{enumerate}[3.1)]
   \item if $\theta$ has the label `H', split $\theta$ in half horizontally;
   \item if $\theta$ has the label `V', split $\theta$ in half vertically;
   \item if $\theta$ has the label `C', subdivide $\theta$ by inserting a cross.
 \end{enumerate}
 \item Output $\mathbb{T}_{k+1}$.
\end{enumerate}
\hrule
\medskip

\begin{rmk}\label{remark1}
In Step 1) of Algorithm 1, the classification is achieved by "flood-fill" through aligned-adjacent relations. Namely, two cells of level $k$ will be classified into one connected group,  1) if they are adjacent, they must be aligned-adjacent; 2) if they are not adjacent, there exists a sequence of aligned-adjacent cells in $\Theta_k$ connecting them. Namely, for two cells $\theta_s, \theta_e \in \Theta_k$, if they are not adjacent, there exist $\hat{\theta}_j^k \in \Theta_k, j=1,\ldots, L$ such that
$\hat{\theta}_1^k =\theta_s$, $\hat{\theta}_L^k =\theta_e$, $\hat{\theta}_j^k$ and $\hat{\theta}_{j+1}^k$ are aligned-adjacent, $j=1,\ldots, L-1$.
\end{rmk}

\begin{rmk}\label{remark2}
In Step 2.1) of Algorithm 1, cells in $G_i^k$ are relabeled based on anisotropic information if their original labels are identical.
Otherwise, cells in $G_i^k$ are relabeled in Step 2.2) based on neighbor relations.
\end{rmk}

\begin{figure}[htbp!]
\centering
\subfloat[$\mathbb{T}_0$ with given labels]{
\includegraphics[width=1.2in]{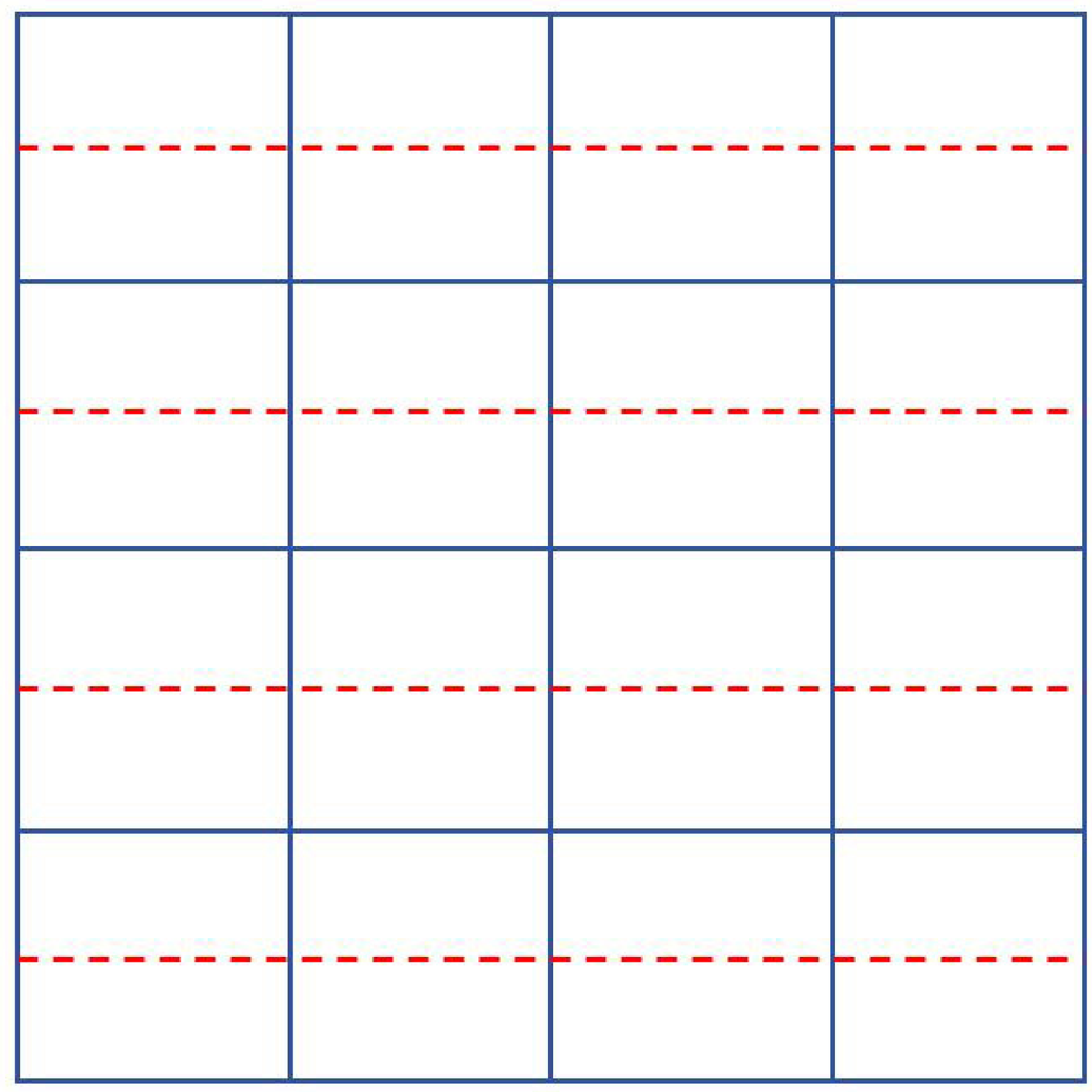}}\hfill
\subfloat[$\mathbb{T}_1$]{
\includegraphics[width=1.2in]{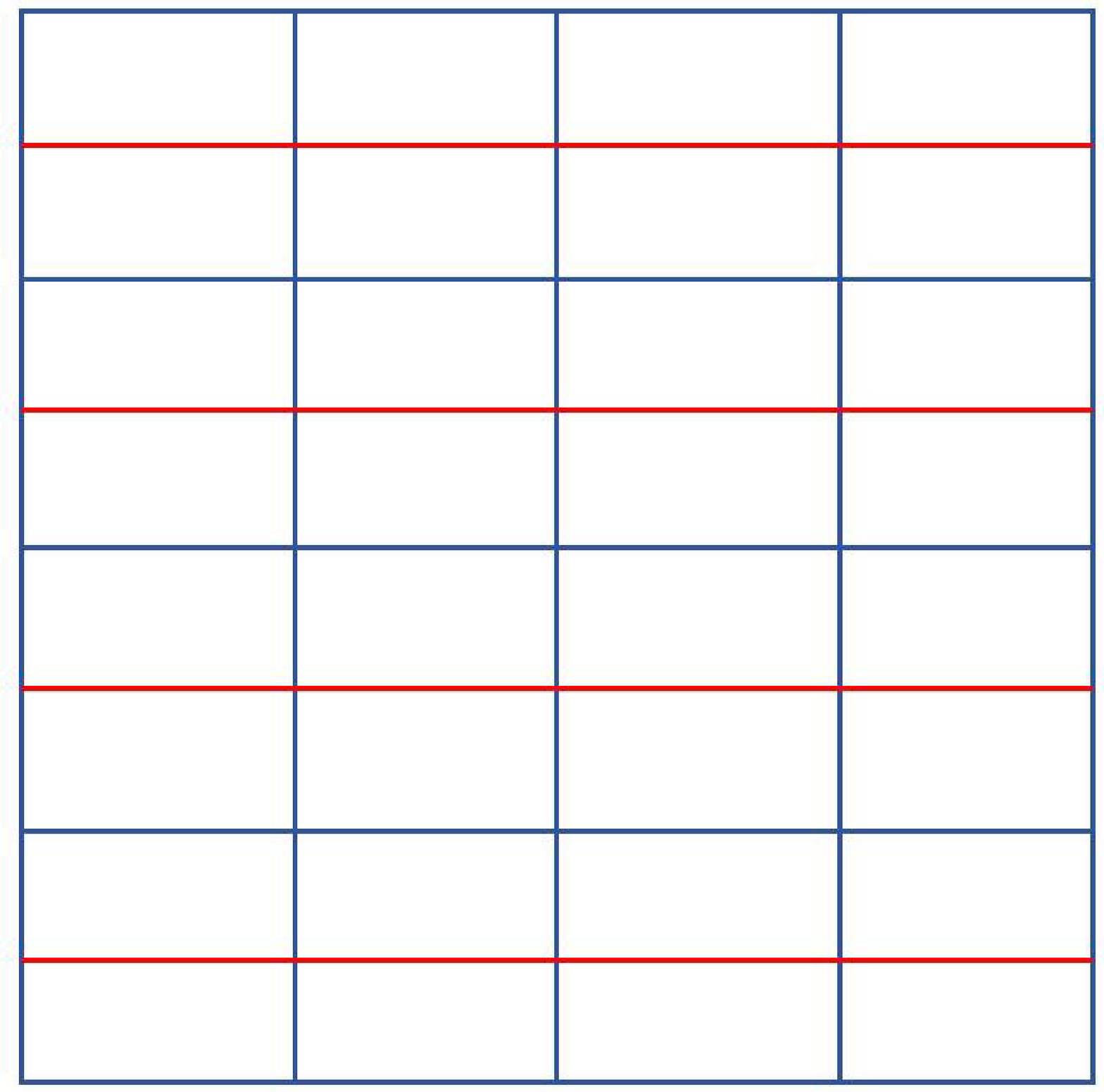}}\hfill
\subfloat[$\Theta_1$ is classified into three connected groups]{
\includegraphics[width=1.2in]{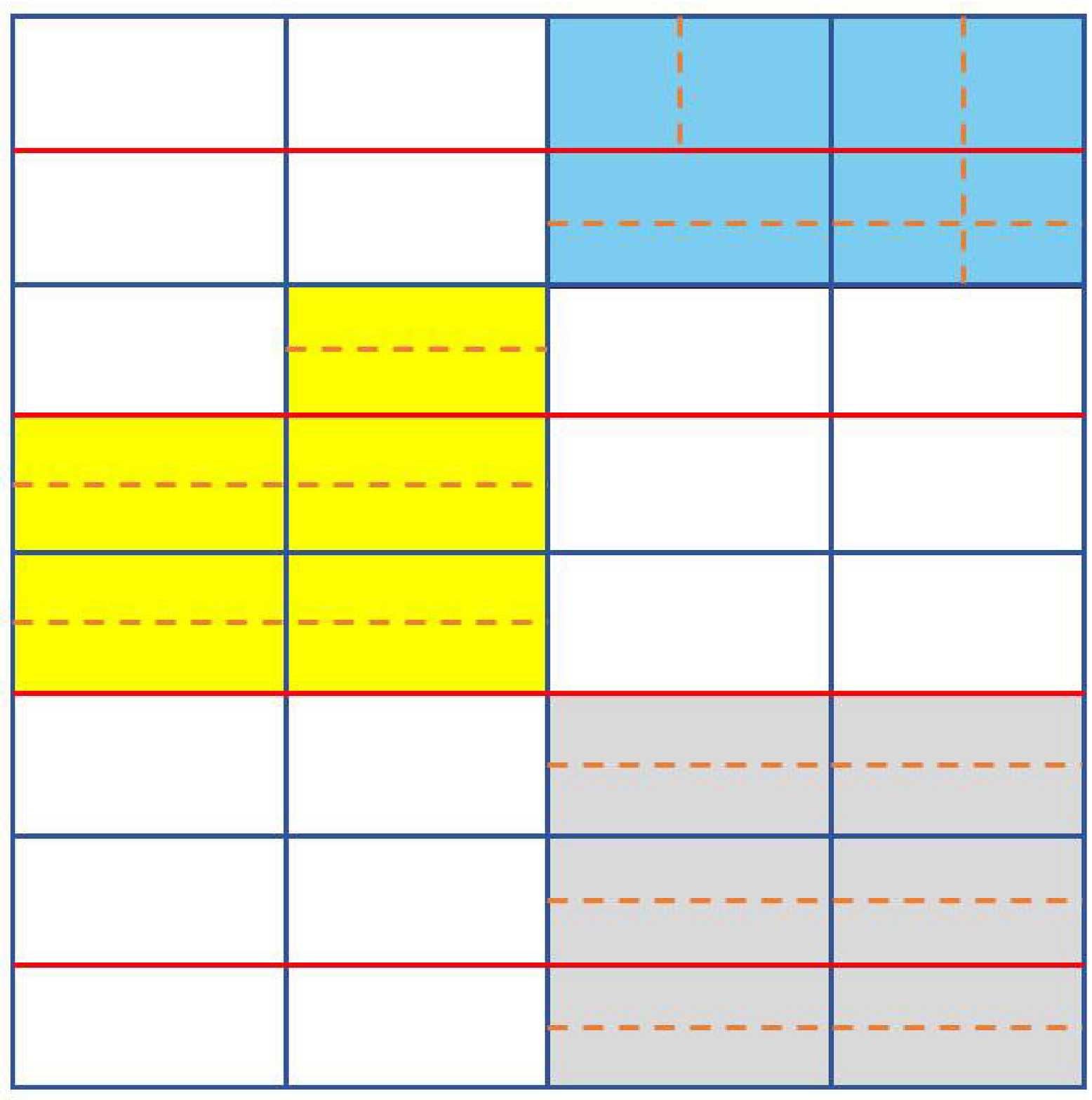}}\hfill
\subfloat[$\mathbb{T}_2$]{
\includegraphics[width=1.2in]{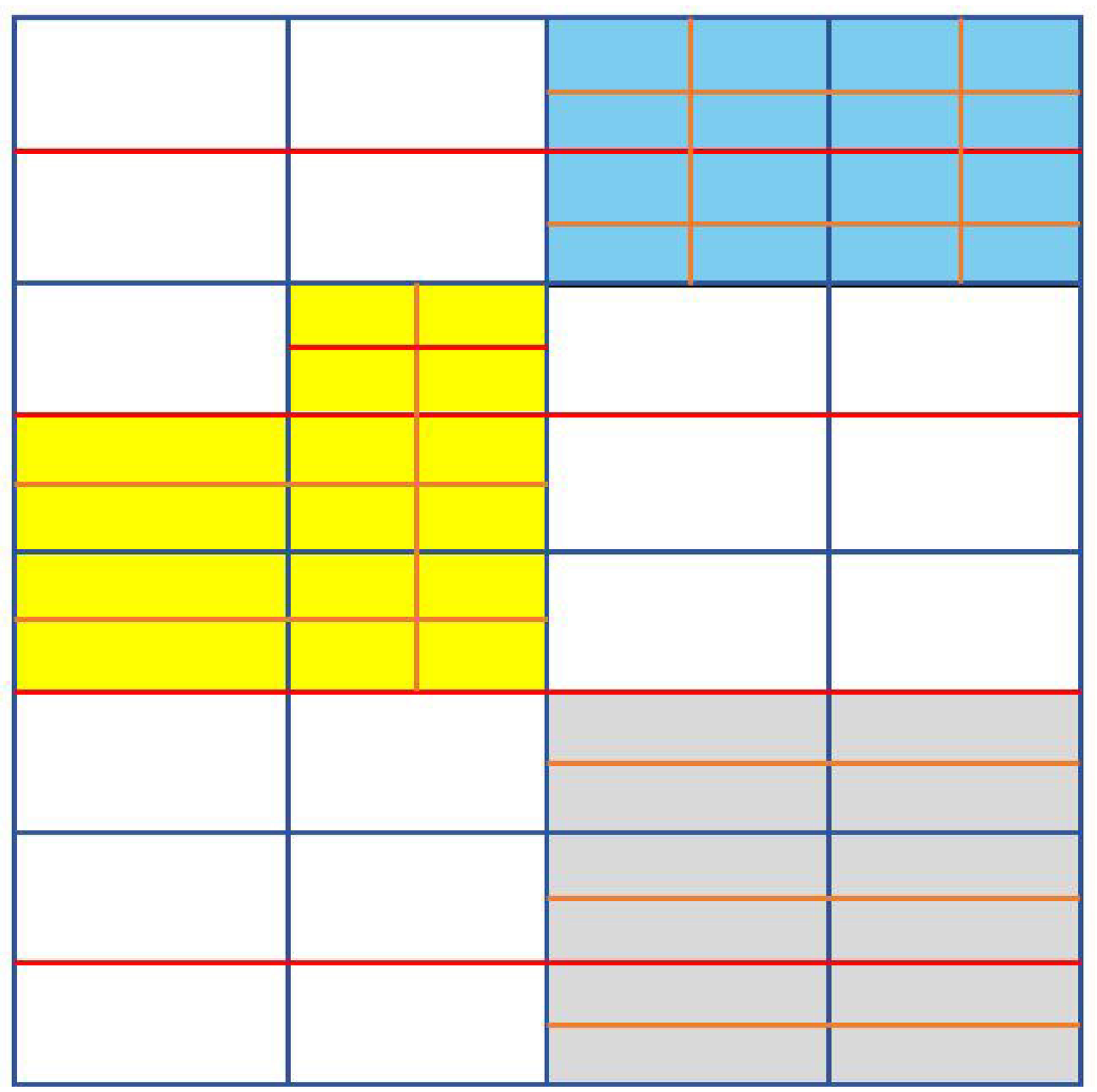}}
\caption{An example of subdivision result based on our strategy. The dashed edges in the cells are candidate subdivision types tend to choose.}
\label{fig:ex-pght}
\end{figure}

An example in \reffig{fig:ex-pght} illustrates our subdivision strategy.
All cells in an initial tensor-product mesh $\mathbb{T}_0$ have labels `H'(see \reffig{fig:ex-pght}(a)). After executing Algorithm 1 once, we get $\mathbb{T}_1$ in
\reffig{fig:ex-pght}(b).
Cells to be subdivided in $\mathbb{T}_1$ are classified into three connected groups (see \reffig{fig:ex-pght}(c)).
If we subdivide the cells directly as their labels indicate,
no basis vertex will appear in the upper right cell in the yellow group,
and T-vertices may change into basis vertices after subdividing cells several times (blue parts).
By Step 2) and Step 3) of our strategy,  we get $\mathbb{T}_2$ in \reffig{fig:ex-pght}(d).

The connected group plays an important role in Algorithm 1. Here we shall provide the properties of the connected groups.
\begin{lem}\label{lemma1}
Suppose $G$ is a connected group. Subdivide all cells in $G$ as in Algorithm 1, and let $G'$ be the set consisting of all new cells. Then:
\begin{enumerate}
  \item There is no T-vertex on the interior edges of the group $G$.
  \item $G'$ forms a new connected group.
  \item There is no T-vertex on the interior edges of the group $G'$.
  \item Let $G_1$ be a subset of $G'$.  Assume $G_1$ is classified into several connected groups (based on Remark \ref{remark1}).
  Then, there is no common edge between any two connected groups.

\end{enumerate}
\end{lem}
\begin{pf}
1. There is no T-vertex on the interior edges of $G$ before subdivision.
Otherwise, three cells around the interior T-vertex, which are in one connected group, are adjacent to each other but not aligned-adjacent,
which contradicts the assumption that they are from one connected group.

2. Consider any two adjacent cells, denoted by $\theta_0$ and $\theta_1$, in $G$ first.
By Remark \ref{remark1},  $\theta_0$ and $\theta_1$ are also aligned-adjacent.
Let the common edge of $\theta_0$ and $\theta_1$ be $\mathbf{e}$.
Assume that $\theta_0$ and $\theta_1$ are horizontally (vertically) adjacent,
then a T-vertex will appear on the common edge $\mathbf{e}$ if
\begin{enumerate}[(i)]
  \item one is subdivided by inserting a cross, and the other is split in half vertically (horizontally), or
  \item one is split in half vertically, and the other is split in half horizontally.
\end{enumerate}
However, by Step 2.1) and Step 2.2) in Algorithm 1, neither of the above cases will occur after subdivision.
Hence, any new cells acquired by subdividing two adjacent cells will be aligned-adjacent if they are adjacent.
If new cells are not adjacent, there exists a sequence of new aligned-adjacent cells connecting them.
Consequently, all these new cells will be classified into the same connected group.

Repeat the above process for any two adjacent cells in $G$. We get that all new cells,
which are generated by subdividing cells in $G$ as Algorithm 1, are in the same connected group.

3. This result is a direct consequence of property 1 and property 2 in Lemma \ref{lemma1}.

4. Note that there is no common cell between any two different connected groups. Otherwise,
these two connected groups will merge into a bigger connected group through the common cells.
Hence, suppose there exists a common edge $\mathbf{e}$ between two different connected groups.
Then, the common edge $\mathbf{e}$ is shared by two cells, which are from different connected groups.
By property 2 in Lemma \ref{lemma1}, $G'$ is a connected group, and it follows that these two cells are aligned-adjacent.
Therefore, these two connected groups should be classified into the same connected group,
which contradicts our assumption that they are classified into different connected groups. \qed
\end{pf}

\begin{thm}\label{thm-algr}
Start from a tensor-product mesh $\mathbb{T}_0$, let $\mathbb{T}_{N}$ be the mesh output by performing  Algorithm 1 $N$ times ($N>0$).
Then, $\mathbb{T}_{N}$ is a modified hierarchical T-mesh. Hence, for each level,
\begin{enumerate}
  \item no T-vertex changes into a crossing vertex;
  \item new basis vertices will appear for each cell that is to be subdivided.
\end{enumerate}
\end{thm}
\begin{pf}
In Algorithm 1, cells are allowed to be subdivided by splitting in half or inserting crosses.
It is straightforward to check that $\mathbb{T}_{N}$ is a modified hierarchical T-mesh.


1. For any connect group $G_i^k$ at level $k$, by property 1 of Lemma \ref{lemma1}, there is no T-vertex on the interior edges of the group $G_i^k$ before subdivision.
Therefore, the change from a T-vertex to a crossing vertex on the interior edges of the group is impossible.

On the other hand, since we start with a tensor product mesh $\mathbb{T}_0$,
it follows by Property 2 of Lemma \ref{lemma1} that there is no common edge for any two different connected groups at level 0.
Note that by Remark \ref{remark-celllevel}, if a cell is not subdivided at level $k$, it will be excluded from subdivision henceforth.
Hence, by property 4 of Lemma \ref{lemma1}, for each time we perform Algorithm 1, any two different connected groups at the same level have no common edge.
Consequently, no T-vertex on the boundary edges of connected groups changes into a crossing vertex at each level.

Therefore,  no T-vertex will change into a crossing vertex if we subdivide cells as Algorithm 1.

2. At each level, when the labels of cells in a connected group $G_i^k$ are identical,
Algorithm 1 will subdivide these cells directly as their labels indicate when new basis vertices will appear for each cell.
The exceptions are
\begin{enumerate}[(a)]
  \item for those cells with label `H', there are no horizontally neighbor cells in $G_i^k$;
  \item for those cells with label `V', there are no vertically neighbor cells in $G_i^k$.
\end{enumerate}
With the case (a) (case (b)), Algorithm 1 will relabel these cells and their neighbor cells in the same row (column) with `C' as Step 2.1),
which ensures the new basis vertices for these cells.

When the labels of cells in a connected group are not identical, Algorithm 1 will subdivide these cells based on neighboring relations.
It is straightforward to verify that for each cell in $G_i^k$, at least one basis vertex will appear on the edge of this cell.
Thus, the appearance of new basis vertices for each cell, which is to be subdivided, is guaranteed in Algorithm 1.\qed

\end{pf}

\subsection{Basis construction}\label{section-modify}
The concept in constructing basis functions relies on the modifying mechanism that operates on the B\'ezier form of basis functions directly.

\begin{figure}[htbp!]
\centering
\subfloat[New B\'ezier ordinates after the split in half horizontally]{
\includegraphics[width=1.2in]{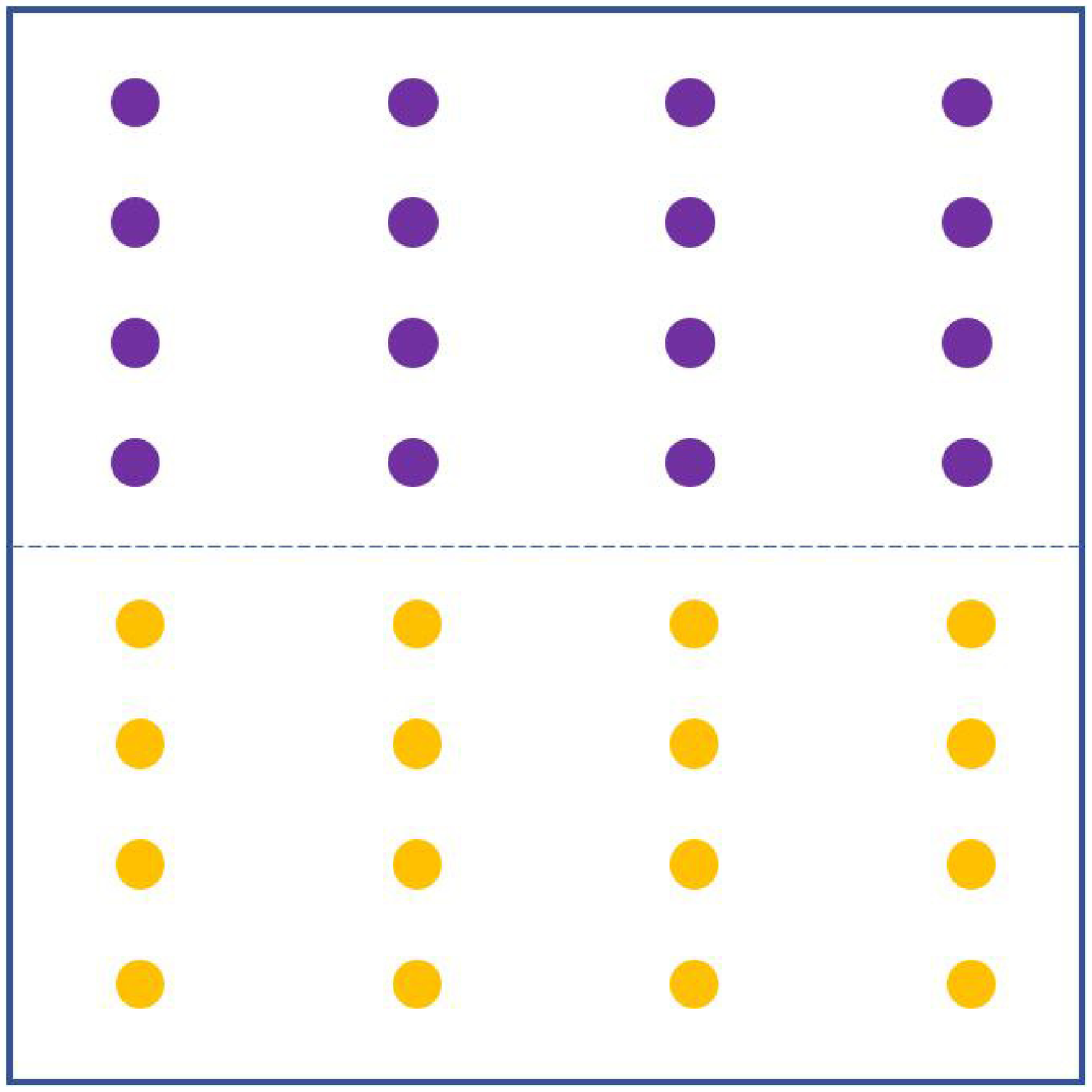}}\hspace{1cm}
\subfloat[New B\'ezier ordinates after the split in half vertically]{
\includegraphics[width=1.2in]{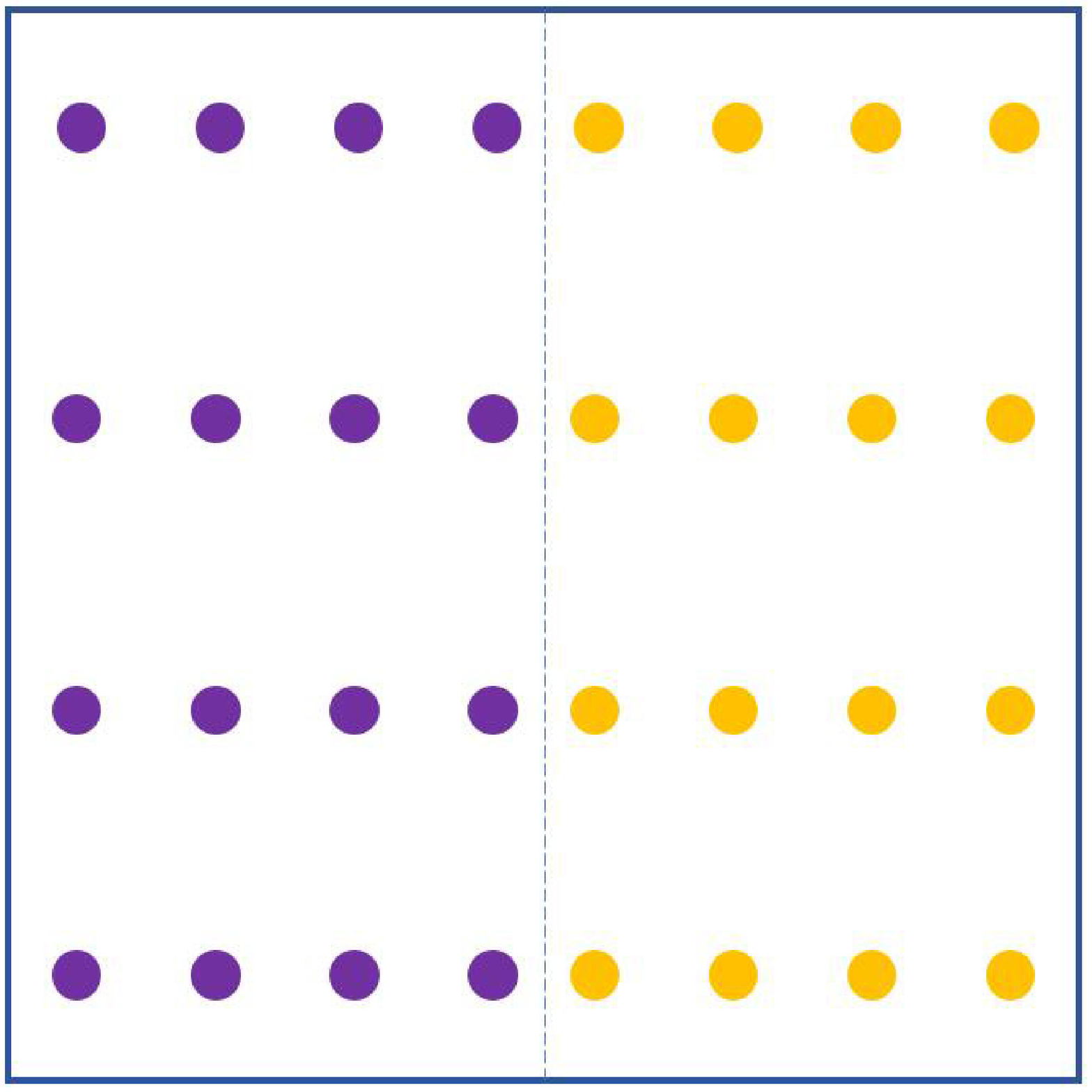}}\hspace{1cm}
\subfloat[New B\'ezier ordinates after cross insertion]{
\includegraphics[width=1.2in]{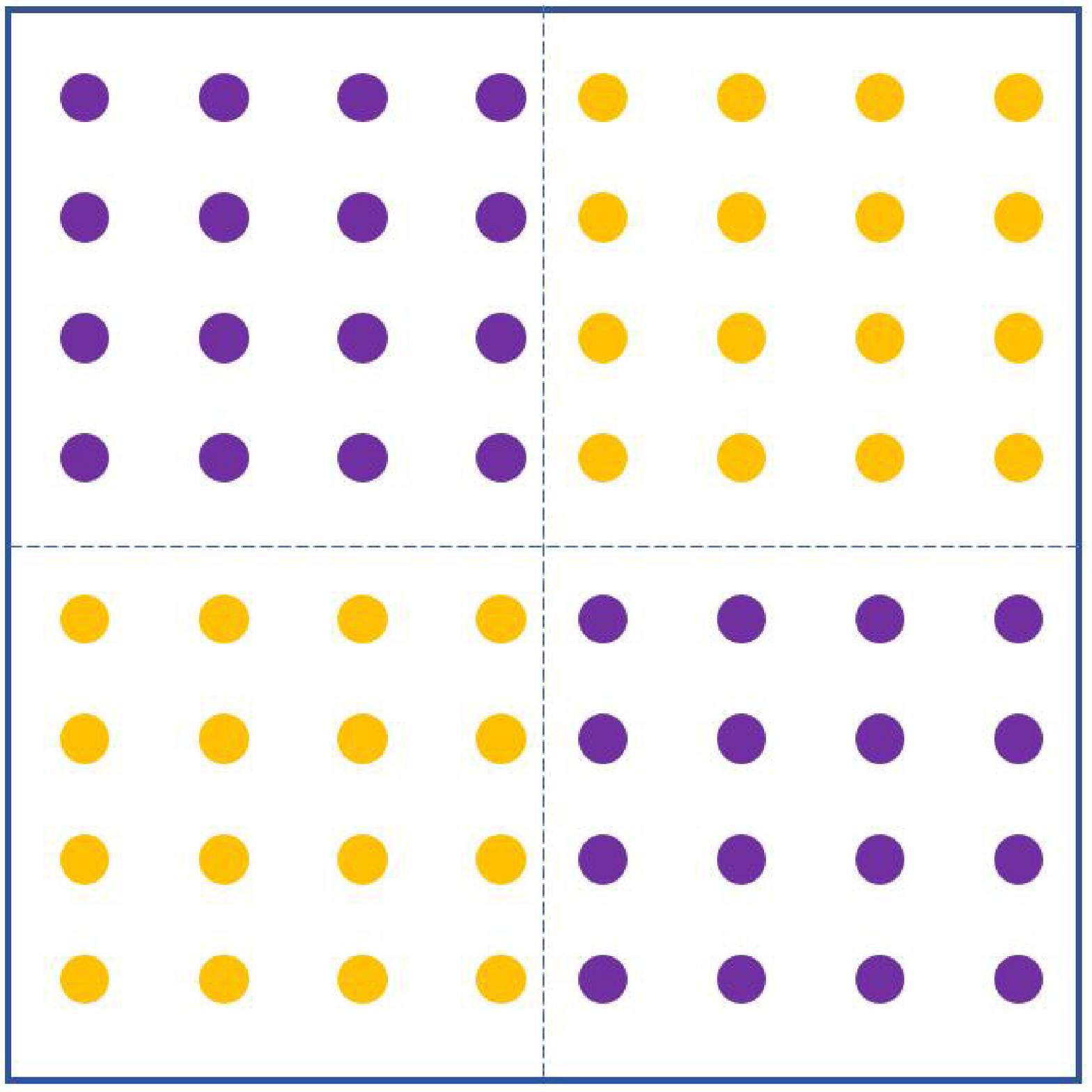}}
\caption{Using de Casteljau's algorithm, new B\'ezier ordinates are generated for three different subdivision types.}
\label{fig:bezier}
\end{figure}

\begin{definition}\label{def-modify}
Given a basis function $b_{i}^{k}(s,t)$ at level $k$, for each cell $\theta$ of $\Theta_k$ with $\theta \subset \mbox{supp} \big(b_{i}^{k}(s,t)\big)$.
  The modification of $b_{i}^{k}(s,t)$  from level $k$ to level $k+1$ $$
 \bar{b}_{i}^{k+1}(s,t) =\mathfrak{M}^{k+1} \big( b_{i}^{k}(s,t)\big)
  $$
  is defined as follows.
  \begin{enumerate}[1)]
  \item Represent $b_{i}^{k}(s,t)$ over $\theta$ by specifying 16 B\'ezier ordinates.
  \item Generate new B\'ezier ordinates by applying de Casteljau's algorithm according to the different subdivision type (see \reffig{fig:bezier} for an illustration).
  \item  Set all of the B\'ezier ordinates that are associated with the new basis vertices to zero.
\end{enumerate}
\end{definition}

\begin{figure}[htpb!]
\centering
\includegraphics[width=1.2in]{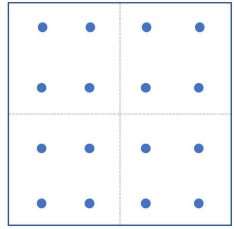}
\caption{The B\'ezier ordinates are associated with four corner vertices}
\label{fig:bez}
\end{figure}

\begin{figure}[htbp!]
\centering
\subfloat[B\'ezier\ ordinates\ over\ a\ cell]{
\includegraphics[width=1.2in]{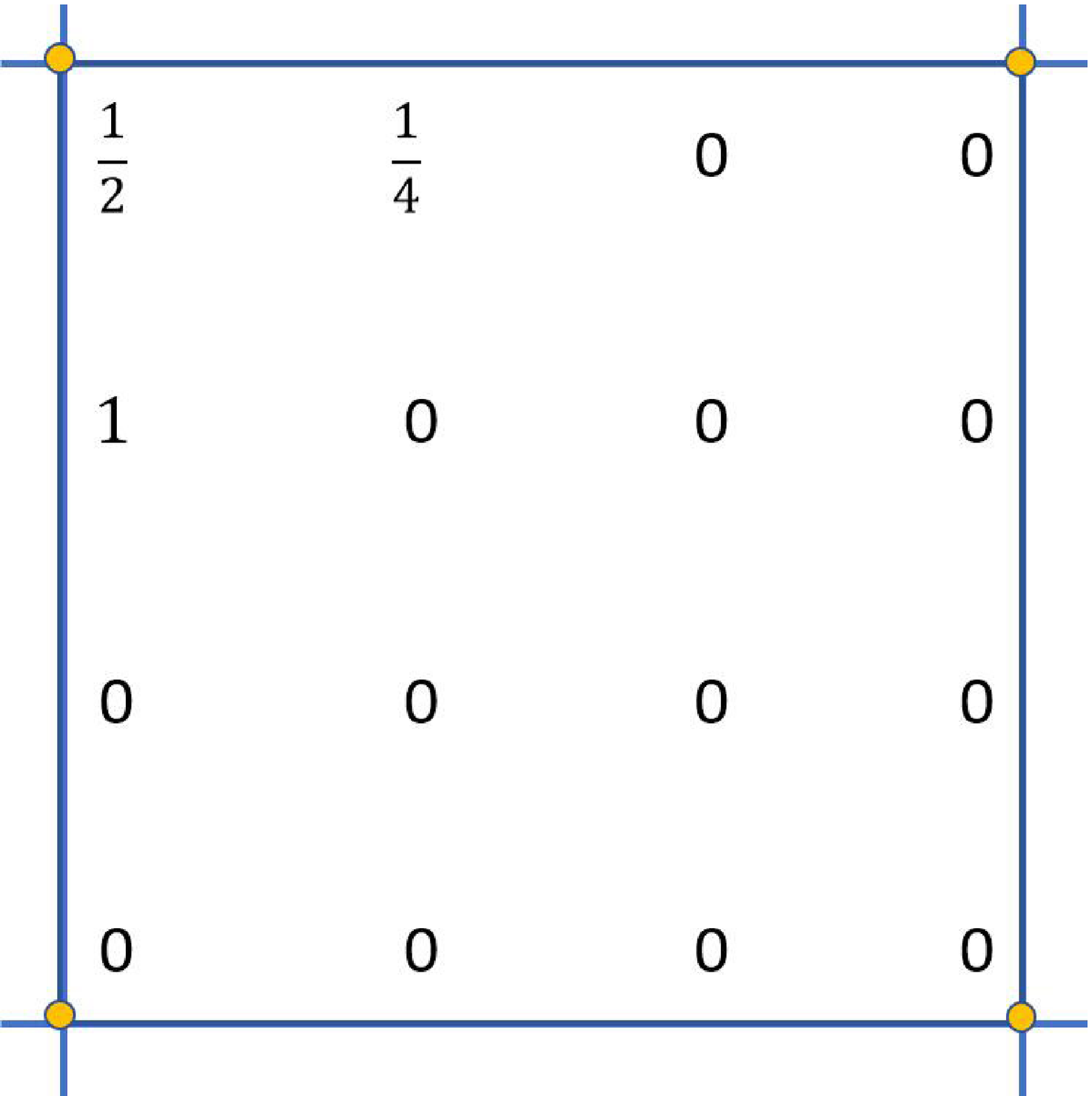}}\hspace{1cm}
\subfloat[The\ function\ is\ subdivided\ into\ two\ subcells,\ after\ adding\ a\ line]{
\includegraphics[width=1.2in]{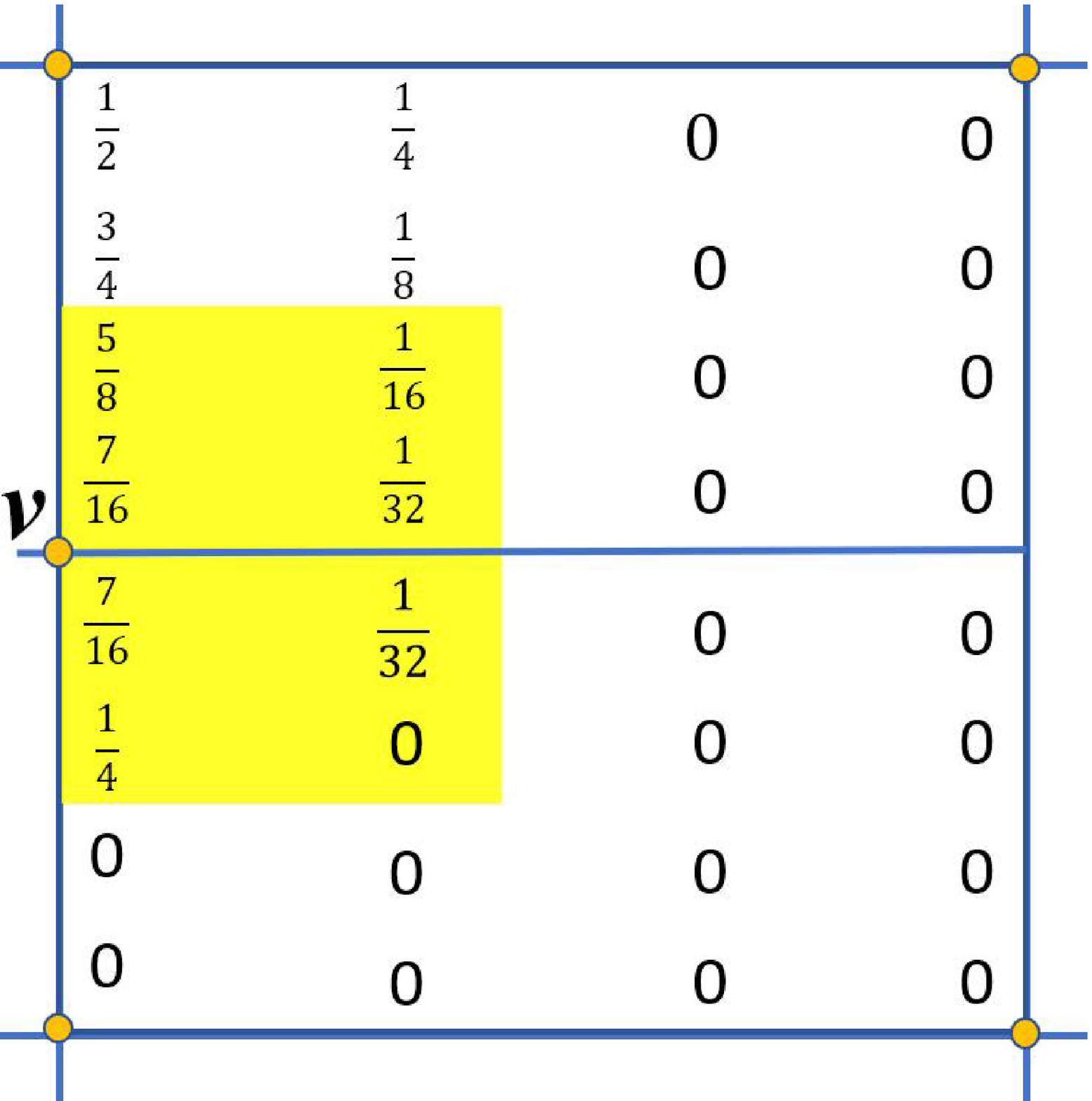}}\hspace{1cm}
\subfloat[Set B\'ezier ordinates associated with new basis vertex to be zero]{
\includegraphics[width=1.2in]{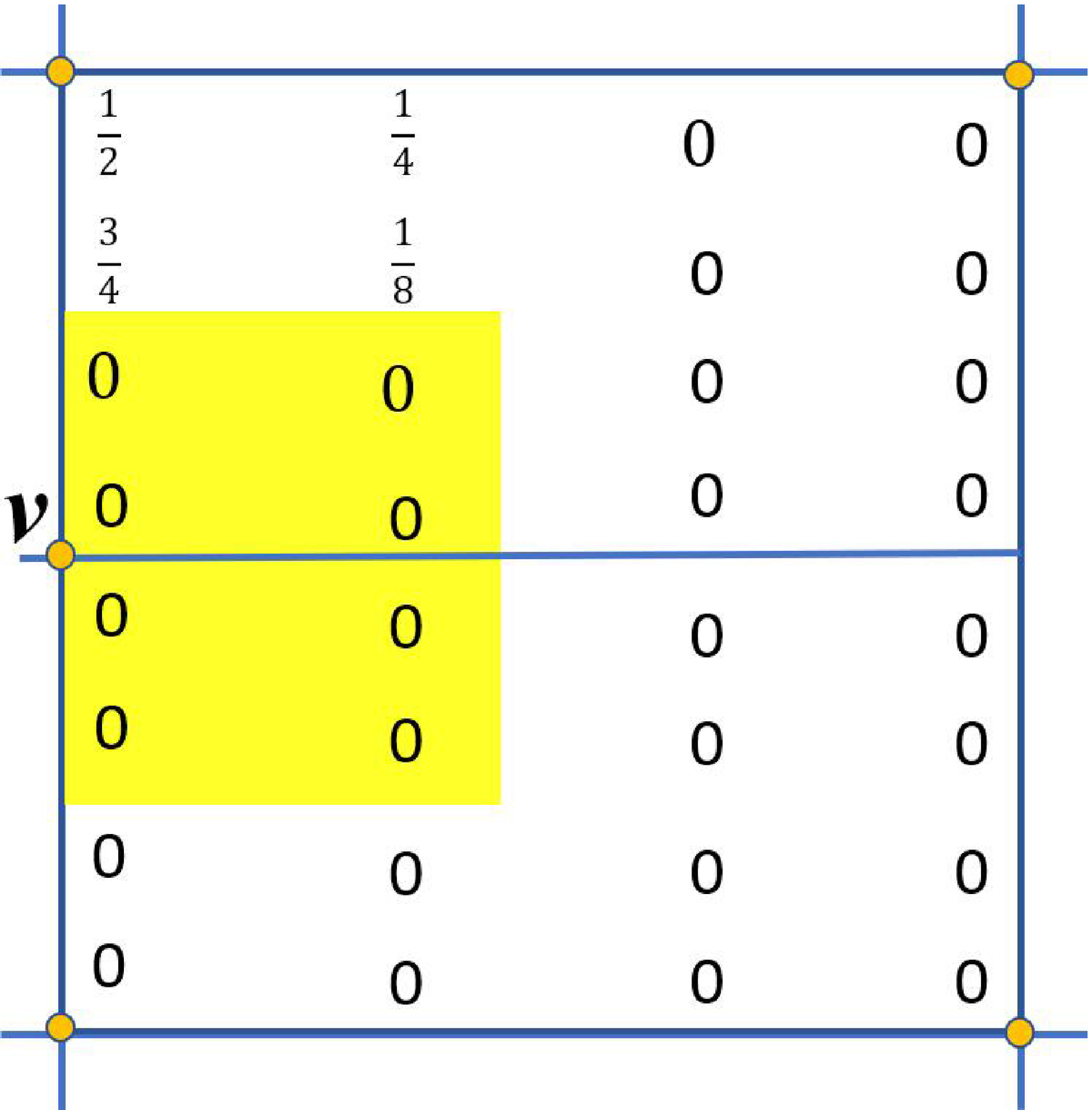}}
\caption{Modification of a basis function, where $\mathbf{v}$ is a new basis vertex.}
\label{fig:modify}
\end{figure}

Note that the function $b_{i}^{k}(s,t)$ remains unchanged if there is no cell $\theta$ in $\Theta_k$ such that $\theta\subset\mbox{supp}\big(b_{i}^{k}(s,t)\big)$.
For each cell that is to be subdivided at level $k+1$, the 16 B\'ezier ordinates are divided into four or two parts according to its subdivision type.
Hence, each part is associated with a cell corner vertex, which is illustrated in \reffig{fig:bez}. By setting all of the B\'ezier ordinates that are associated with the new basis vertices to zero (see \reffig{fig:modify} for an illustration), the basis function $b_{i}^{k}(s,t)$ at level $k$ is modified into a basis function $\bar{b}_{i}^{k}(s,t)$ at level $k+1$.
In addition, the B\'ezier ordinates around the basis vertex  (located in neighbor cells) are reset simultaneously. Hence, the conditions of $C^1$ continuity still hold and  $\bar{b}_{i}^{k+1}(s,t) \in \mathcal{S}(3,3,1,1,\mathbb{T}_{k+1})$.


Next, we shall discuss adding new basis functions at level $k+1$.

By the definition of modified hierarchical T-mesh, a new interior basis vertex appears at level $k+1$ if
\begin{enumerate}[(a)]
  \item a cell at level $k$ is subdivided by cross insertion, which introduces a crossing vertex at the center of the cell;
  \item the common edge of two aligned-adjacent cells is split after subdivision, which introduces a crossing vertex on the common edge.
\end{enumerate}
For the case (a), the neighboring four cells around the new interior vertex are just four subcells of the original cell.
For the case (b), the neighboring four cells around the new interior vertex are also aligned-adjacent.
Therefore, for each new basis vertex, if it is an interior vertex, its neighboring four cells form a $2\times 2$ tensor-product mesh. On the other hand, new boundary basis vertices are introduced by subdividing of the boundary cells at level $k$. Therefore, for each new boundary basis vertex, its neighboring two cells are two aligned-adjacent subcells of the original boundary cell at level $k$. Hence, these two subcells also form a tensor-product mesh.
Thus, for each new basis vertex $\mathbf{v}_{i}^{k+1}$, four new basis functions associated with the basis vertex $\mathbf{v}_{i}^{k+1}$ can be constructed.



We are now ready to construct the basis functions of modified PHT-splines.

\begin{definition}\label{def-mpht-basis}
  Given a modified hierarchical T-mesh $\mathbb{T}$ in the 2D plane. A set of basis functions $b_i(s,t)$ for the spline space $\mathcal{S}(3,3,1,1,\mathbb{T})$ is recursively constructed as follows:
  \begin{enumerate}[\ \ 1. ]
  \item Initialization: $\mathcal{S}^0 =\{b_i^0(s,t):$ tensor-product bicubic B-spline basis function over $\mathbb{T}_0$\}.
  \item Recursive case:
  $S^{k+1} = \mathcal{S}_{A}^{k+1} \cup \mathcal{S}_{B}^{k+1}$ for $k=0,\ldots, N-1$, where
  \begin{align*}
  &\mathcal{S}_A^{k+1} = \{\mathfrak{M}^{k+1}(b_i^k(s,t)):  b_i^k(s,t) \in \mathcal{S}^{k}\},&\\
  &\mathcal{S}_B^{k+1} = \{b_i^{k+1}(s,t):  \mbox{new basis functions at level } k+1 \}.&
  \end{align*}
  \item $\mathcal{S}=\mathcal{S}^{N}$
\end{enumerate}
\end{definition}

Next we shall discuss the properties of the basis functions constructed according to the above definition.
\begin{thm}\label{thm-6}
Let $\mathbb{T}$ be a modified hierarchical T-mesh. The modified PHT-spline basis $\mathcal{S}$ is constructed according to Definition \ref{def-mpht-basis}. Then:
\begin{enumerate}[\ \ a)]
  \item The functions in $\mathcal{S}$ are linearly independent;
  \item $\mathcal{S}$ spans the modified PHT-spline space $\mathcal{S}(3,3,1,1, \mathbb{T})$.
  \item $\mathcal{S}$ forms a partition of unity.
\end{enumerate}
\end{thm}
\begin{pf}
Assume that there are $n$ basis vertices in $\mathbb{T}$, denoted by $\mathbf{v}_i, i=0,1,\ldots,n-1$.
Let $b_{4i+j}(s,t), j=0,1,2,3$ be the four basis functions that associated with $\mathbf{v}_i$.
By the construction of $\mathcal{S}$, it follows that there are $4n$ basis functions in $\mathcal{S}$.

a) Let  $f(s,t)$ be a linear combination of all functions in $\mathcal{S}$,
$$
f(s,t)=\sum_{j=0}^{4n-1} c_j b_j(s,t).
$$
The linearly independence of $\mathcal{S}$ can be proven if
$$
f(s,t) \equiv 0 \ \Rightarrow \  c_j =0, j=0,\ldots, 4n-1.
$$

Define a linear operator $\mathfrak{L}(\cdot)$ as
\begin{equation}\label{eq-linear-lop}
\mathfrak{L}(f(s,t))= \big(f(s,t), f_s(s,t), f_t(s,t), f_{st}(s,t)\big).
\end{equation}
Substituting $\mathbf{v}_i$ into \eqref{eq-linear-lop} yields
$$
\mathfrak{L}(f(\mathbf{v}_i))= \sum_{j=0}^3 c_{4i+j}\mathfrak{L}(b_{4i+j}(\mathbf{v}_i))=  (c_{4i}, c_{4i+1}, c_{4i+2}, c_{4i+3} )\cdot \mathbf{B}(\mathbf{v}_i) = \mathbf{0},
$$
where $c_{4i+j}$ is the coefficient for $b_{4i+j}(s,t)$ and
\begin{align*}
  \mathbf{B}(s,t)=&(\mathfrak{L}b_{4i}(s,t), \mathfrak{L}b_{4i+1}(s,t), \mathfrak{L}b_{4i+2}(s,t), \mathfrak{L}b_{4i+3}(s,t)).&
\end{align*}
Suppose that the basis vertex $\mathbf{v}_i=(s_i,t_i)$ is of level $k, k=0,1,\ldots,N$.
Then its four neighbor cells at level $k$ form a rectangle with four vertices
$(s_i -3\Delta s_i, t_i-3\Delta t_i)$, $(s_i +3\Delta s_{i+1}, t_i-3\Delta t_i)$, $(s_i -3\Delta s_{i}, t_i+3\Delta t_{i+1})$, $(s_i +3\Delta s_{i+1}, t_i+3\Delta t_{i+1})$ (from left to right and bottom to top). Hence $\mathbf{B}(\mathbf{v}_i)$ can be expressed by
\begin{equation}\label{eq-B-matrix}
\mathbf{B}(\mathbf{v}_i)=\left(
  \begin{array}{cccc}
    (1-\lambda)(1-\mu) & \lambda(1-\mu) & \lambda \mu & \mu(1-\lambda) \\
    \lambda(1-\mu) & \alpha(1-\mu) & -\beta \mu & \alpha\beta \\
    (1-\lambda)\mu & \alpha \mu & \beta (1-\lambda) & \alpha\beta \\
    \lambda\mu & \alpha \mu & \beta \lambda & \alpha\beta \\
  \end{array}
\right),
\end{equation}
where $\alpha = \frac{1}{\Delta s_i +\Delta s_{i+1}}, \beta=\frac{1}{\Delta t_i +\Delta t_{i+1}}$, $\lambda=\alpha \Delta s_{i}$, $\mu=\beta \Delta t_{i}$. It is straightforward to check that $\mathbf{B}(\mathbf{v}_i)$ is invertible. Hence $(c_{4i}, c_{4i+1}, c_{4i+2}, c_{4i+3} ) =\mathbf{0}$.

Iterating through each basis vertex in $\mathbb{T}$, and finally we have $c_i=0, i=0,\ldots,4n-1$.  Hence the functions in $\mathcal{S}$ are linearly independent.

b) Since $b_i(s,t) \in \mathcal{S}(3,3,1,1,\mathbb{T})$, $i=0,\ldots, 4n-1$ and dim $\mathcal{S}(3,3,1,1,\mathbb{T})=4n$, it follows directly by a) that $\mathcal{S}$ spans the modified PHT-spline space $\mathcal{S}(3,3,1,1, \mathbb{T})$.

c) From level $k$ to level $k+1$, we can define
$$
f_1(s,t)= \sum_{j=0}^{4V_{k+1}-1} b_j^{k+1}(s,t),
f_2(s,t)= \sum_{j=0}^{d_k-1} \bar{b}_j^{k+1}(s,t),
$$
where $b_j^{k+1}(s,t), j=0,\ldots 4V_{k+1}-1$ are new basis functions associated with new basis vertices at level $k+1$, and  $\bar{b}_j^{k+1}(s,t), j=0,\ldots d_k-1$ are basis functions after modification according to Definition \ref{def-modify} above. Check the B\'ezier ordinates of $f_1(s,t)$ and $f_1(s,t)$ in every cell. For $f_1(s,t)$, the B\'ezier ordinates associated with new basis vertices are one, while all the others are zero; for $f_2(s,t)$, the B\'ezier ordinates associated with new basis vertices are zero, while all of the others are one. Thus
$
f_1(s,t)+f_2(s,t)\equiv 1.
$
\qed
\end{pf}
In addition, it is straightforward to check that the basis functions in $\mathcal{S}$ have other properties, such as nonnegativity, and local support.

\begin{rmk}
By setting certain B\'ezier ordinates to zero (as done in Definition \ref{def-modify}), we modify the basis functions, while THB-splines use similar techniques, which truncates the basis functions for splines with a more general setting \cite{giannelli-2012-thb}.
\end{rmk}


\begin{definition}
Let $\mathbb{T}$ be a modified hierarchical T-mesh, and $b_{j}(s,t)$, $j = 0, 1, \ldots, d$ be the basis functions constructed as Definition \ref{def-mpht-basis}. A modified PHT-spline surface over $\mathbb{T}$ is defined by
\begin{equation}\label{eq-spline-surface}
\mathbf{S}(s,t) = \sum_{j=0}^{d} \mathbf{C}_j b_j (s,t), \ (s,t) \in [0,1] \times [0,1],
\end{equation}
where $\mathbf{C}_j, j=0,1,\ldots,d$, are control points.
\end{definition}

Similar to PHT-spline surfaces in  \cite{deng2008}, modified PHT-spline surfaces have beneficial properties such as convex hull, affine invariant, and local support.

\section{Fitting open meshes}\label{section-fit-open-mesh}
Surface fitting is a fundamental task in computer aided geometric design and computer graphics that has been discussed in many papers (see \cite{chinvate1995} for a review of the literature). To show the potential of modified PHT-splines, we present a scheme to fit open mesh models based on modified PHT-spline surfaces in this section.

Given an open mesh $\mathcal{M}$ with vertices $\mathbf{P}_j$, $j = 0, 1, \ldots, M-1$. Let the corresponding parameter values of $\mathbf{P}_j$ be $(s_j, t_j)$, $j=0,1, \ldots,M-1$,  which are calculated from some parameterization of the mesh
(uniform weights are taken in barycentric mapping for the current paper).
The parameter domain is $[0,1]\times [0,1]$. For more details about mesh parameterization, we refer readers to a specific review on this topic in \cite{floater2004}.

\subsection{Outline of the fitting process}
The surface fitting scheme repeats the sequential steps 2, 3 and 4 until the fitting error in each cell is less than some tolerance $\varepsilon$.
\begin{enumerate}
  \item Initialize a tensor product mesh $\mathbb{T}_0$. Set level $k=0$.
  \item Keep unchanged the control points associated with the old basis functions. Compute the control points for the new basis functions on the $k$th level mesh $\mathbb{T}_k$ to obtain a modified PHT-spline surface $\mathbf{S}_k(s,t)$ (in the beginning, every basis function is new). See the following subsection \ref{sub-compute-ctrl} for details.
  \item Find the cells of level $k$ whose fitting errors are greater than $\varepsilon$, and denote them as $\Theta_k$.  The fitting error over the cell $\theta$ is defined to be
      $\max_{(s_j,t_j)\in\theta} \|\mathbf{P}_j - \mathbf{S}_k(s_j,t_j)\|$.
  \item Label each cell in $\Theta_k$ with marks based on discrete curvature information (see the following Subsection \ref{sub-label-cell} for details). Subdivide these cells as Algorithm 1. Set $k=k+1$.
\end{enumerate}

\subsection{Compute control points}
\label{sub-compute-ctrl}
To get a modified PHT-spline surface that fit the given open mesh $\mathcal{M}$, one needs to determine the control points in \eqref{eq-spline-surface}.
As a standard way to evaluate the control points, we may solve the following least square optimization problem
$$
\min \sum_{i=0}^{M-1}\|\mathbf{S}(s_i, t_i) - \mathbf{P}_i \|^2
= \min_{\mathbf{C}_j} \sum_{i=0}^{M-1}\| \sum_{j=0}^{d} \mathbf{C}_j b_j (s_i,t_i) - \mathbf{P}_i\|^2,
$$
which needs to solve $\mathbf{C}_j$ globally.

Similar to \cite{deng2008,speleers-2016-thb-gm}, an efficient method based on local geometric information can be provided as follows.

Based on the proof of Theorem \ref{thm-6}, for each basis vertex $\mathbf{v}_i$ and its associated four basis functions $b_{4i+k}(s,t), k=0, 1, 2, 3$,
\begin{equation}\label{eq-ctrlpts-eva}
\mathfrak{L}(\mathbf{S}(\mathbf{v}_i)) = \sum_{j=0}^{d} \mathbf{C}_j \mathfrak{L} b_j (\mathbf{v}_i) = \sum_{k=0}^3 \mathbf{C}_{4i+k} \mathfrak{L}b_{4i+k}(\mathbf{v}_i) = \mathbf{C}\cdot \mathbf{B}(\mathbf{v}_i),
\end{equation}
where $\mathbf{C}=(\mathbf{C}_{4i}, \mathbf{C}_{4i+1}, \mathbf{C}_{4i+2}, \mathbf{C}_{4i+3})$ is a $3\times 4$ matrix and $\mathbf{B}(\mathbf{v}_i)$ is defined as \eqref{eq-B-matrix}.
On the other hand, by fitting the set of points around $\mathbf{v}_i$ with a quadratic surface, $\mathfrak{L}(\mathbf{S}(\mathbf{v}_i))$ can be approximated. Hence by \eqref{eq-ctrlpts-eva}, the control points $\mathbf{C}$ can be estimated by
$
\mathcal{L}(\mathbf{S}(\mathbf{v}_i)) \cdot \mathbf{B}^{-1}.
$

\subsection{Label the cells}\label{sub-label-cell}
To capture the anisotropic information of the fitting surface, discrete curvature information of $\mathbf{S}_k(s,t)$ over those cells is evaluated and used as follows.

For any cell $\theta \in \Theta_k$, choose $l$ parametric points $(\overline{s}_{j}, \overline{t}_{j}) \in \theta$, $j=1,\ldots,l$.
Curvatures along the $s$ and $t$ directions at these points $\mathbf{S}_k(\overline{s}_{j}, \overline{t}_{j})$
are
 $$
 {\kappa_s}_j = \frac{\|\frac{\partial \mathbf{S}_k}{\partial s} \times \frac{\partial^2 \mathbf{S}_k}{\partial s^2}\|}{\| \frac{\partial \mathbf{S}_k}{\partial s} \|^3} \bigg|_{(\overline{s}_{j}, \overline{t}_{j})}, \ \ \
 {\kappa_t}_j = \frac{\|\frac{\partial \mathbf{S}_k}{\partial t} \times \frac{\partial^2 \mathbf{S}_k}{\partial t^2}\|}{\| \frac{\partial \mathbf{S}_k}{\partial t} \|^3} \bigg|_{(\overline{s}_j, \overline{t}_j)}, \ j=1,\ldots, l.
 $$
 Set $K_s = \frac{1}{l}\sum_{j=1}^{l} {\kappa_{s}}_j$ and $K_t = \frac{1}{l} \sum_{j=1}^{l} {\kappa_{t}}_j$.

Assume that $K_t \neq 0$, then the ratio $\rho = K_s/K_t$ is chosen to characterize the anisotropic feature. When $\rho >\delta$ or $\rho <1/\delta$, the fitting surface is considered to change sharply along one direction but stay flat along the other direction, which is assumed to possess the anisotropic feature over the cell $\theta$. Hence when $\rho >\delta$,  label the cell $\theta$ with `V'.  Similarly, label the cell $\theta$ with `H' when $\rho <1/\delta$. For the other cases, label the cell $\theta$ with `C'. Note that the choice of the threshold $\delta$ may influence the label of anisotropic information on each cell. Therefore, $\delta$ can be assigned depending on the problem being addressed.

\subsection{Numerical Examples}
\begin{figure}[htbp!]
\centering
\includegraphics[width=3cm,height=3cm]{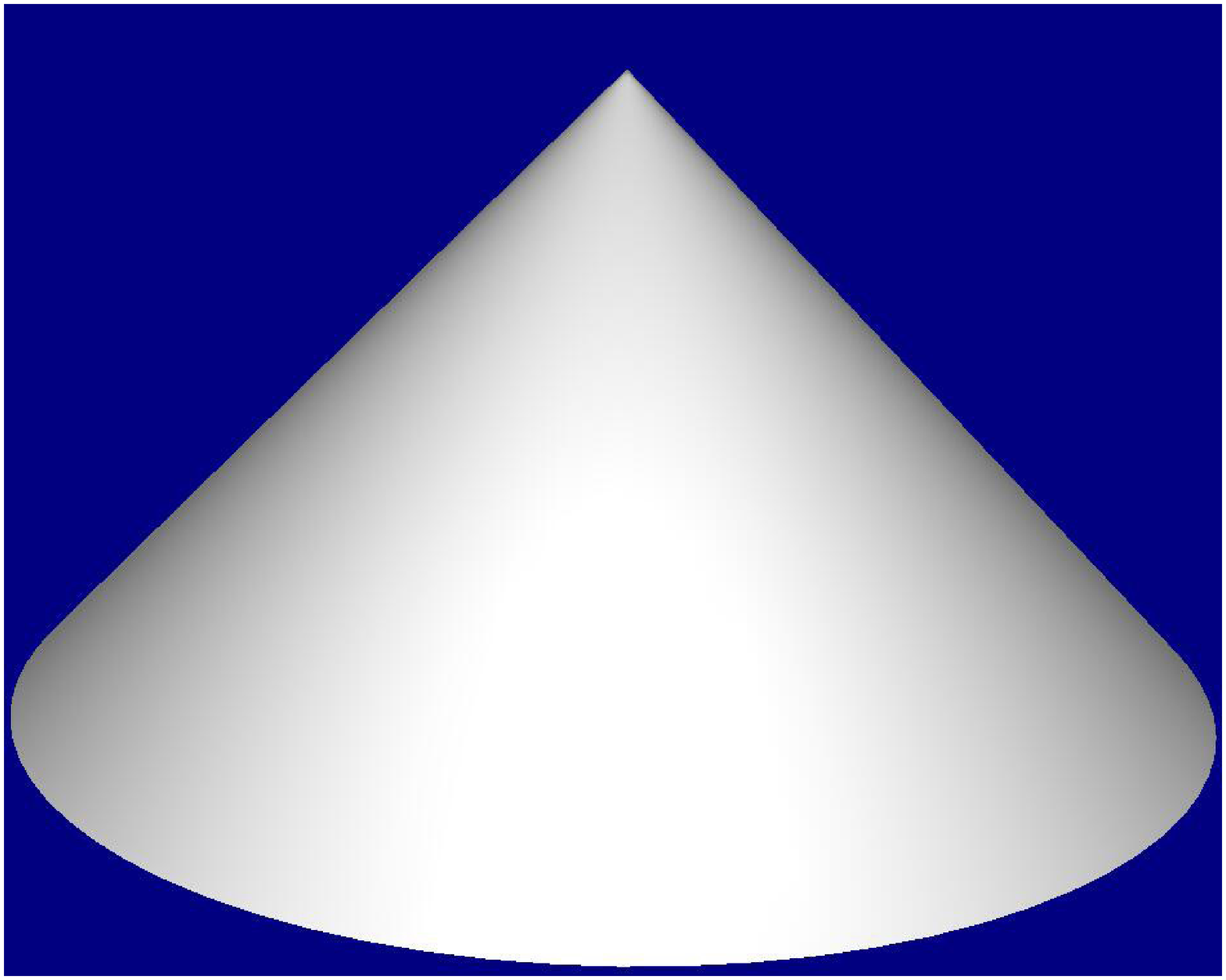}\hfill
\includegraphics[width=3cm,height=3cm]{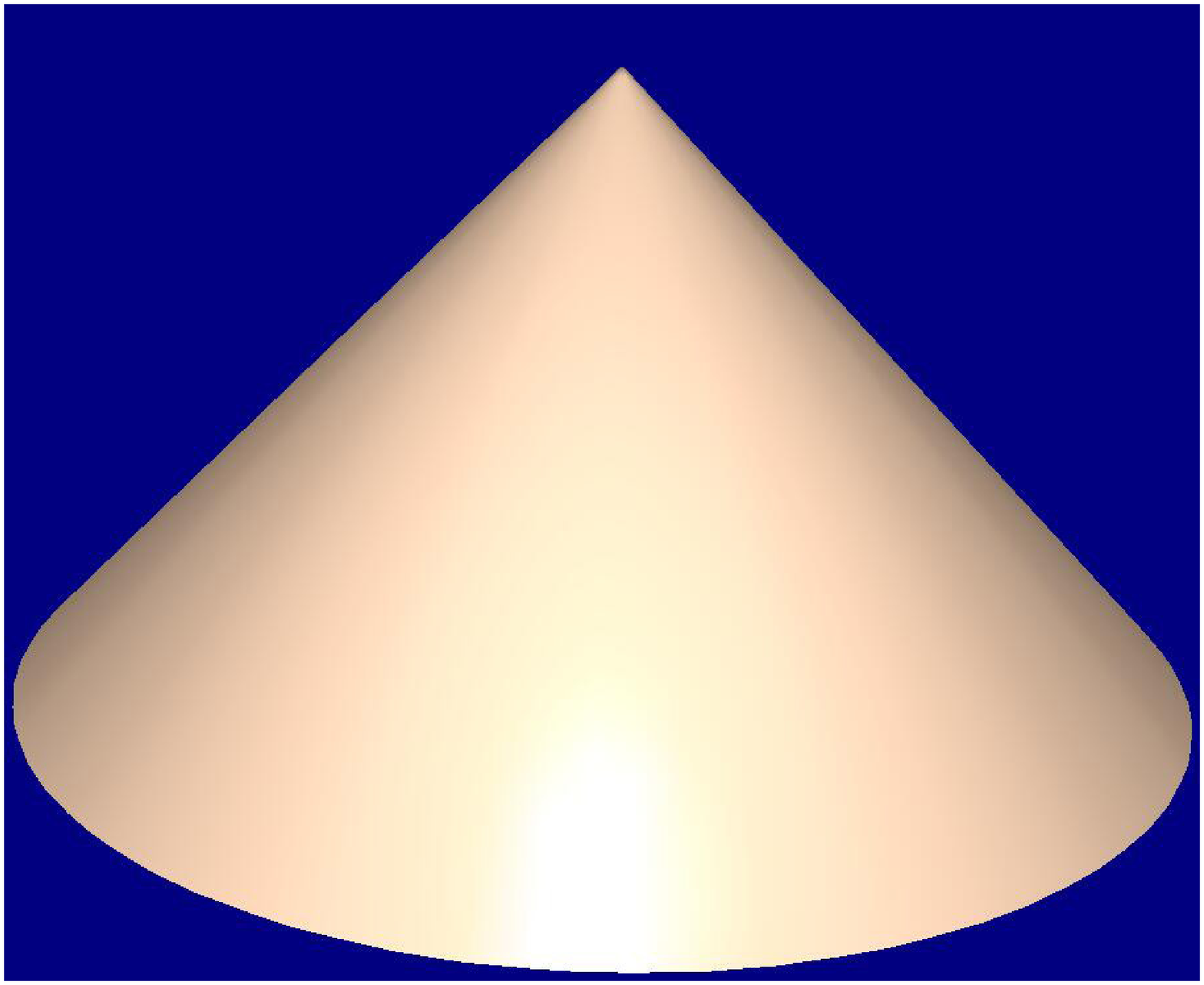}\hfill
\includegraphics[width=3cm,height=3cm]{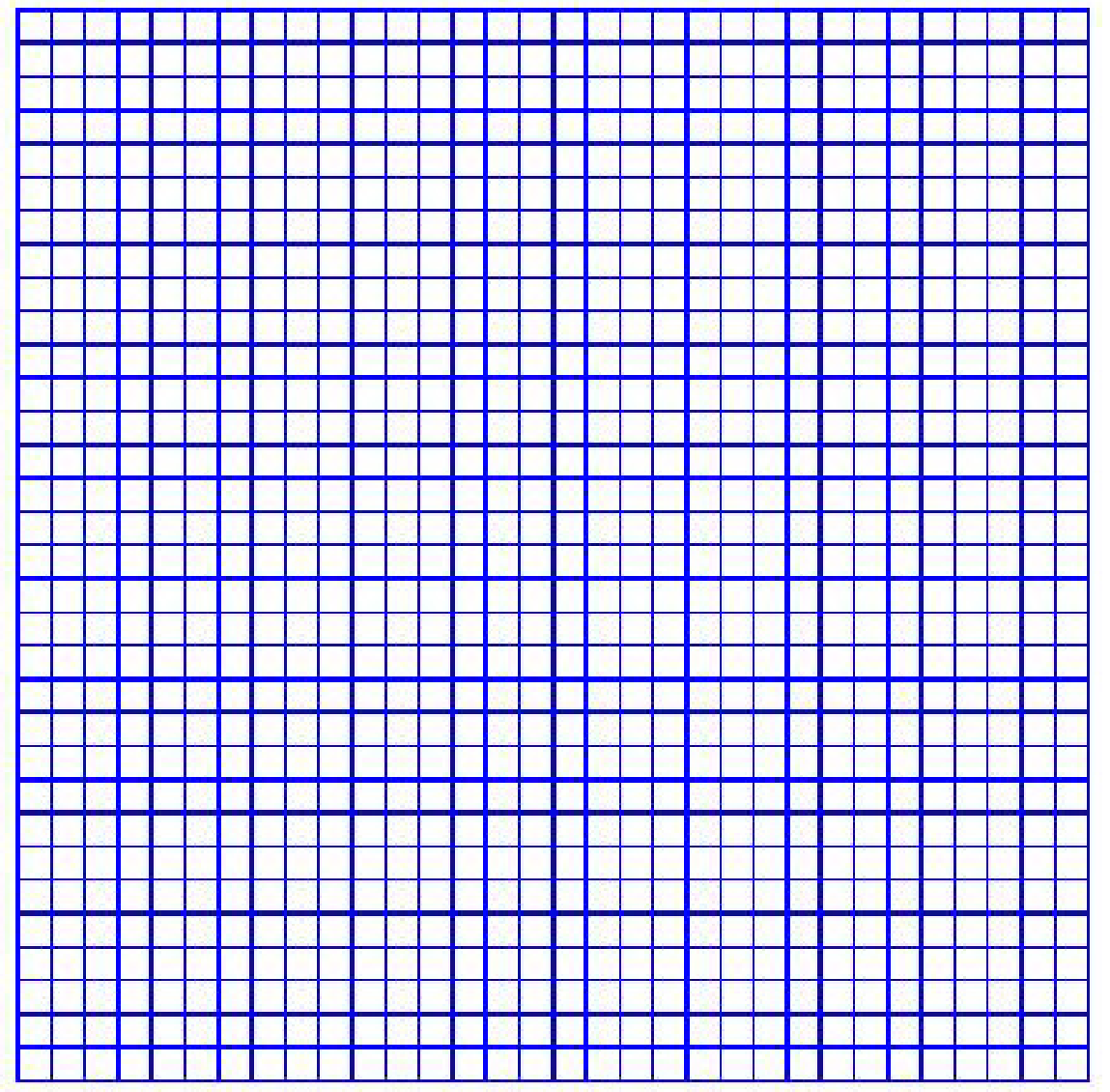}\hfill
\includegraphics[width=3cm,height=3cm]{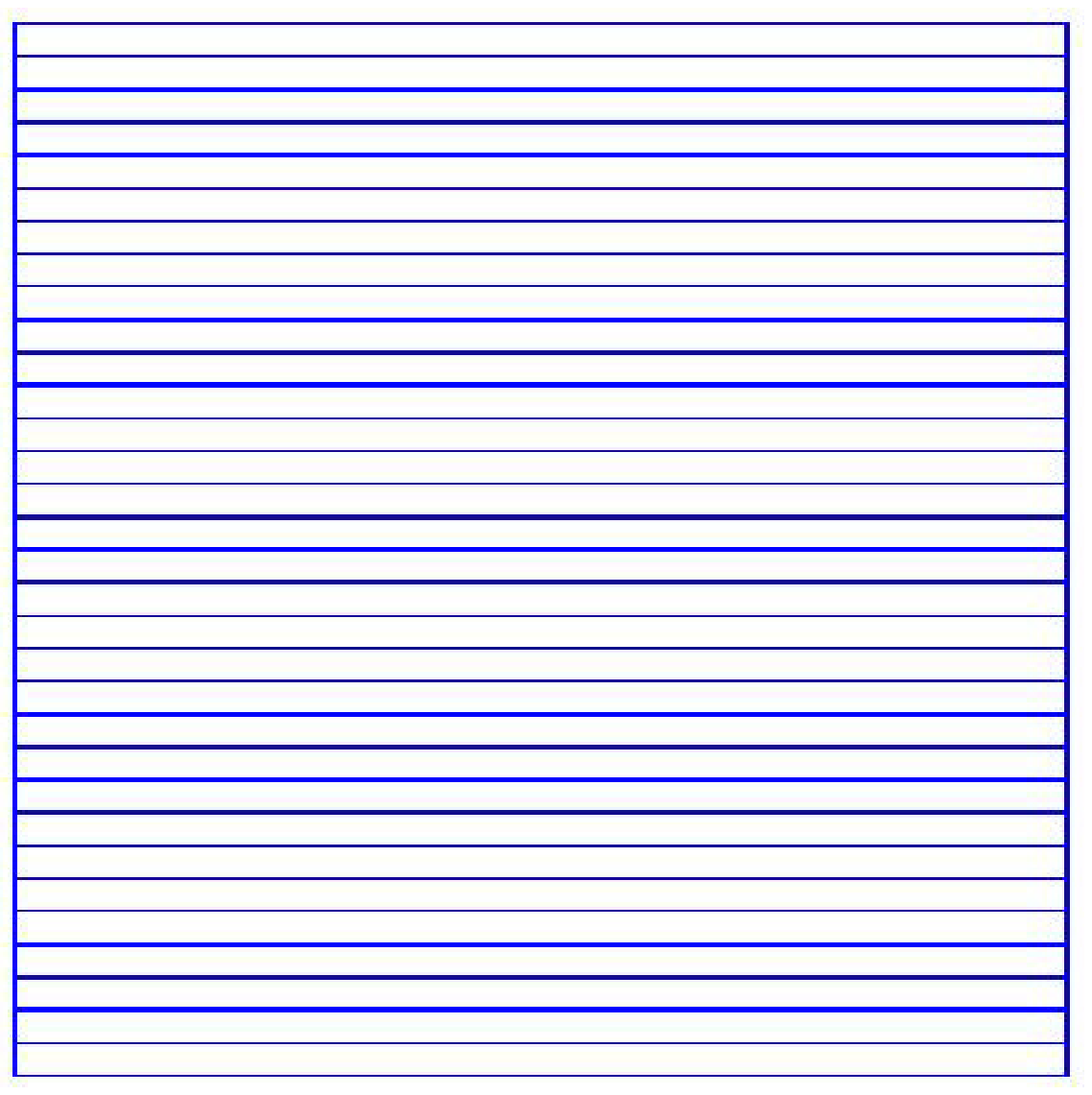}\vspace{0.5cm}\\
\includegraphics[width=3cm,height=3cm]{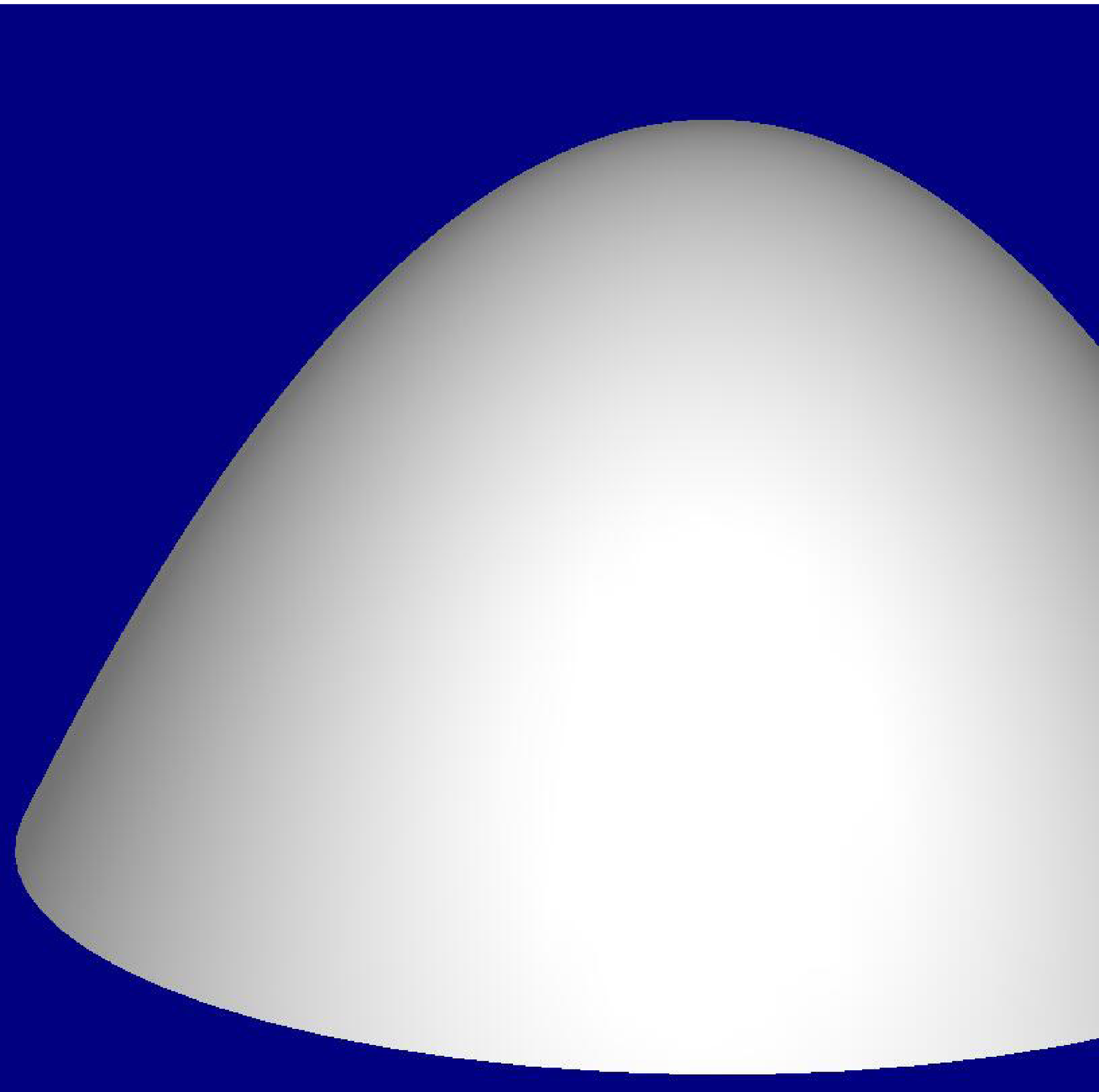}\hfill
\includegraphics[width=3cm,height=3cm]{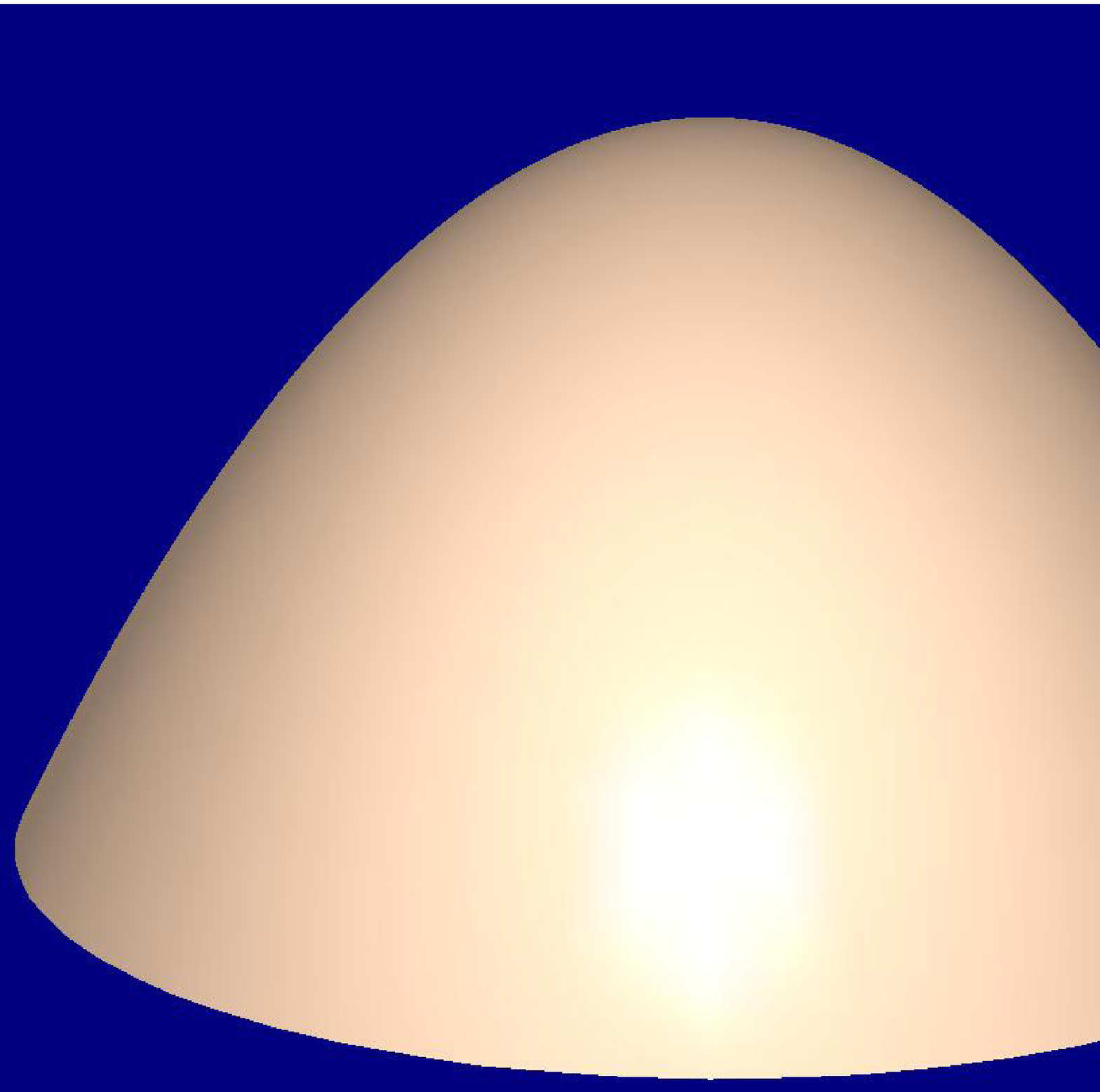}\hfill
\includegraphics[width=3cm,height=3cm]{cone1.eps}\hfill
\includegraphics[width=3cm,height=3cm]{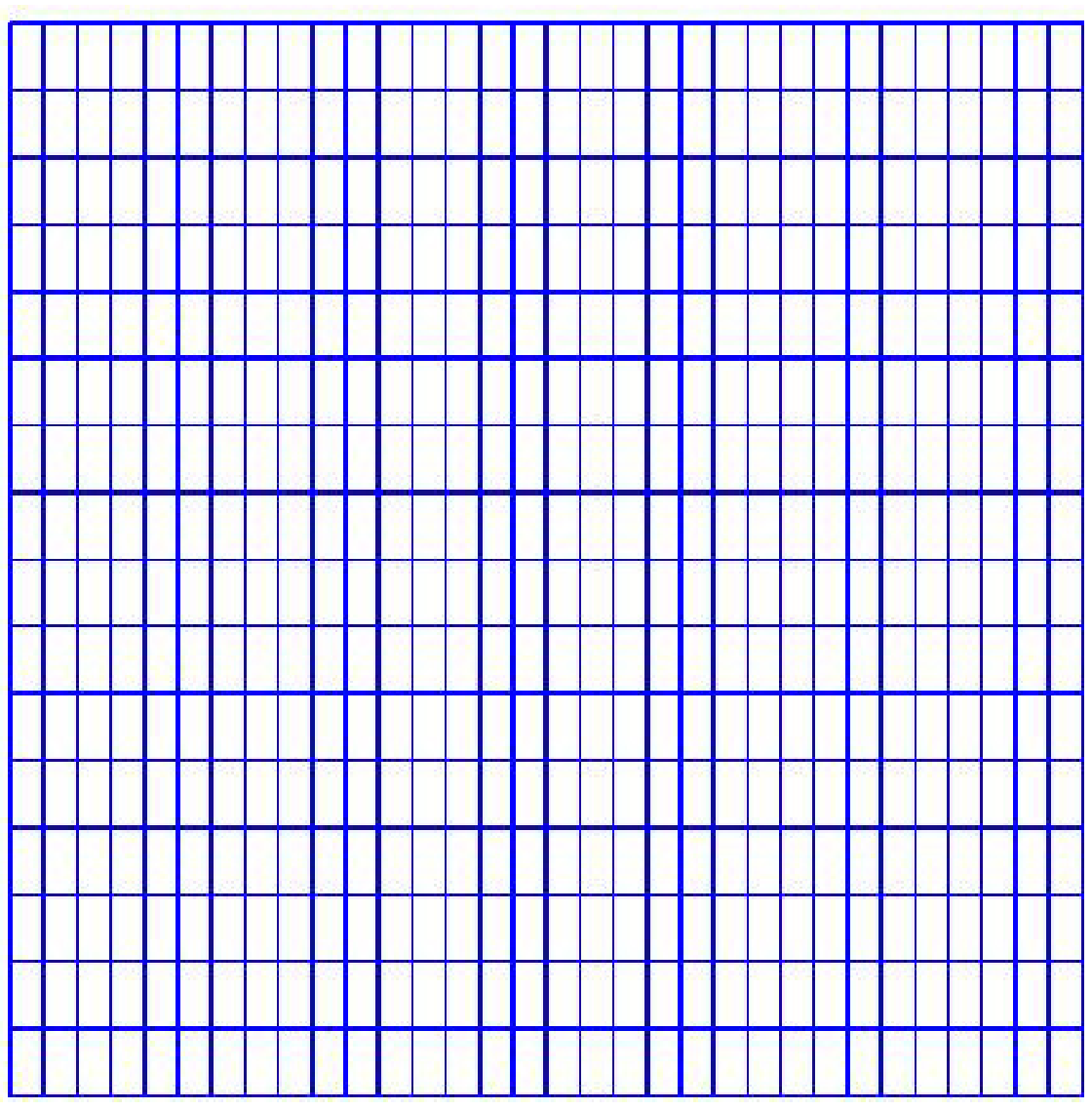}\vspace{0.5cm}\\
\includegraphics[width=3cm,height=3cm]{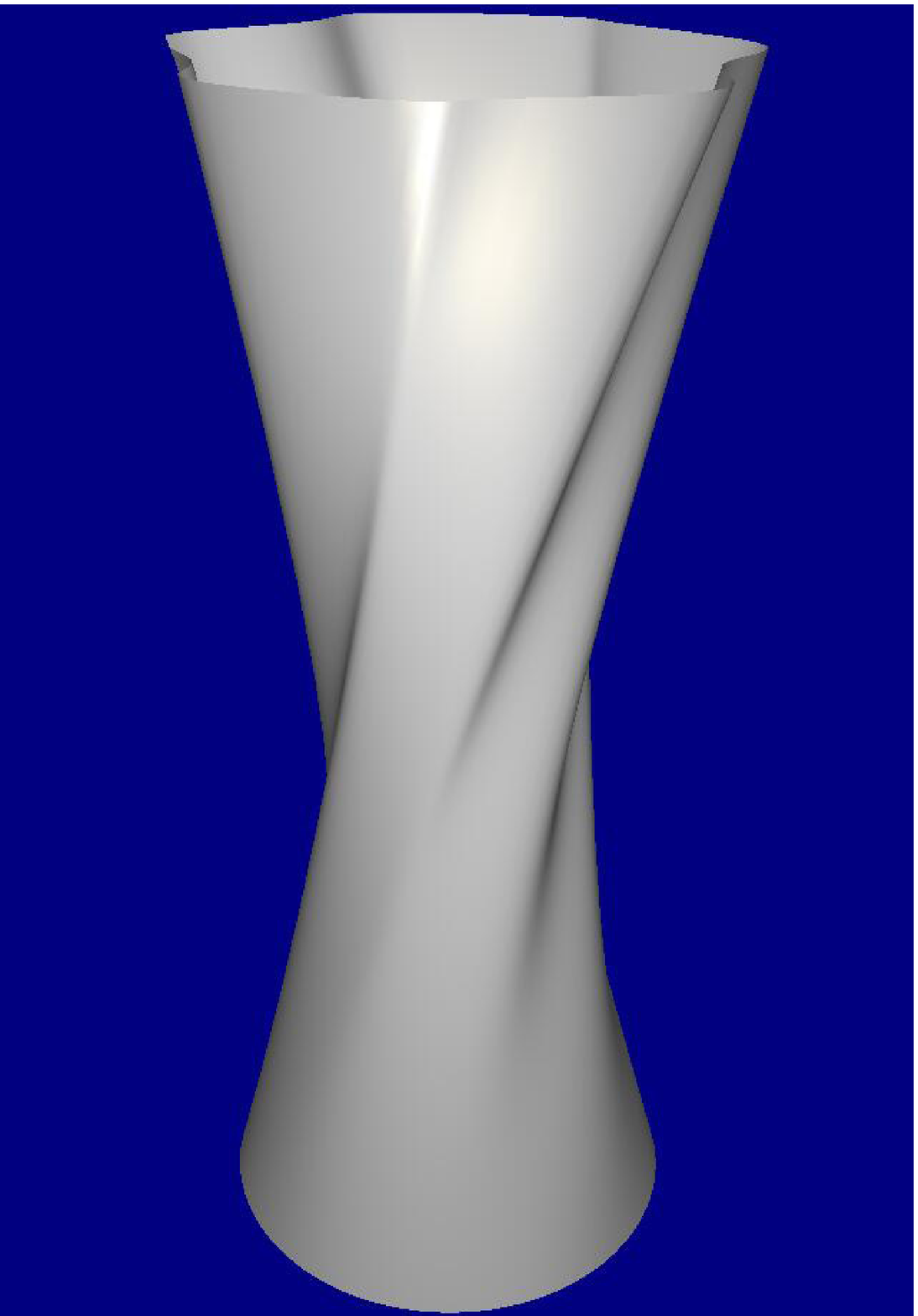}\hfill
\includegraphics[width=3cm,height=3cm]{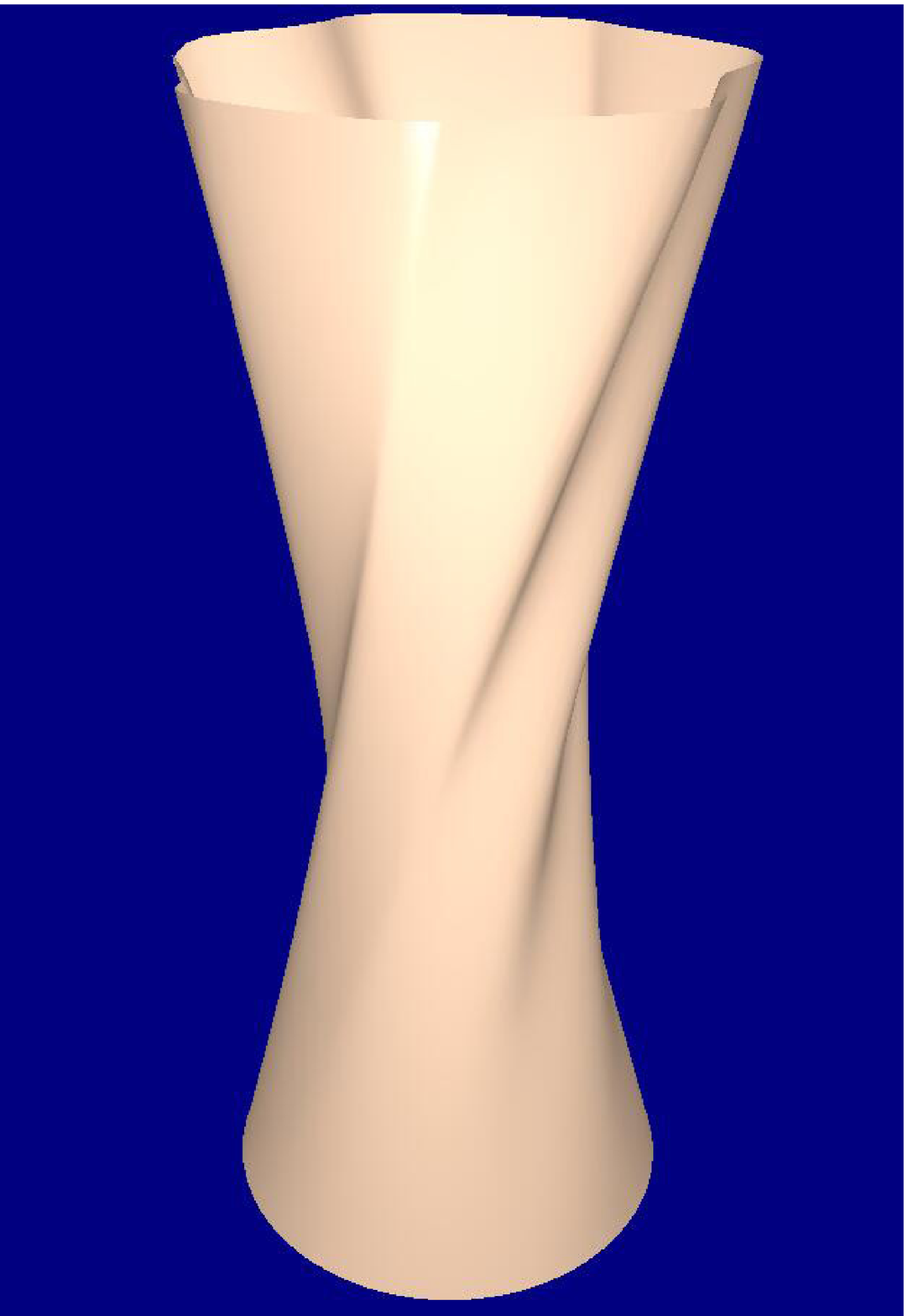}\hfill
\includegraphics[width=3cm,height=3cm]{cone1.eps}\hfill
\includegraphics[width=3cm,height=3cm]{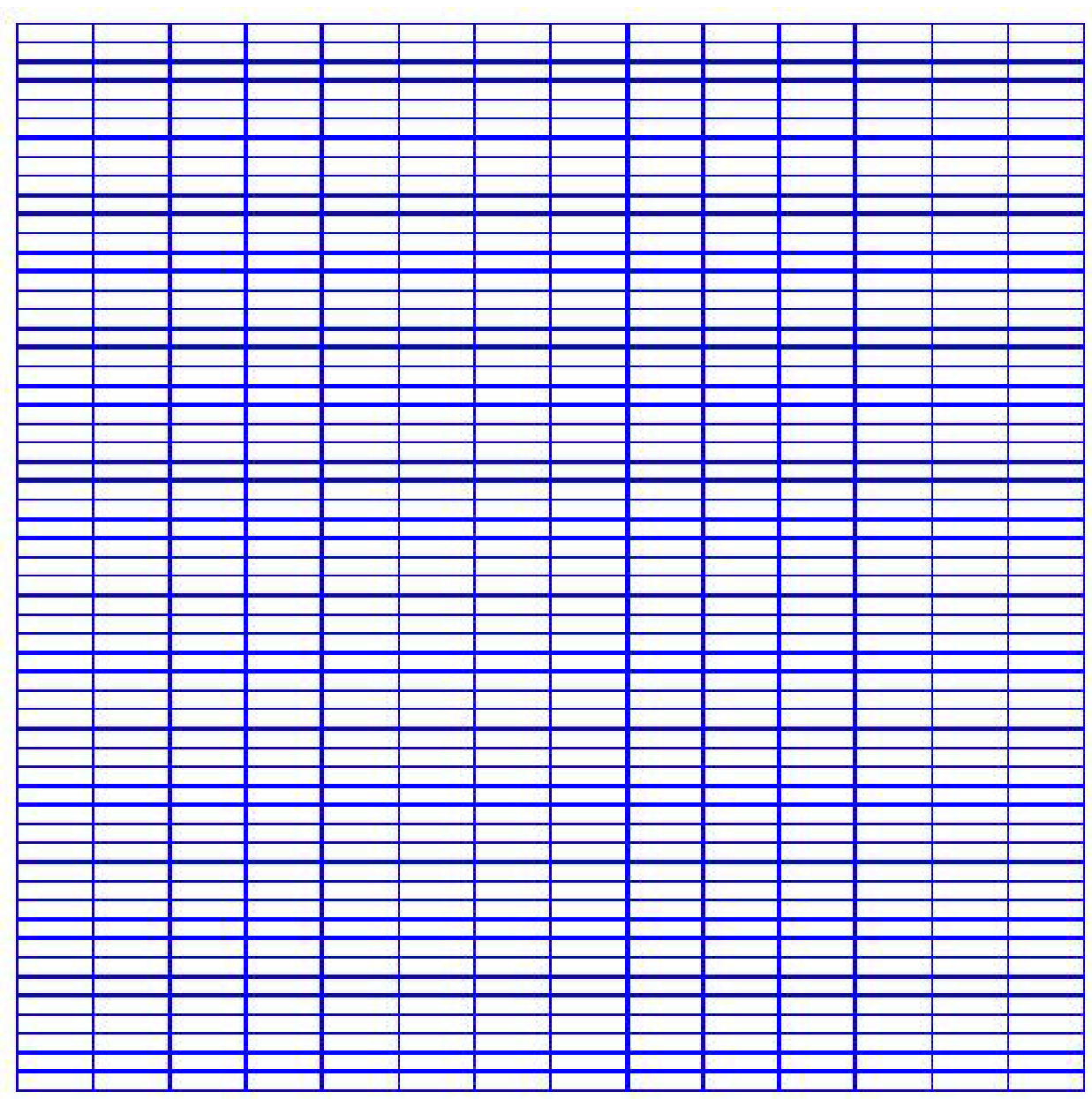}\vspace{0.2cm}\\\hspace{-0.25cm}
\subfloat[Original model]{
\includegraphics[width=3cm,height=3cm]{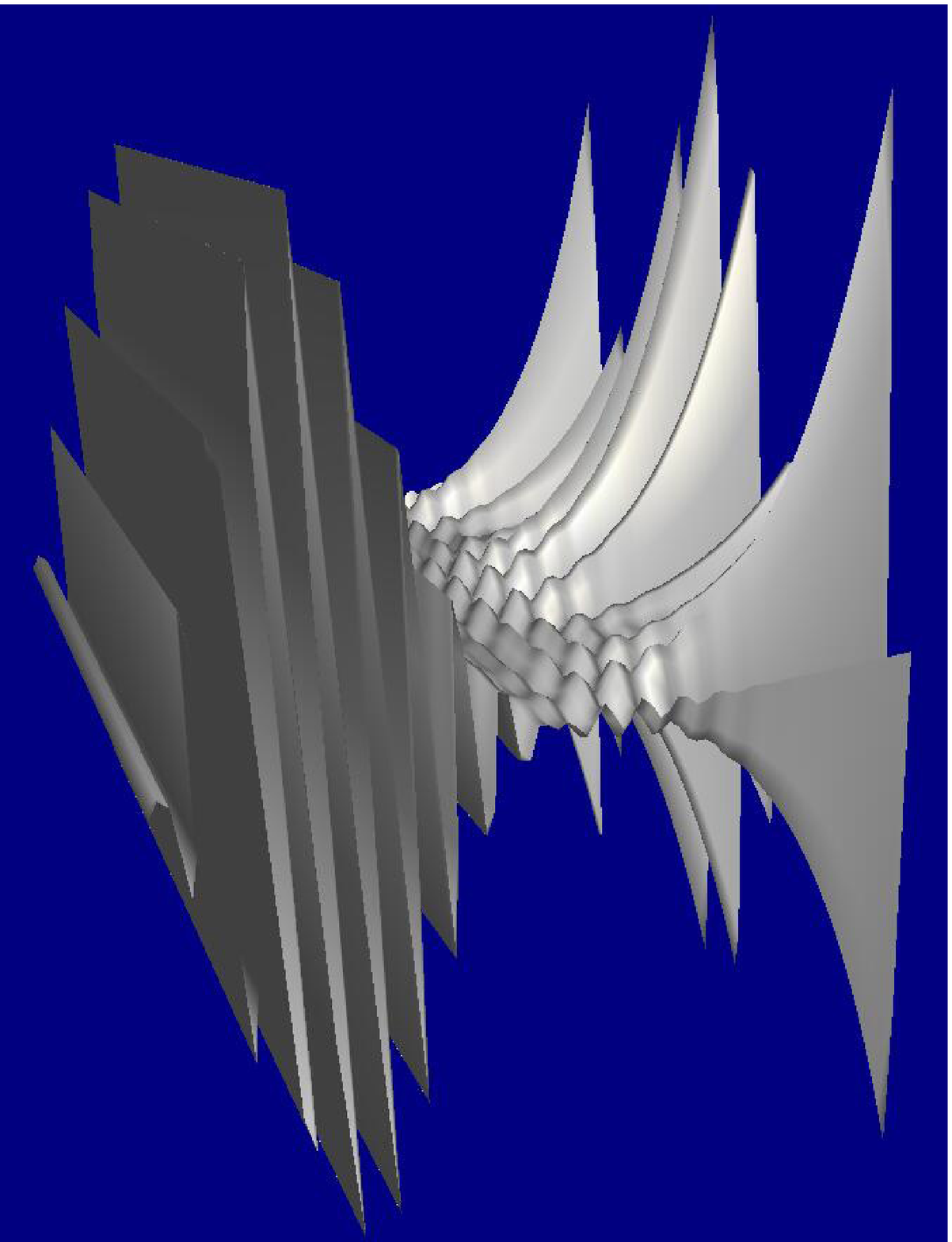}}\hfill
\subfloat[Fitting results]{
\includegraphics[width=3cm,height=3cm]{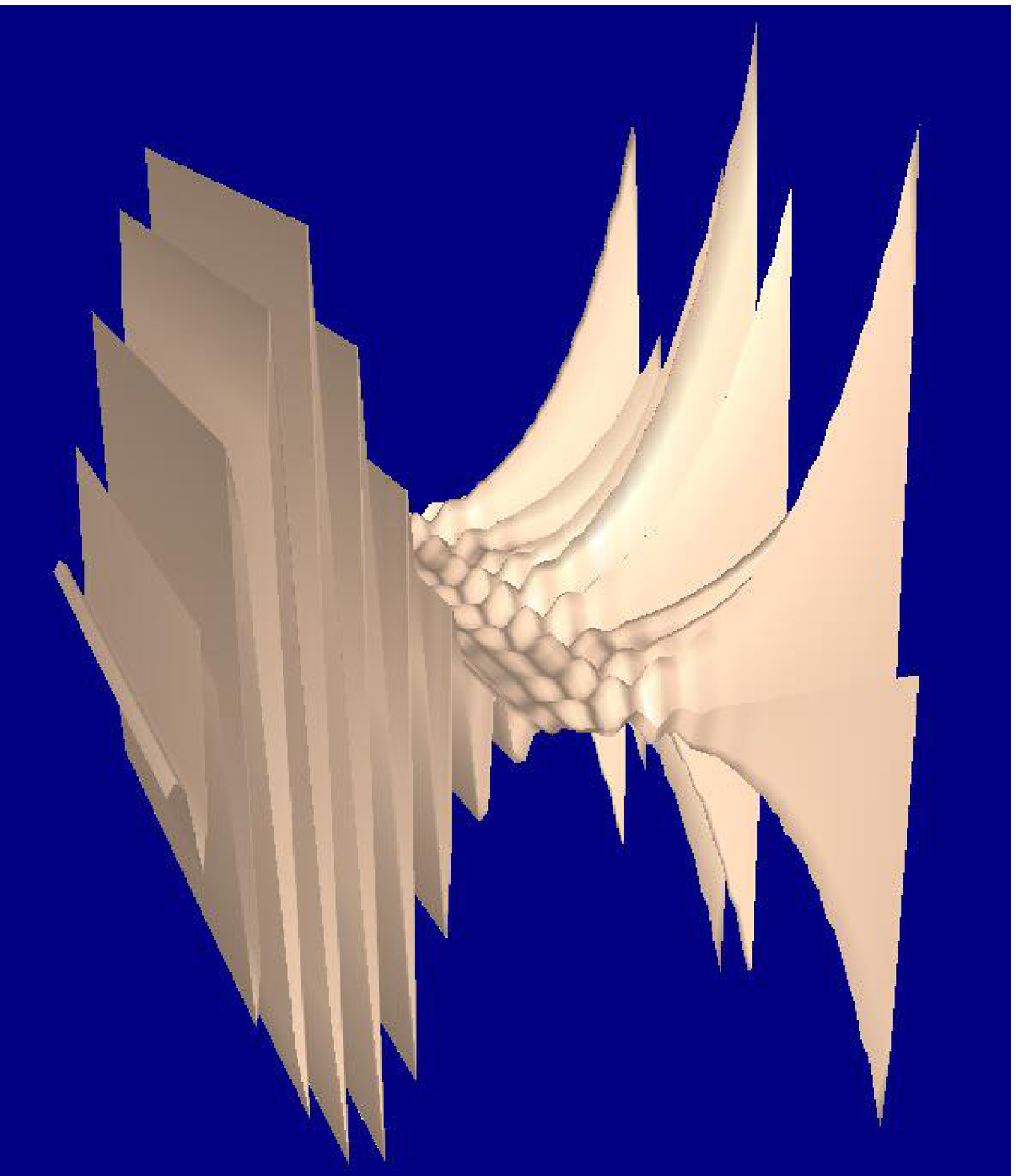}}\hfill
\subfloat[by PHT-splines]{
\includegraphics[width=3cm,height=3cm]{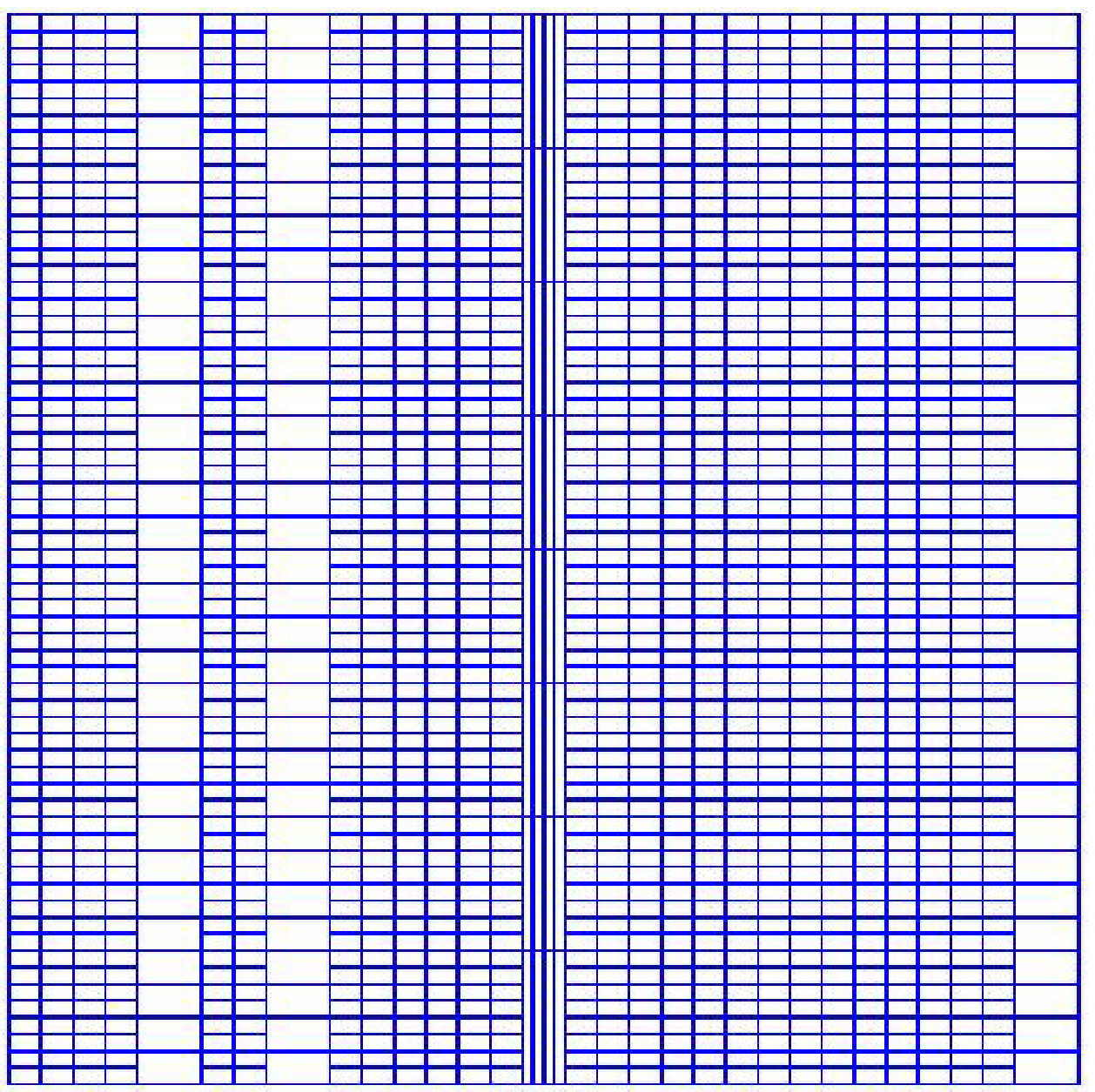}}\hfill
\subfloat[by modified PHT-splines]{
\includegraphics[width=3cm,height=3cm]{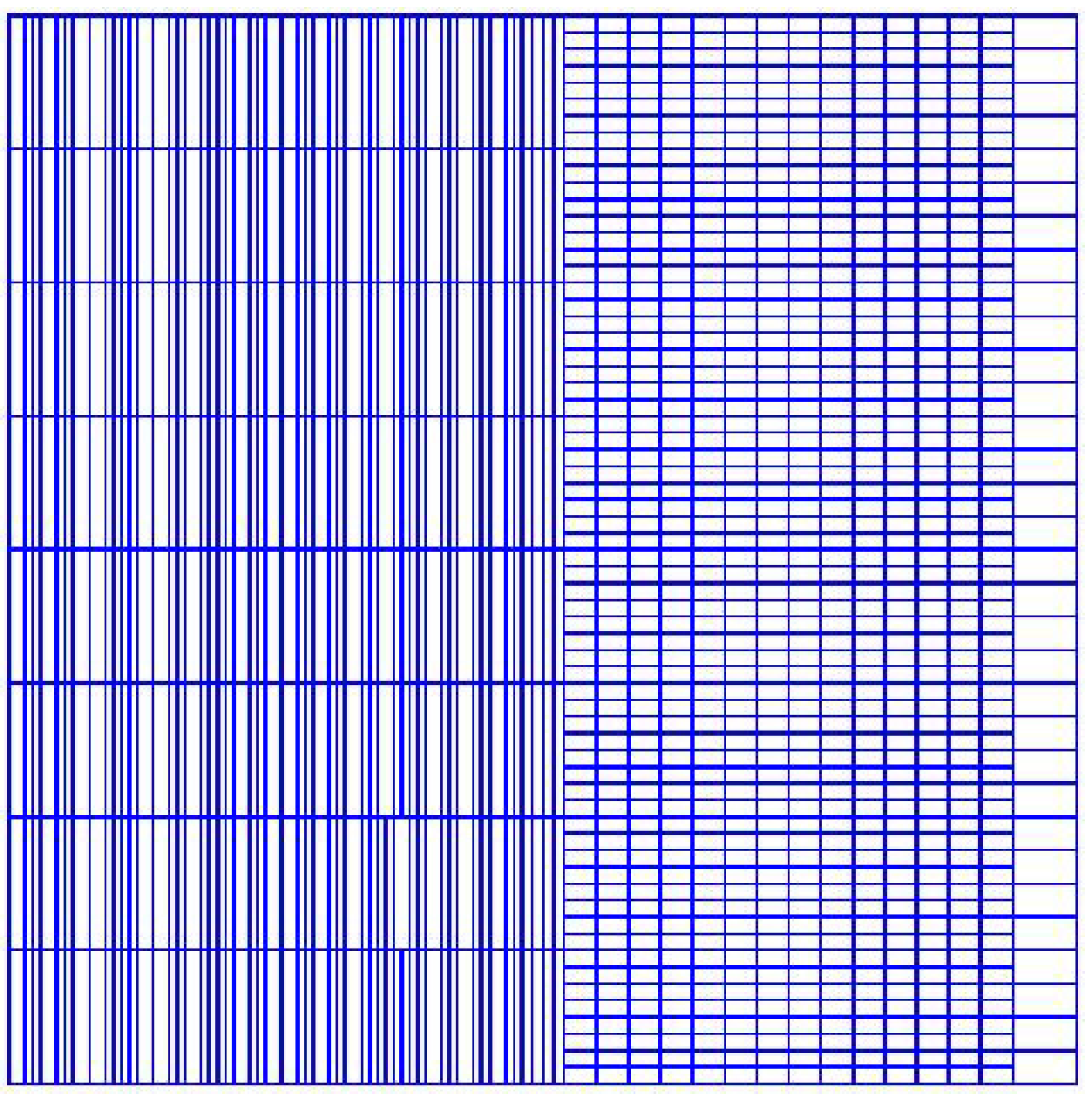}}
\caption{ Four examples of fitting open meshes by PHT-splines and modified PHT-splines.}
\label{fig:fitting}
\end{figure}

\reffig{fig:fitting} illustrates four examples, in which fitting open mesh models with modified PHT-splines are shown. The first two examples are open meshes generated by sampling points on a conic surface and a paraboloid, respectively. The fittings of these two open meshes work as prototypes to verify the efficiency of our modified PHT-splines when dealing with open meshes with anisotropic features. The third model is a practical model. The fourth model is referred to as Example 6 in \cite{engleitner-2017-aniso}.  We sampled $101\times 101$ data points on a uniform grid from the function
\begin{align*}
f(u,v) &= 0.1(B_{0,7}(u)(\sin(120u)\sin(2\pi u))\\
&+ B_{1,7}(u)(2\sin(120u)\sin(2\pi u))\\
&+ B_{7,7}(u)(2-2(1+0.4\sin(60v))|\cos(2\pi v)|)) ,
\end{align*}
with $B_{i,7}(u)$ being the $i$-th Bernstein polynomial of degree $7$.

For all four examples, the initial T-meshes are set to $2\times 2$ uniform tensor-product meshes for the first two examples and $5\times 5$ tensor-product meshes for the latter two examples on the square domain $[0,1]\times[0,1]$.
The tolerance of the fitting error is $\epsilon=0.1\%$, which refers to the size of the corresponding model's bounding box. The threshold $\delta =0.5$ is chosen to characterize anisotropic features for all of these four examples. The statistical data for comparison are listed in Table \ref{table-compare}.

\begin{table}[htbp!]
\caption{Comparison of results}
\centering
\begin{tabular}{c|ccc|ccc}
\hline
\   \   & \multicolumn{3}{c|}{original method}&\multicolumn{3}{|c}{new method} \\
\hline
Level 0 & \#cp & max error &average error& \#cp & max error&average error \\
\hline
Level 1    & 4356 & 0.0174641 &0.00161962& 520 & 0.00486232 &0.00160813\\
\hline
Level 2& 4356 & 0.0174641 & 0.00162178 & 2244 & 0.0174388 & 0.00219253\\
\hline
Level 3&4356&0.0222689& 0.00222854&3420&0.0209686& 0.00248527\\
\hline
Level 4& 6952 & 0.0546684&0.01776762 &    5552&0.0336198&0.0193418\\
\hline
\end{tabular}
\label{table-compare}
\end{table}

\section{IGA based on modified PHT-splines}\label{section-iga}
In this section, we attempt to use our modified PHT-splines in an isogeometric method to solve elliptic partial differential equations.

\subsection{Model problem}
Suppose the model problem is an elliptic partial  differential equation defined as
\begin{align}\label{eq-model}
&-\Delta u = f \mbox{ in } \Omega,&\nonumber\\
&u = 0 \mbox{ on } \Gamma_D,&\\
&\frac{\partial u}{\partial \mathbf{n}} = g_N \mbox{ on } \Gamma_N\nonumber&
\end{align}
where $\Omega \in \mathbb{R}^2$ is a connected, bounded domain with a Lipschitz-continuous boundary $\Gamma = \Gamma_D \bigcup \Gamma_N$, $\Gamma_D \bigcap \Gamma_N = \emptyset$, $\mathbf{n}$ is the outward unit normal to $\Gamma_N$,  $\Gamma_D$ is assumed to be closed relative to $\Gamma$ and has a positive length, while $f$ and $h$ are square-integrable on $\Omega$ and $\Gamma_N$, respectively.

\subsection{Discretization}
The framework of isogeometric analysis based on modified PHT-splines is as follows. For more details about isogeometric analysis based on splines, see \cite{Thann-2013-pht,pingwang2011} and references therein.

Suppose the parametrization $\mathbf{G}$ of the physical domain $\Omega$ is defined by
\begin{align*}
\mathbf{G}: & \Omega_0  \rightarrow \mathbb{R}^2&\nonumber\\
&(s,t)\mapsto (x,y)= \sum_{i=1}^{d}
\mathbf{P}_i b_i(s,t),
\end{align*}
where $\Omega_0 =[0,1]\times [0,1]$, $\mathbf{P}_i \in \mathbb{R}^2$, $b_i(s,t)$ is a modified PHT-spline function,
$d$ is the number of basis functions.

Let $V(\Omega) = H^1(\Omega)$ be the underlying Hilbert space for both the space of test functions and the solution space. The function space $H^1(\Omega)$ is defined by
$$
H^1(\Omega) =\{v\in L^2(\Omega): |\nabla v|\in L^2(\Omega)\}
$$
The weak form solution of problem \eqref{eq-model} is to find
$u \in V= \{ v \in H^1(\Omega): v|_{\Gamma_D} =0\}$
such that
$$
a(u,v)= \langle f, v\rangle, \forall v\in V
$$
where $a(,)$ is a bilinear form and $\langle, \rangle$ is a linear functional defined by
\begin{align*}
&a(u,v) = \int_{\Omega} \nabla u \cdot \nabla v d\Omega, &\\
&\langle f,v\rangle = \int_{\Omega} fv d\Omega +\int_{\Gamma_N} g_N v d\Gamma&
\end{align*}

Based on Galerkin's principle, finite dimensional function space $V_h$ is set up to solve the following problem:

Find $u^h \in V^h$ such that
\begin{equation}\label{eq-weak-form-approx}
a(u^h, v^h) = \langle f, v^h\rangle, \forall v^h\in V^h,
\end{equation}
where $V^h \subset V$ is
$$
V^h= \mbox{span}\{\psi_i(x,y) = b_i \circ \mathbf{G}^{-1}, \psi_i(x,y)|_{\Gamma_D}=0, i=1,2,\ldots, n\}.
$$
The approximation solution $u^h$ can be written as
$$
u^h(x,y)=\sum_{i=1}^{n} c_i \psi_i,
$$
where $c_i, i=1,\ldots, n$ are coefficients need to be determined.
Define the stiffness matrix $\mathbf{A}$ by
$$
\mathbf{A}= (a_{ij})_{i,j=1}^{n}, \mbox{ with } a_{ij} = a(\psi_i,\psi_j),
$$
and the load vector $\mathbf{F}$ by
$$
\mathbf{F}= (F_i)_{i=1}^{n}, \mbox{ with } F_{i} = \langle F, \psi_i \rangle.
$$
Thus problem \eqref{eq-weak-form-approx} is equivalent to the following linear system
$$
\mathbf{A} \mathbf{c} = \mathbf{F},
$$
where $\mathbf{c}= (c_1, c_2,\ldots, c_n)$ is the coefficient vector.

\subsection{Solving in the adaptive process with modified PHT-splines}\label{sub-iga-label}
Consider the flexibility of the T-meshes and local refinement algorithm for modified PHT-splines, the adaptive procedure consists of the following successive loops
$$
\mbox{Solve} \rightarrow \mbox{Estimate \& Mark} \rightarrow \mbox{Refine}.
$$
The essential part of the loops is the estimate and mark step. Follow the well-developed way in IGA \cite{kleiss-2012-iga,pingwang2011}, the posteriori error is considered here. It should be noted that error estimators for anisotropic refinement have been discussed in the literature on finite element methods \cite{Siebert-1996-error}. The posteriori error on a cell is
\begin{equation}\label{eq-post-error}
\eta_\theta^2 = h_\theta^2 \|\Delta u_h + f\|_{L^2(\theta)}^2+\|g_N-\frac{\partial u_h}{\partial \textbf{n}}\|_{L^2(\partial\theta)}^2 h_\theta,
\end{equation}
where $h_{\theta}$ is the diameter of cell $\theta$ and $\Delta = \frac{\partial^2}{\partial x^2}+ \frac{\partial^2}{\partial y^2}$.
The posteriori error on the modified hierarchical T-mesh $\mathbb{T}$ is the sum of the posteriori errors over all the cells, that is,
$$
\eta_\tau = \sqrt{\sum_{\theta\in \mathbb{T}}\eta_{\theta}^2}
$$

The strategy for the estimate and mark, and refine steps in the loops is as follows.
\begin{enumerate}[Step 1.]
  \item Find the cells of the current level that need to be subdivided, which are indicated by the the posteriori error estimator \eqref{eq-post-error}. Here we use a prescribed threshold to decides the cells.
  \item For each cell that to be subdivided, choose $l$ parametric points $(\overline{s}_j, \overline{t}_j)$, $j=1,\ldots, l$ in the cell and evaluate
  $$
  \hat{\kappa}_{sj} = \frac{\|\frac{\partial u^h(G(s,t))}{\partial s} \times \frac{\partial^2 u^h(G(s,t))}{\partial s^2}\|}{\| \frac{\partial u^h(G(s,t))}{\partial s} \|^3} \bigg|_{(\overline{s}_{j}, \overline{t}_{j})}, \
  \hat{\kappa}_{tj} = \frac{\|\frac{\partial u^h(G(s,t))}{\partial t} \times \frac{\partial^2 u^h(G(s,t))}{\partial t^2}\|}{\| \frac{\partial u^h(G(s,t))}{\partial t} \|^3} \bigg|_{(\overline{s}_{j}, \overline{t}_{j})},
  $$
  Set $\hat{K}_s =\sum_{j=1}^{l}{\hat{k}_{sj}}/l, \hat{K}_t =\sum_{j=1}^{l}{\hat{k}_{tj}}/l$. If $\hat{K}_t\neq 0$, the ratio  $\hat{K}_s/\hat{K}_t$ is adopted to describe the anisotropic features of this cell. If $\hat{K}_s/\hat{K}_t > \delta_0$ or $\hat{K}_s/\hat{K}_t <\delta_1$, where $\delta_0>\delta_1$ be two prescribed values, the solution $u^h(G(s,t))$ is assumed to change sharply along one direction but stay flat along the other direction. In addition, we can treat the case $\hat{K}_t= 0$ similarly.
  \item Use Algorithm 1 to refine the mesh.
\end{enumerate}

One of the important steps of constructing modified PHT-splines is to label cells in T-meshes, which indicates the subdivision type the cells tend to choose.
Similar to the method we used for fitting open mesh in subsection \ref{sub-label-cell}, for each cell that needs to be subdivided, the discrete second-order partial derivatives along two parametric directions at parametric points are computed and the ratio $\rho$ is prescribed to characterize the anisotropic feature of the solution $u^h$.

The boundary conditions in IGA for modified PHT-splines are imposed as follows.
1)In the case of a homogenous boundary condition on $\Gamma_D$, for each boundary vertex that is mapped into a point in $\Gamma_D$,
find its associated four basis functions. Set the coefficients of those functions that not vanish on $\Gamma_D$ to zero.
2)In the case of a nonhomogeneous boundary condition on $\Gamma_D$, for each boundary vertex that is mapped into a point in $\Gamma_D$, find its associated four basis functions and solve the coefficients of these functions by minimizing errors (the least square method may be an option).
3)In the case of a Neumann boundary condition on $\Gamma_D$, for each boundary vertex that is mapped into a point in $\Gamma_D$, find its adjacent vertices. Construct the equations for the coefficients based on the information of directional derivatives.

\subsection{Numerical experiments}


Here, an IGA problem based on modified PHT-splines is illustrated. For convenience, DOF is used as the abbreviation for degree of freedom. A cell $\theta$ is marked for refinement, if $\eta_\theta >0.0001$. In addition, $\delta =0.5$ is chosen as the threshold to characterize the anisotropic feature.

\begin{figure}[htbp!]
\centering
\includegraphics[width=2.2in]{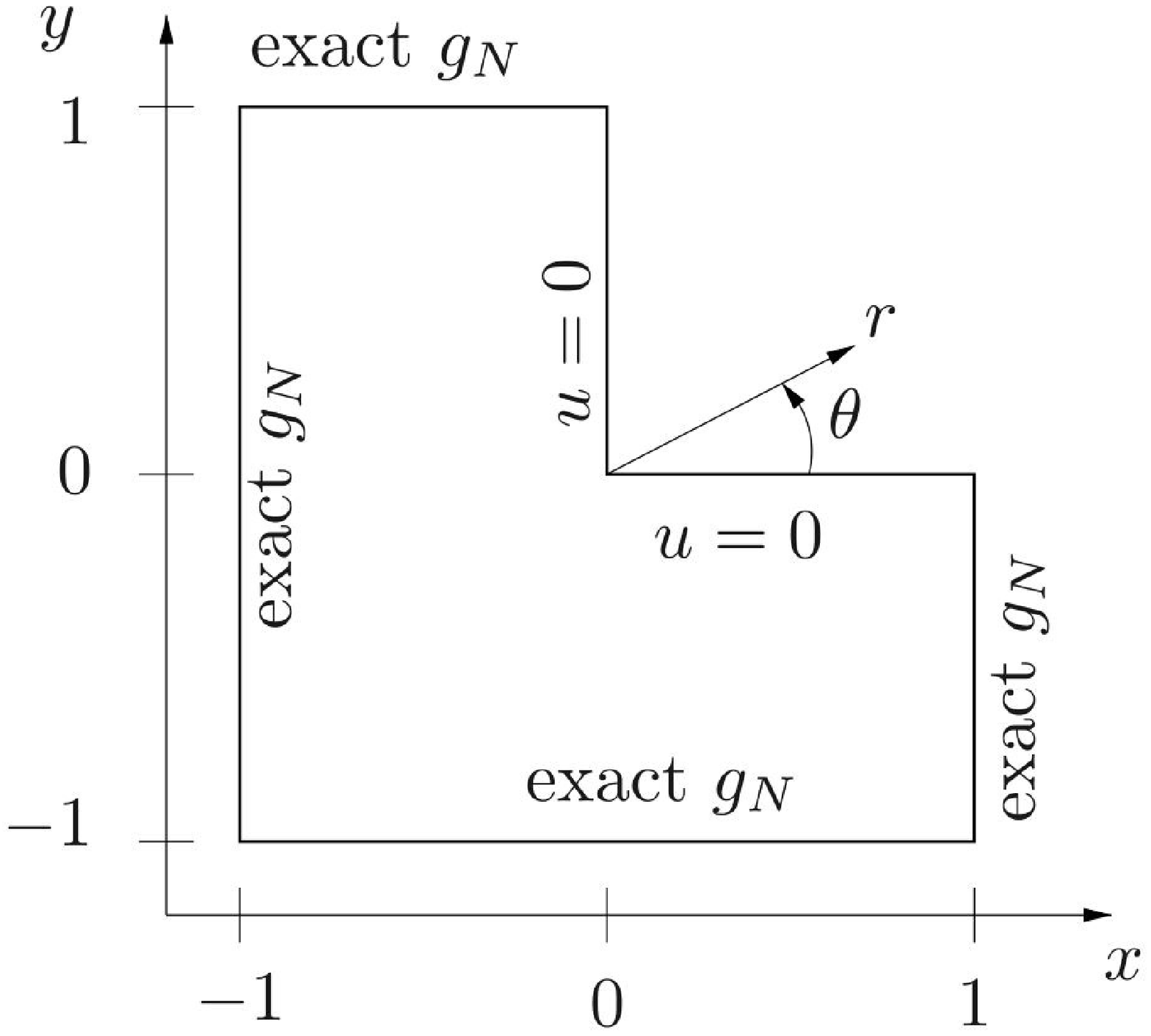}
\caption{Problem setting}
\label{fig:Lshape22}
\subfloat[Initial mesh]{
\includegraphics[width=1.5in]{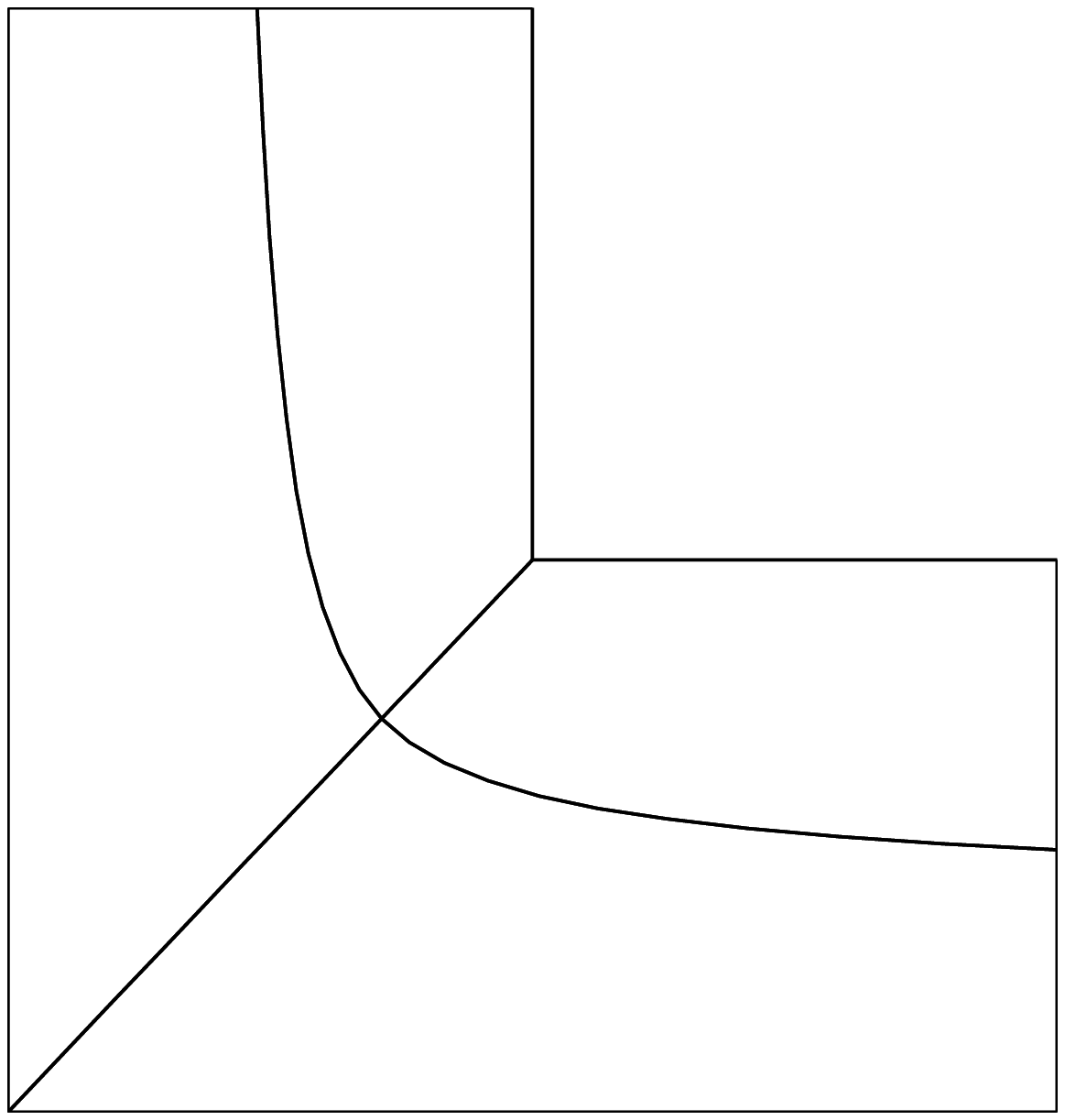}}\hfill
\subfloat[Level $3$]{
\includegraphics[width=1.5in]{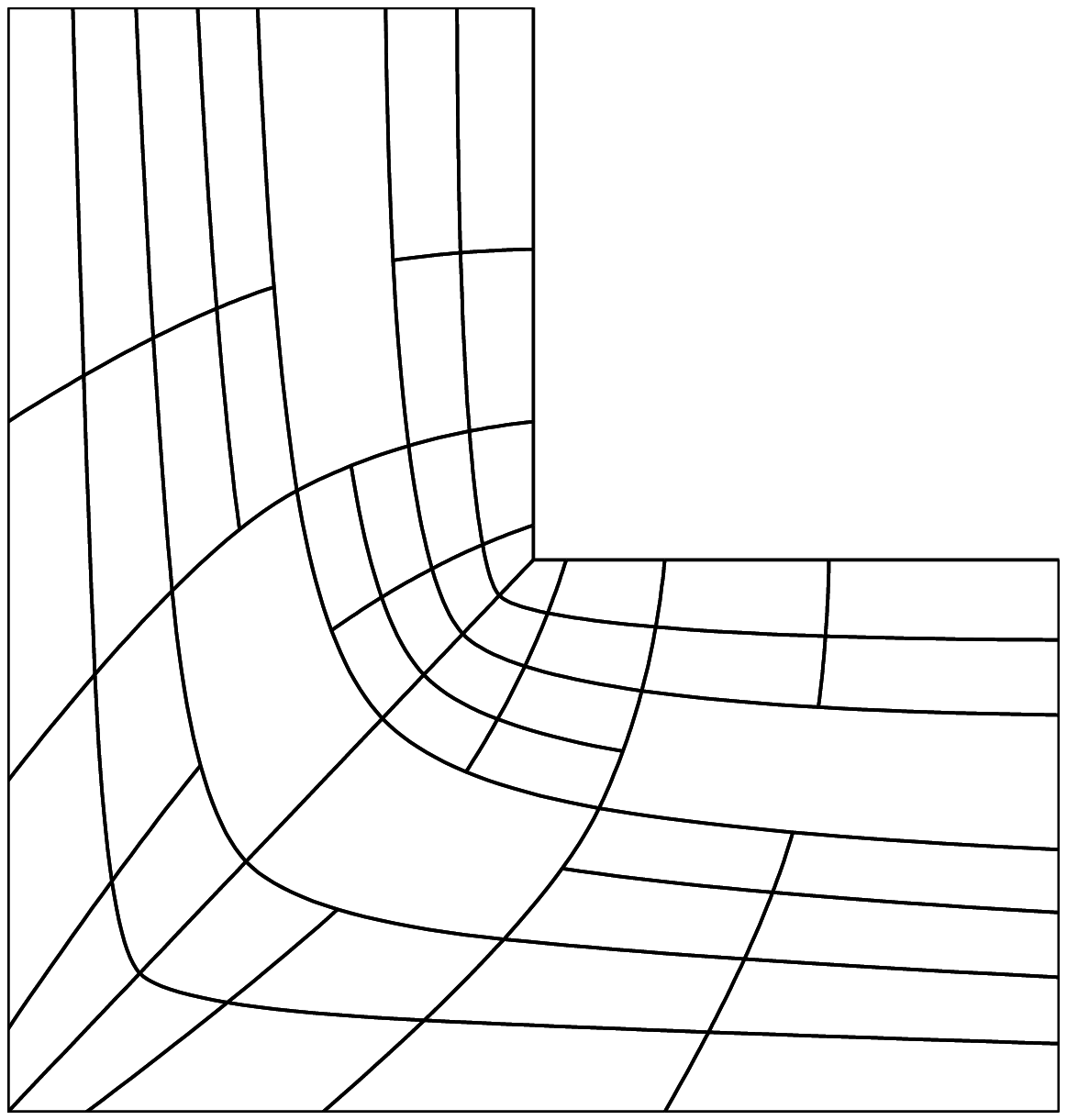}}\hfill
\subfloat[Level $5$]{
\includegraphics[width=1.5in]{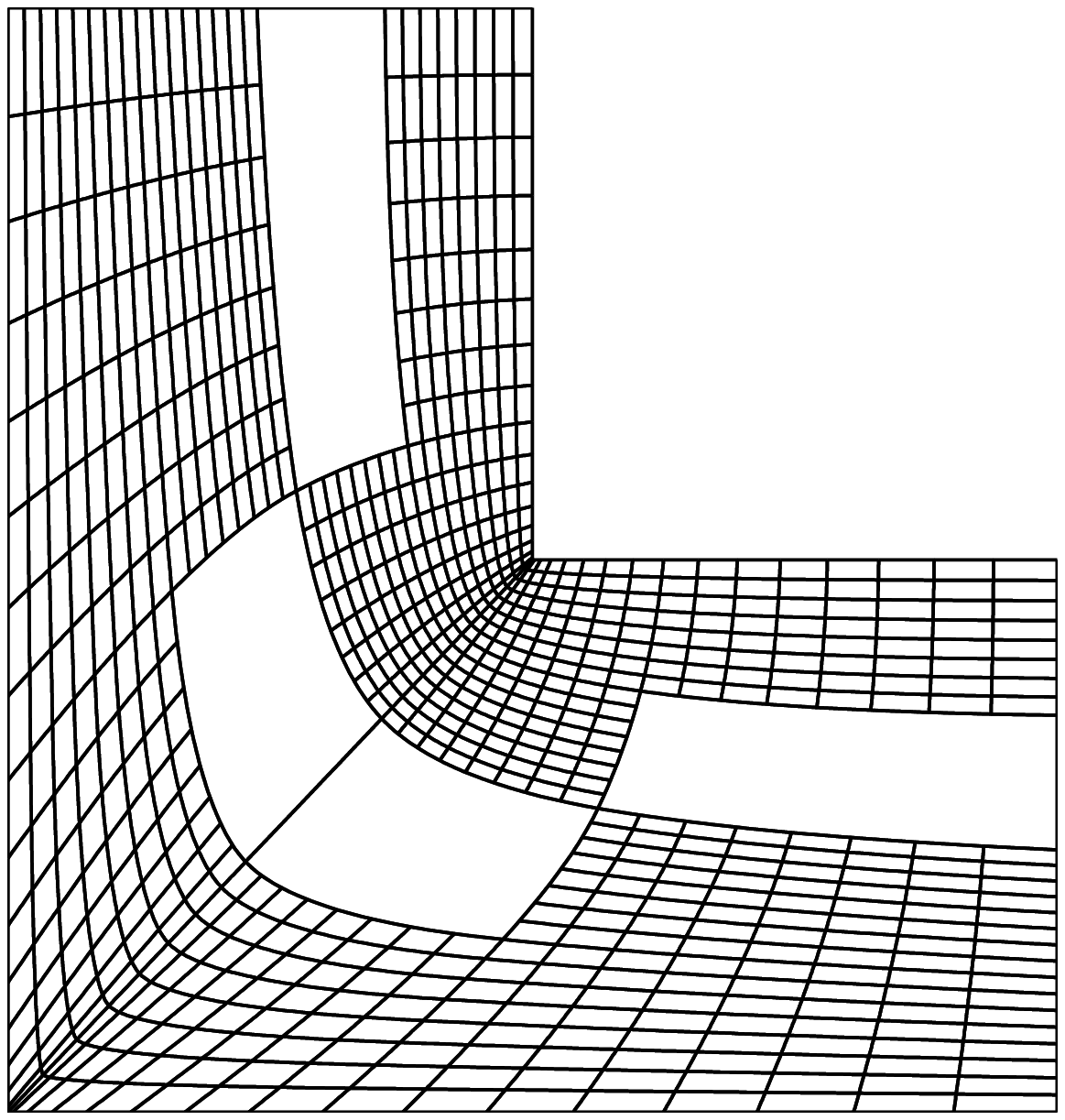}}\hfill
\caption{Meshes on the physical domain by PHT-splines}
\label{fig:Lmeshc}
\subfloat[Initial mesh]{
\includegraphics[width=1.5in]{igampht1.eps}}\hfill
\subfloat[Level $3$]{
\includegraphics[width=1.5in]{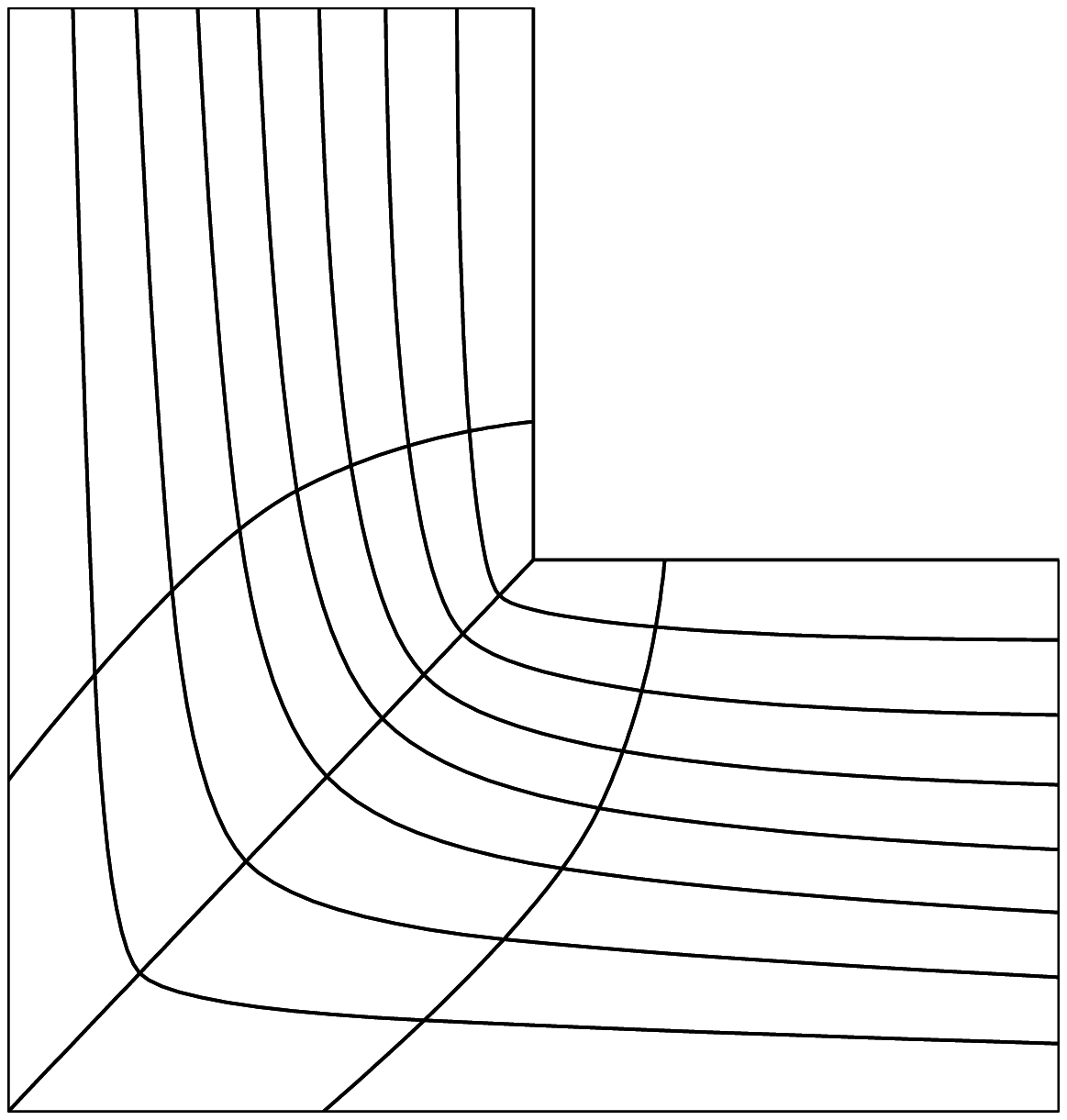}}\hfill
\subfloat[Level $5$]{
\includegraphics[width=1.5in]{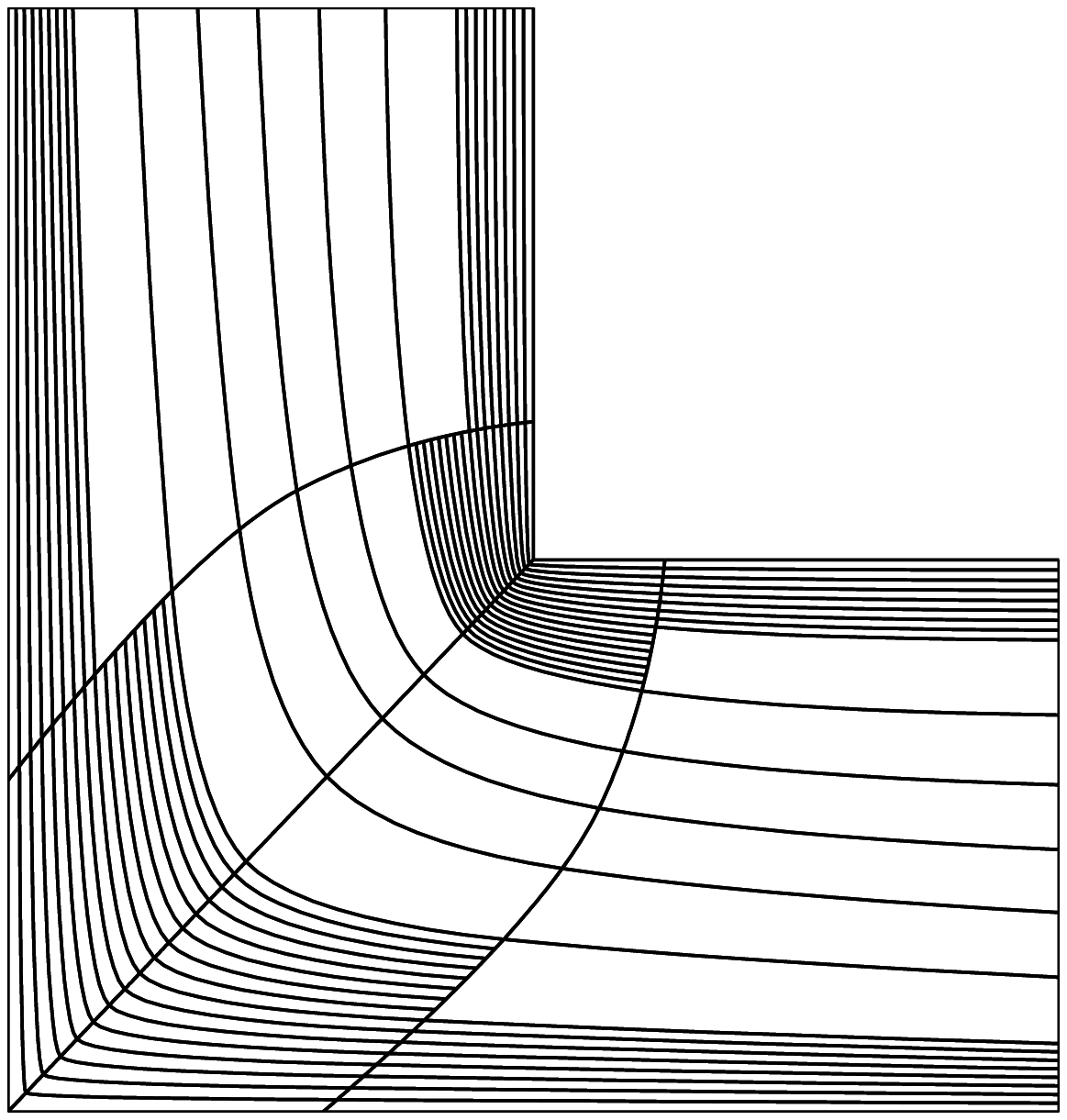}}\hfill
\caption{Meshes on the physical domain of by modified PHT-splines.}
\label{fig:Lmesh}
\end{figure}

The example was also considered by D{\"o}rfel et al \cite{dorfel-2009-iga} and Kleiss et al.  \cite{kleiss-2012-iga}. We solve the Laplace equation \eqref{eq-model} on the $L$-shaped domain
$\Omega = (-1,1)^2\backslash[0,1]^2$, with boundaries $\Gamma_D = \{0\}\times\{0,1\}\cup[0,1]\times\{0\}$, and $\Gamma_N = \partial\Omega\backslash\Gamma_D$, that is illustrated in \reffig{fig:Lshape22}. Double control points are used to model the corners at $(0,0)$ and $(-1,-1)$, the complete geometry data is referred to Appendix A.2 in \cite{kleiss-2012-iga}.

The following function $u$ solves \eqref{eq-model} and is used as our exact solution:
\begin{align*}
  u:\ \ &\mathbb{R}^+ \times (0,2\pi]\rightarrow \mathbb{R}\\
  (r,\theta)&\mapsto r^{\frac{2}{3}}\sin\big( \frac{2\theta-\pi}{3}\big).
\end{align*}
Hence, $g_D$ and $g_N$ in \eqref{eq-model} are determined by the exact solution $u$.

We start from a $3\times3$ tensor-product mesh. The refined meshes on the physical domains solved by PHT-splines and  modified PHT-splines  are shown in \reffig{fig:Lmeshc} and \reffig{fig:Lmesh} respectively.

\begin{figure}[htbp!]
\centering
\includegraphics[width=2.6in]{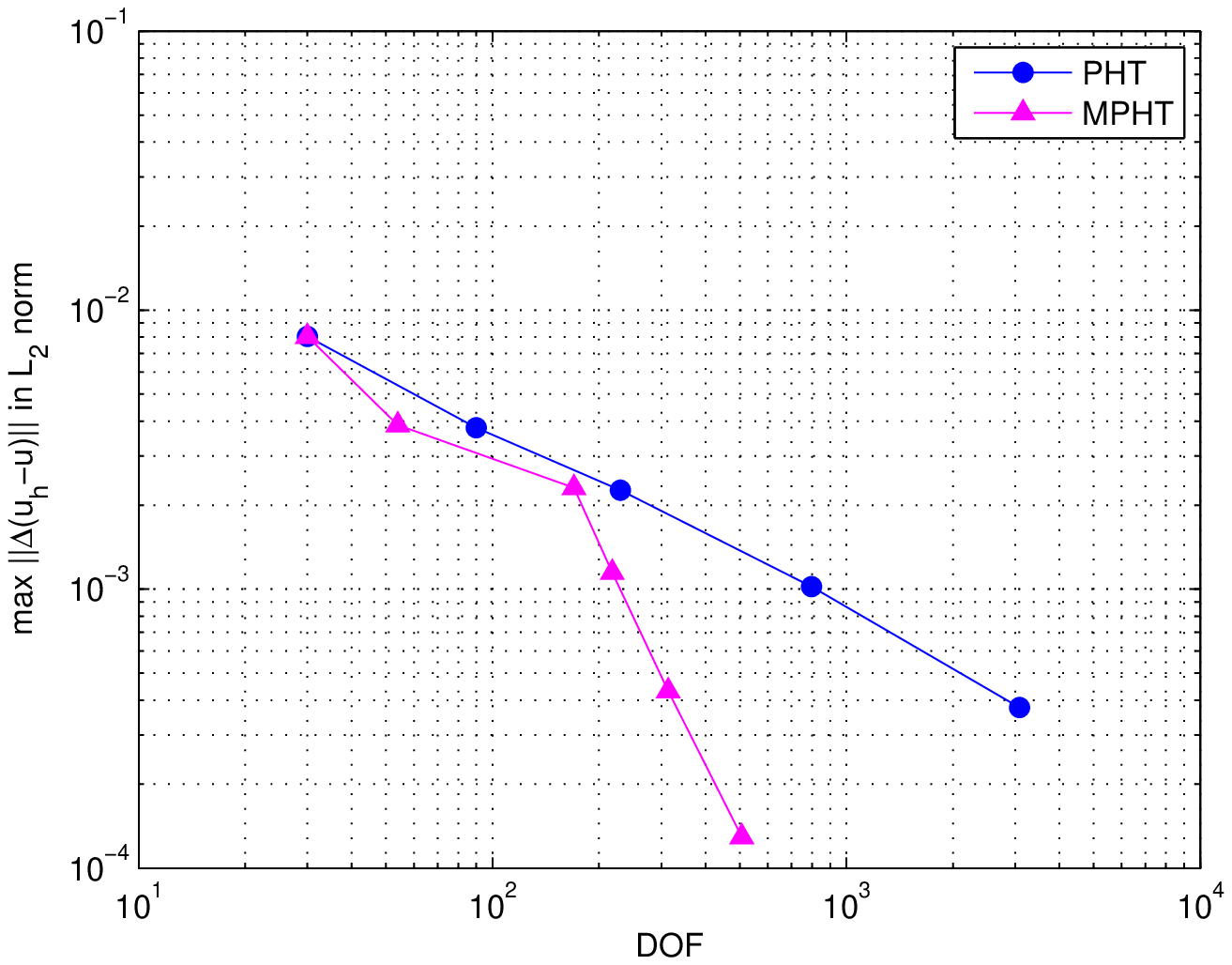}\hfill
\includegraphics[width=2.6in]{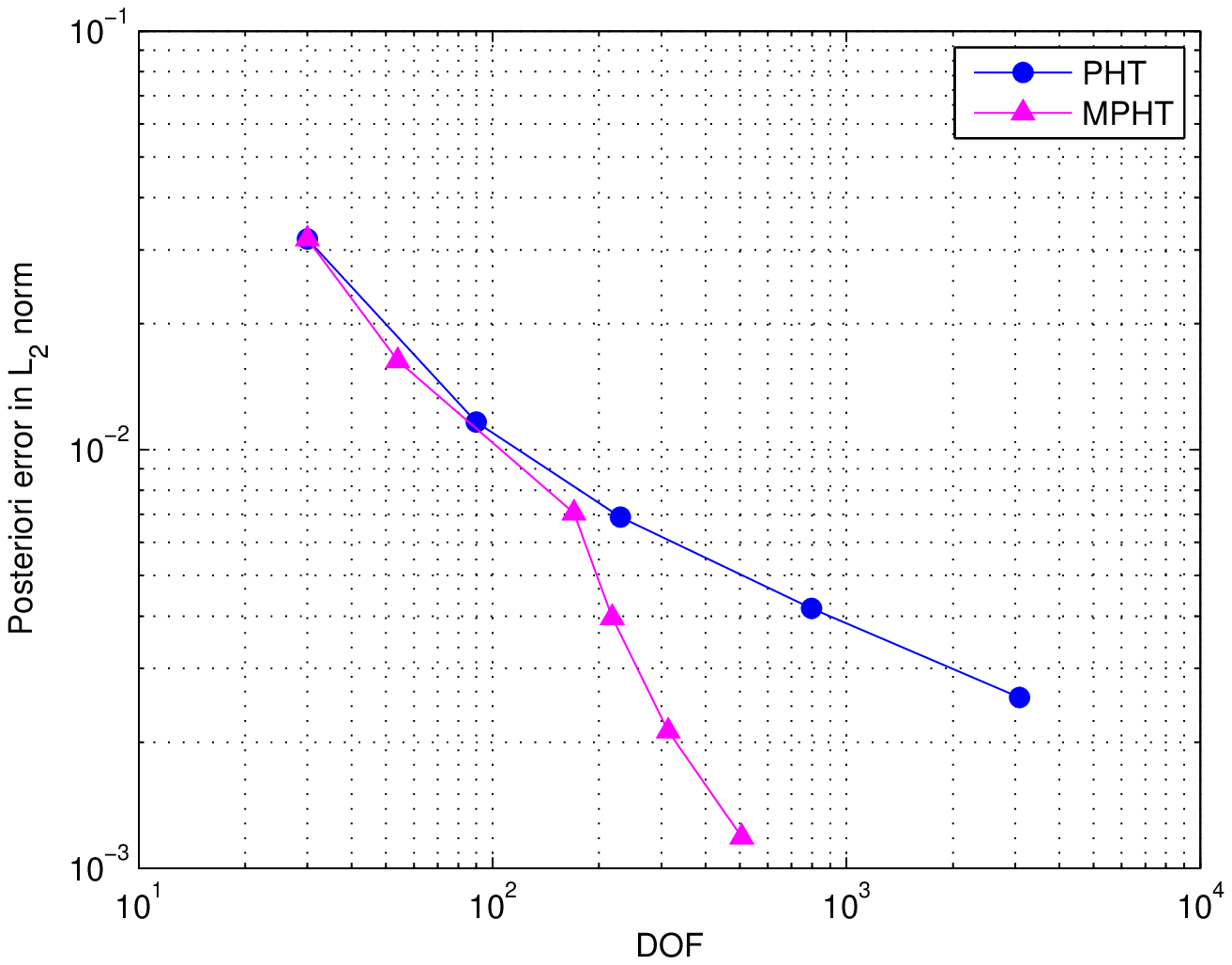}
\caption{Convergence results of PHT-splines and modified PHT-splines.}
\label{fig:iga-compare}
\end{figure}

A comparison of the solution of the PHT-splines and modified PHT-splines is provided in \reffig{fig:iga-compare}. The comparison shows that we can achieve the better accuracy with the same convergence rate.

%
%
%
%
%
%

\section{Conclusion and future work}
\label{section-conclusion}
In this paper, we extend PHT-splines to modified PHT-splines, a special type of polynomial spline defined over modified hierarchical T-meshes. Simultaneously, the basis functions of modified PHT-splines that we constructed inherit the beneficial properties of PHT-splines, i.e., nonnegativity, partition of unity, local support, and linear independence.  A refinement strategy based on neighborhood relations and labels of cells is also presented, in which the labels can be estimated from the geometric information in practical problems.

Modified PHT-splines have been applied to fitting open meshes and isogeometric analysis for the elliptic partial differential equation. Numerical results show that modified PHT-splines have advantages when applied to problems with anisotropic features.

It should be pointed out that only cells of lasted level are marked for refinement in our scheme. 
This assumption may limit our splines to be used in further applications. Especially, for the problems of isogeometric analysis, there is no guarantee that only cells of lasted level need to be marked in some practical problem, although our refinement scheme has good performance in many numerical examples.
In the future, we will focus on improving the algorithm when it is applied to isogeometric analysis.
In addition, we will further explore the applications of modified PHT-splines in isogeometric analysis. To exactly represent common
geometric objects such as circles, cylinders, spheres, and ellipsoids, rational forms of modified PHT-splines may be required.
Furthermore, the generalization of modified PHT-splines in three-dimensional space is worthy of consideration.





\section*{Acknowledgements}
We would like to thank the anonymous referees for providing us with constructive comments and suggestions.

\bibliographystyle{model1-num-names}



\end{document}